\documentclass[opre,nonblindrev]{mystyle} % current default for manuscript submission
\usepackage{mathtools}
\usepackage{lscape}
\usepackage{placeins}
\usepackage{multirow}
\usepackage{setspace}
\usepackage{bm}
\usepackage{fixmath}
\usepackage{pbox}
\usepackage{mathrsfs}
\usepackage{tikz}
\usetikzlibrary{positioning,chains,fit,shapes,calc}
\OneAndAHalfSpacedXI% current default line spacing
\let\footnote=\endnote

\usepackage{tikz} 
\usetikzlibrary{arrows} 
\usetikzlibrary{arrows.meta}
\usetikzlibrary{shapes.geometric}
\usetikzlibrary{fit}
\usetikzlibrary{backgrounds}
\usepackage{natbib}
 \bibpunct[, ]{(}{)}{,}{a}{}{,}%
 \def\bibfont{\small}%
 \def\bibsep{\smallskipamount}%
 \def\bibhang{24pt}%
\usetikzlibrary{automata, positioning}
\usepackage{setspace}
\usepackage{color} 
\usepackage{hyperref}
\usepackage{enumitem}

\usepackage{algorithm}

\usepackage[noend]{algpseudocode}
\TheoremsNumberedThrough     % Preferred (Theorem 1, Lemma 1, Theorem 2)
\ECRepeatTheorems

\EquationsNumberedThrough    % Default: (1), (2),\dots
\newcommand{\mb}{\mathbf}
\newcommand{\mbs}{\boldsymbol}
\newcommand{\mU}{\mathcal{U}} % for the uncertainty set.
\newcommand{\lbar}{\overline}
\newcommand\agg[1]{\widetilde{#1}}

\usepackage{lscape}
\usepackage{booktabs}
\begin{document}
%%%%%%%%%%%%%%%%
\definecolor{myblue}{RGB}{80,80,160}
\definecolor{mygreen}{RGB}{80,160,80}
\definecolor{myred}{RGB}{200,40,40}

%\setcounter{page}{1}
%\MONTH{xxxxx}% (month or a similar seasonal id)
%\DOI{10.1287/xxxx.0000.0000}%

\RUNAUTHOR{Bandi, Han, and Nohadani}

\RUNTITLE{Robust Periodic-Affine Policies}

% \TITLE{Bundling Information Goods of Decreasing Value}
\TITLE{Sustainable Inventory with Robust Periodic-Affine Policies and Application to Medical Supply Chains}

\ARTICLEAUTHORS{%
\AUTHOR{Chaithanya Bandi}
\AFF{Kellogg School of Management, Northwestern University, Evanston, IL 60208\\\EMAIL{c-bandi@kellogg.northwestern.edu}} 
\AUTHOR{Eojin Han and Omid Nohadani}
\AFF{Industrial Engineering \& Management Sciences, Northwestern University, Evanston, IL 60208\\\EMAIL{eojinHan2020@u.northwestern.edu,  nohadani@northwestern.edu}} 
} % end of the block

\ABSTRACT{
\textbf{Abstract}:
We introduce a new class of adaptive policies called \emph{periodic-affine policies}, that allows a decision maker to optimally manage and control large-scale newsvendor networks in the presence of uncertain demand without distributional assumptions. 
These policies are data-driven and model many features of the demand such as correlation, and remain robust to parameter mis-specification. 
We present a model that can be generalized to multi-product settings and extended to multi-period problems.
This is accomplished by modeling the uncertain demand via sets.
In this way, it offers a natural framework to study competing policies such as base-stock, affine, and approximative approaches with respect to their profit, sensitivity to parameters and assumptions, and computational scalability.
We show that the periodic-affine policies are sustainable, i.e. time consistent, because they warrant optimality both within subperiods and over the entire planning horizon.
This approach is tractable and free of distributional assumptions, and hence, suited for real-world applications.
We provide efficient algorithms to obtain the optimal periodic-affine policies and demonstrate their advantages on the sales data from one of India's largest pharmacy retailers.

%
% by %
%
}

\maketitle
%%%%%%%%%%%%%%%%%%%%%%%%%%%%%%%%%%%%%%%%%%%%%%%%%%%%%%%%%%%%%%%%%%%%%%

%%%%%%%%%%%%%%%%%

\section{Introduction}
\label{sec:intro}

Despite the physicians' diagnostic matching of patients to drugs, the heterogeneity in patients' illness, drug's efficacy, potential side effects, and varying length of treatment lead to sizable uncertainty in drug's demand~\citep{crawford2005uncertainty}.
Retailers are mandated to service level guarantees, and overstocking drugs is neither economical nor practical since they are perishable.
Such healthcare problems affect a wide section of the population and have large societal implications.
In this context, newsvendor models offer a natural framework and are used for decision making.

Practical solutions to such problems are critical to a broad range of industries. 
In particular, pharmaceutical companies with a large turnover are interested in optimal inventory management.
GlaxoSmithKline spends over \$4.5 billion each year on manufacturing and supplying products.
Johnson \& Johnson spends approximately \$30 billion annually in leveraging its purchasing power to set sustainability expectations beyond its operations.
Similarly, companies like Teva Pharmaceuticals, Pfizer, and Merck spend millions of dollars to ensure the safety and supply of their products, even though they have manufacturing units in multiple locations.
Therefore, any variation in inventories can lead to multiple disturbances in the system.
A pharmacy's inventory represents its single, largest investment.
As discussed in \cite{pmref2}, in a common pharmacy, cost of goods sold accounts for approximately 68\% of total expenditures.
For every 1\% change in costs of goods, profits may increase or decrease by more than 20\% (see \cite{pmref1}).
Thus, the sheer magnitude of dollars involved makes seemingly minor inefficiencies in purchasing and inventory control matter of great importance to both cash flow and profitability. 

The challenges of such networks are multifold.
Real-world settings are typically high-dimensional with multiple products and multiple stages of decision-making.
These settings also suffer from substantial uncertainties in demand.
Modeling such demand uncertainty is challenging because demand is often not stationary or its uncertainty can depend on previous decisions~\citep{nohadani2018}.

In this work, we consider a newsvendor network with uncertain and correlated demand.
Using the paradigm of robust optimization, we model such demand to reside in uncertainty sets and provide tractable formulations and associated algorithms for sustainable policies.
To gain insight from a real-world setting, we apply the results to a major online pharmacy retailer in India, where a prohibitively large penalty occurs when customers' demand is not satisfied.  
This company carries over 163 different brands, and the sales grow at about 23\% per year.
Their distribution network spans the entire country through fixed retail locations and online platforms.
The decision makers of this company observe a sizable uncertainty in demand over the course of the year (in addition to seasonality) and significant correlations amongst various product categories. 
In close collaboration with this company's managers, we seek to design optimal implementable policies to control their inventory levels in their network.

Our contributions are:
\begin{itemize}
\item  \emph{Modeling}: We provide a distribution-free description of uncertainty in demand using two types of sets.
Independent demands are modeled via budget constraints.
We also incorporate correlated demands using a factor model approach.
The inventory control problem is then cast as a multi-stage robust optimization problem.
As a result, a novel solution concept of \emph{periodic-affine policy} is provided for newsvendor networks with time-dependent and potentially correlated demand uncertainty.
\item \emph{Algorithms}: We provide a tractable algorithm that provides periodic-affine policies.
These policies decompose the overall problem into a more tractable formulation than affine policies.
\item \emph{Application}: We analyze the sales data of a pharmacy retailer in India for the fourth quarter of 2016.
This entails 1.5 million transactions for 228 different products.
We construct the demand uncertainty set for the 20 most-popular products, comprising 80\% of all transactions.
Our numerical experiments show that even for the single-station case, the computational burden for the optimal periodic-affine policies is significantly reduced over affine policies (by 100$\times$ for a 15-period problem), making the proposed approach practical for real-world and large-sized problems.
Moreover, the periodic-affine policy improved the cost effectiveness of the operation by 19\% over a base-stock policy for realistic penalty costs.
\end{itemize}

\subsection{Literature review}
The seminal work of \cite{Arrow1951} introduced the multistage periodic review inventory model, where the inventory is reviewed once every period and a decision is made to place an order, if a replenishment is necessary.
The $(s,S)$ inventory policy establishes a lower (minimum) stock point $s$ and an upper (maximum) stock point $S$.
When the inventory level drops below $s$, an order is placed ``up to $S$."
The $(s,S)$ ordering policy has been proven optimal for simple stochastic inventory systems.
\cite{Scarf1960} proved that base-stock policies are optimal for a single installation model.
\cite{Clark1960} extended the result to serial supply chains without capacity constraints and showed that the optimal ordering policy for the multi-echelon system can be decomposed into decisions based on the echelon inventories.
\cite{Karlin1960} and \cite{Morton1978} showed that base-stock policies are optimal for single-state systems with non-stationary demands.
\cite{Federgruen1986} generalized the analysis to a single-stage capacitated system, and \cite{Rosling1989} extended the analysis of serial systems to assembly systems.
For more work, refer to \cite{Langenhoff1990, Sethi1997, Muharremoglu2008, Huh2008}. 

Simulation optimization has attempted to take advantage of the availability of computational resources and the power of simulation for evaluating functions. 
For a comprehensive overview of commonly used simulation optimization techniques, we refer the reader to the survey by \cite{Fu2005}.
\cite{Fu1994b, Glasserman1995, Fu1997b} and \cite{Kapiscinski1999} have developed various gradient-based algorithms to study inventory systems. 
These methods are practical whenever the input variables are continuous and their success depends on the quality of the gradient estimator.

On the other hand, \textcolor{black}{\cite{Scarf1958, Kasugai1961, Gallego1993, Graves2000}} developed distribution-free approaches to inventory theory. 
\cite{Bertsimas2006} took a robust optimization approach to inventory theory and showed that base-stock policies are optimal in the case of serial supply chain networks. 
\cite{Bienstock2008} presented a family of decomposition algorithms aimed at solving for the optimal base-stock policies using a robust optimization approach. 
\cite{Rikun2011} extended the robust framework introduced by \cite{Bienstock2008} to compute optimal $(s,S)$ policies in supply chain networks and compared their performance to optimal policies obtained via stochastic optimization. 
\cite{ben2004} extended the robust optimization framework to dynamic settings and explored the use of disturbance-affine policies by allowing the decision maker to adjust their strategy leveraging the information revealed over time.  
\cite{Bertsimas2006} and \cite{Bienstock2008} studied the performance of base-stock policies, and \cite{ben2005retailer}, \cite{Kuhn2011}, and \cite{Bertsimas2010} investigated polices that are affine in prior demands under a robust optimization lens.
Within the robust optimization framework, affine policies have gained much attention due to their tractability; depending on the class of the nominal problem, the optimal policy can be solved via linear, quadratic, conic or semidefinite programs (see \cite{Lofberg2003, Kerrigan2004}). 
Empirically, \cite{ben2005retailer} and \cite{Kuhn2011} have reported that affine policies have excellent performance and in many instances optimal.

Another approach is distributionally robust optimization which assumes that the uncertainties follow a distribution within a prespecified set of distributions.
Such sets can be based on moment constraints \citep{delage2010distributionally}, phi-divergences \citep{ben2013robust}, or Wasserstein metric \citep{esfahani2015data} to allow tractable reformulations.
This approach typically yields less conservative solutions than deterministic robust optimization solutions.
For multi-stage problems, \cite{van2016distributionally} proposed a tractable framework for distributionally robust linear feedback policy for discrete time linear control systems with quadratic objective functions.

%%%%%%%%Pharma lit review%
In the context of pharmaceutical systems, \cite{guerrero2013joint} provided a near-optimal base-stock policy for two-echelon distribution networks with multiple products, where every sink node is replenished by a single supplier. 
They provided a Markov chain formulation and a heuristic algorithm for Poisson distributed and independent demands.
For a combined setting of a pharmaceutical compony and a hospital, \cite{uthayakumar2013pharmaceutical} developed a two-echelon supply chain model to determine the optimal lot size, lead time, and total number of deliveries between the pharmaceutical compony and a hospital.
Using Lagrange multipliers, they provided decision tools for optimal costs while ensuring required service levels.
In a two-level pharmaceutical supply chain, \cite{baboli2011replenishment} studied a specific product with a constant demand rate and numerically showed that the overall cost is improved when pharmacies and hospitals are centralized.% and co-ordinated.%, rather than a decentralized system.

%%%%%%%%%%%%%%%%5

\subsubsection*{Notation.}
Lowercase italic is used to denote scalars; lowercase bold is used to denote vectors, and uppercase bold is used to denote matrices. 
Sets are in calligraphic.
Section specific notation is introduced where needed.
All proofs are relegated to the Online Appendix.

\section{Model}
\label{sec:model}
We consider \textcolor{black}{a newsvendor network} in which inventories are reviewed periodically and unfulfilled orders are backlogged. 
For simplicity, we assume zero lead times throughout the network; however, our framework can be adapted to systems with non-zero lead times. 
We consider a $T$-period time horizon and, within each period, events occur in the following order: (1) the ordering decision is made at the beginning of the period, (2) demands for the period occur and are filled or backlogged depending on the available inventory, (3) the stock availability is updated for the next period.

\begin{itemize}
		\item $\mathcal{N}$ 	:	 Set of all installations where ordering decisions are made (source nodes) with $|\mathcal{N}| = m$
				\item $\mathcal{S}$        :	 Set of all installations with external demand (sink nodes) with $|\mathcal{S}| = n$
		\item$\mathcal{L}$ 		: Set of all links (edges) within the inventory network with $|\mathcal{L}| = p$
		\item $\mathcal{N}_{{k}}$ 	: Set of source nodes supplying stock to a sink node ${k} \in \mathcal{S}$
	\item	$\mathcal{S}_{v}$	: Set of sink installations that are fed from a source node $v \in \mathcal{N}$
\item $s_{t}^{v}$	:	Amount of order at the beginning of period $t$ at a source $v \in \mathcal{N}$
	\item	$d_{t}^{{k}}$	: Demand observed at a sink ${k}\in \mathcal{S}$ throughout time period 
	\item	$x_{t}^{\ell}$	:	 Stock delivered along a link $\ell\in\mathcal{L}$ at time $t$
	\item	$u_{t}^{s,v}$ 	:	 Stock available after the period $t$ at a source node $v \in \mathcal{N}$
	 \item	$u_{t}^{d,{k}}$		:	 Backorders after the period $t$ at a sink node ${k} \in \mathcal{S}$. 
		\end{itemize}

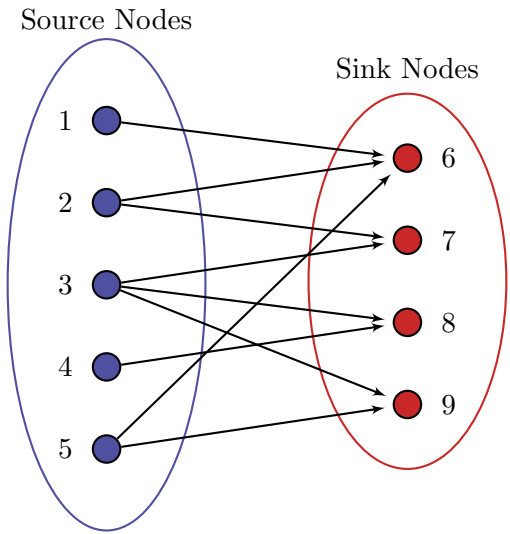
\begin{figure}[h!]
	\centering
\vspace{-1.2cm}
\begin{tikzpicture}[thick, %yscale=2,
  every node/.style={draw,circle},
  fsnode/.style={fill=myblue},
  ssnode/.style={fill=myred},
  every fit/.style={ellipse,draw,inner sep=-2pt,text width=2cm},
  >={latex'[length=6mm, width=6mm]},shorten >= 3pt,shorten <= 0pt
]

% the vertices of U
\begin{scope}[start chain=going below,node distance=7mm]
\foreach \i in {1,2,...,5}
  \node[fsnode,on chain] (f\i) [label=left: \i] {};
\end{scope}

% the vertices of V
\begin{scope}[xshift=4cm,yshift=-0.5cm,start chain=going below,node distance=7mm]
\foreach \i in {6,7,...,9}
  \node[ssnode,on chain] (s\i) [label=right: \i] {};
\end{scope}

% the set U
\node [myblue,fit=(f1) (f5),label={[yshift=-1.05cm]$\text{Source Nodes}$}] {};
% the set V
\node [myred,fit=(s6) (s9),label={[yshift=-0.8cm]$\text{Sink Nodes}$}] {};

% the edges
\draw[->] (f1) -- (s6);
\draw[->] (f2) -- (s6);
\draw[->] (f2) -- (s7);
\draw[->] (f3) -- (s7);
\draw[->] (f3) -- (s8);
\draw[->] (f3) -- (s9);
\draw[->] (f4) -- (s8);
\draw[->] (f5) -- (s9);
\draw[->] (f5) -- (s6);
\end{tikzpicture}

	\label{networks}
	\caption{\small{Example of a nine-installation network with $n=4$ sink nodes and $m=5$ source nodes.}}
\end{figure}

To track the system's operation, we capture information about the stock available and the stock ordered at source installations at the beginning of each time period as well as the demand at each sink installation throughout each time period. 
Specifically, assuming zero initial input and demands, we can express the dynamics of inventory levels and backlogged demands for $t=1, \ldots, T$ as 
\begin{equation}
\label{echelon-inventories}
\begin{aligned}
&u_{t}^{s,v} = u_{t-1}^{s,v} + s_{t}^{v} - \!\!\!\!\!\!\sum_{ \ell=(v,{k}), {k} \in \mathcal{S}_{v}} \!\!\!\!\!\!x_{t}^{\ell} = \sum_{\tau=1}^{t} s_{\tau}^{v} - \!\!\!\!\!\!\sum_{\ell=(v,{k}), {k} \in \mathcal{S}_{v}} \sum_{\tau=1}^{t} x_{\tau}^{\ell}~~~~~~\forall~ v \in \mathcal{N} \\
&u_{t}^{d,{k}} = u_{t-1}^{d,{k}} + d_{t}^{{k}} -\!\!\!\!\!\! \sum_{ \ell=(v,{k}), v \in \mathcal{N}_{{k}}} \!\!\!\!\!\!x^{\ell}_{t} = \sum_{\tau=1}^{t} d_{\tau}^{{k}} - \!\!\!\!\!\!\sum_{\ell=(v,{k}), v \in \mathcal{N}_{{k}}} \sum_{\tau=1}^{t} x_{\tau}^{\ell}~~~~\forall~ {k} \in \mathcal{S}. \\
\end{aligned}
\end{equation}
Note that the ordering quantities ${s_{t}^{v} = s_{t}^{v}(\pi, \mb{d})}$, and therefore the amount of available stock ${u_{t}^{s,v} = u_{t}^{s,v}(\pi, \mb{d})}$ and backorders ${u_t^{d,{k}} = u_t^{d,{k}}(\pi, \mb{d})}$, are functions of the ordering policy $\pi$ and the demand $\mb{d}$.

The high-dimensional nature of modeling demand uncertainty probabilistically and the complex dependence on random variables underscore the difficulty of analyzing and optimizing the expected total cost.
Instead, we propose a framework that builds upon the robust optimization paradigm.

%%%%%%%%%%%%%%%%%
\subsection{Robust Newsvendor Network Formulation}
\label{sec:prob_form}

To describe our framework, we first introduce a robust approach to single-period models.
Our models are based on assumption that we have the following cost and revenue structure:
\begin{itemize}
\item $c_S^{v}$:	 Purchasing cost per unit at the source node ${v} \in \mathcal{N}$ 
	\item	$c_H^{v}$:  	 Holding cost per unit for the leftover stock at the source node ${v} \in \mathcal{N}$ 
	\item	$c_P^{k}$: Penalty cost per unit for the unsatisfied demand at the sink node ${k} \in \mathcal{S}$ 
	\item	$r^\ell$: Revenue by satisfying a unit demand occurred at the sink node ${k}$ via $\ell = ({v},{k}) \in \mathcal{L}$. 
	\end{itemize}
The goal of the decision maker is to order a proper amount of products $\{s_{v} : {v} \in \mathcal{N}\}$ and to process network activities $\{x_{\ell} : \ell \in \mathcal{L} \}$ to satisfy the customer demand at the sink nodes, so that the firm maximizes an overall profit.
If we denote $\mU$ as a demand uncertainty set, then a single-period problem is formulated as a two-stage robust optimization problem
\begin{equation}
\begin{aligned}
&\max_{s_{v}\geq 0}~ \Bigg[ -\sum_{{v} \in \mathcal{N}} c_S^{v} s_{v} +~ 
\min_{\mb{d} \in \mU}~\max_{x_\ell\geq 0} \bigg[
\sum_{\ell \in \mathcal{L}} r_{\ell} x_{\ell} - \sum_{{k} \in \mathcal{S}} c_P^{k} \Big( d_{k} -~ \sum_{\mathclap{\ell=({v},{k}), {v}\in \mathcal{N}_{k}}} x_{\ell} ~\Big) - \sum_{{v} \in \mathcal{N}} c_H^{v} \Big( s_{v} -~ \sum_{\mathclap{\ell=({v},{k}), {k} \in \mathcal{S}_{v}}} x_{\ell} ~\Big) \bigg] \Bigg] \\
&\text{s.t. }~~~~ \sum_{\ell=({v},{k}), {k} \in \mathcal{S}_{v}} \!\!\!\!\!\!x_{\ell} \leq s_{v} ~~\forall~ {v} \in \mathcal{N}, \sum_{\ell = ({v},{k}), {v} \in \mathcal{N}_{k}} \!\!\!\!\!\!x_{\ell} \leq d_{k} ~~\forall~ {k} \in \mathcal{S} ,
\end{aligned}
\label{eq:nominal_long}
\end{equation}
where the constraints are network constraints and affect the inner maximization problem.
Note that the order quantities ${\{s_{v}:{v}\in\mathcal{N}\}}$ are ``here-and-now'' decisions; it must be placed before demands are realized, while the network activities $\{x_\ell: \ell \in \mathcal{L}\}$ are ``wait-and-see'' solutions and assigned after demands are observed.

\subsubsection*{Notation.}
To simplify (\ref{eq:nominal_long}), we define $\mb{c}_S \in \mathbb{R}_+^m$, $\mb{c}_H \in \mathbb{R}_+^m$, $\mb{c}_P \in \mathbb{R}_+^n$ and $\mb{r} \in \mathbb{R}_+^p$ as cost and revenue vectors, and define $\mb{R}_S \in \mathbb{R}_+^{m \times p}$ and $\mb{R}_D \in \mathbb{R}_+^{n \times p}$ as matrices that describe the two constraints, respectively.
Decision variables and uncertain demands are $\mb{s} \in \mathbb{R}_+^m$, $\mb{x} \in \mathbb{R}_+^p$, and $\mb{d} \in \mathbb{R}_+^n$.
We obtain
%Then (\ref{eq:nominal_long}) is simplified as 
\begin{equation}
\begin{aligned} 
	&\underset{\mb{s} \geq \mb{0}}{\text{max}} ~\Bigg[ -(\mb{c}_{S}+\mb{c}_{H})^{\top} \mb{s} ~+~ 
	\underset{\mb{d} \in \mU}{\text{min}}~\underset{\mb{x} \ge\mb{0}}{\text{max}} 
	\big[ \mb{v^\top}\mb{x}-\mb{c}_{P}^{\top} \mb{d} \big] \Bigg] \nonumber \\
	&\text{s.t. }~~~~ \mb{R}_S \mb{x} \leq \mb{s}, ~~~~\mb{R}_D \mb{x} \leq \mb{d}, \label{eq:single_obj} 
\end{aligned}
\end{equation}
with $\mb{v} = \mb{r} + \mb{R}^{\top}_S \mb{c}_H + \mb{R}^{\top}_D \mb{c}_P$.

\subsection{Modeling Demand Uncertainty} 

For the sake of simplicity, we assume that there is no demand seasonality and that the demand realizations are light-tailed in nature (i.e., the demand variance is finite). 
For each sink installation $k \in \mathcal{S}$, we denote the demand mean by $\mu_{k}$ and the demand standard deviation by $\sigma_{k}$.
Our framework also captures correlation among the demand, where we denote $\mbs{\Sigma} \in \mathbb{R}^{n\times n}$ as the nominal covariance matrix.
Note that all these values can be inferred from historical data.
Instead of describing the demand as a random variable, we describe the demand and its correlation by using budget uncertainty sets \citep{bertsimas2004price} and a factor-based approach \citep{bandi2012tractable}.
Such sets do not require any distributional assumption other than first two moments, and consequently, they are robust to the distribution choice.

We capture the correlations via the covariance matrix $\mbs{\Sigma}$ with rank $l \leq n$.
This means, there exist $\mb{A}$ and $ \lambda_1, \ldots, \lambda_l > 0 $ that satisfy ${\mbs{\Sigma} ~=~ \mb{A} \cdot \textbf{diag}(\lambda_1^2,\ldots,\lambda_l^2) \cdot \mb{A^{\top}}}$.
%Using these, we define the uncertainty set as follows.

\begin{definition}[Single-period Uncertainty set]
	\label{uncertainty-set-demand}
	The uncertainty set for correlated demands at sink nodes $\mb{d} = (d_1,\ldots, d_n)$ with variability parameters $\Gamma,~\Gamma_B \geq 0$ is  
	\begin{align}
	\mU = \left\{ \mb{d} \in \mathbb{R}_+^n ~\Bigg|~  \mb{d} = \mbs{\mu} +  \mb{A} \cdot \mb{z},~
	\sum_{i=1}^{l} \left| \frac{z_i}{\lambda_i} \right| \leq \Gamma,~ \left| \frac{z_i}{\lambda_i} \right| \leq \Gamma_{B} ~~\forall~ i=1,\ldots,l \right\}.
	\label{def:unc_set_corr}
	\end{align}
\end{definition}
Note that in this definition, $\Gamma$ and $\Gamma_B$ control the degree of conservatism.
The first constraint in $\mU$ captures correlation, and the others are budget constraints which limit the absolute deviation from its nominal value.
While $\mU$ is data driven, it also captures previous results on the effect of mean and standard deviation on the profit in newsvendor networks.
In particular, $\mU$ recovers the insightful properties in~\cite{mieghem2002newsvendor}, as proposed in the following.

\begin{proposition}
	For a single-period robust newsvendor network with the uncertainty set~$\mU$, the worst-case profit increases in $\mu_i$ and decreases in $\lambda_i$.
	\label{prop:van_mieghem}
\end{proposition}

This proposition shows that our framework generalizes the structural properties from stochastic networks without distributional assumptions.
We extend our model to multi-period cases in the next section.

\section{Multi-period Robust Newsvendor Networks}
\label{subsec:multiperiod_formulation}

To extend the single-period models into dynamic cases, we consider a decision maker who has multiple processing points of $T$ periods.
We assume that all parameters $\mb{c}_S, \mb{c}_H, \mb{c}_P, \mb{r}$, and matrices $\mb{R}_S, \mb{R}_D$ remain constant over the time horizon.
As in Section \ref{sec:model}, on-hand input stocks at source nodes and unsatisfied demand at sink nodes are backlogged to the next periods.
We also assume that the demands are correlated over sink nodes, but independent over time, with nominal mean vector $\mbs{\mu}_t$ and covariance matrix $\mbs{\Sigma}_t$ for each time period $t=1,\ldots,T$.

\subsubsection*{Notation.}
Order quantities at time $t$ are denoted by $\mb{s}_t$, customer demands by $\mb{d}_t$, and network activities by $\mb{x}_t$.
Single-station quantities are denoted by $s_t, d_t, x_t$.
Aggregated amount of orders up to time $t$ are denoted by $\agg{\mb{s}}_t$, customer demands by $\agg{\mb{d}}_t$, and network activities by $\agg{\mb{x}}_t$ ($\agg{s}_t, \agg{d}_t, \agg{x}_t$ for single-station).
Inventory levels and backlogged demands after time $t$ are denoted by $\mb{u}_t^\mb{s}$ and $\mb{u}_t^\mb{d}$.
Finally, $\mb{D}_{[t_1:t_2]} = (\mb{d}_{t_1},\ldots,\mb{d}_{t_2}) \in \mathbb{R}_+^{n \times (t_2 - t_1+1)}$ contains every realized demand from time $t_1$ to $t_2$.
Other quantities such as $\mb{S}_{[t_1:t_2]}$ and $\mb{X}_{[t_1:t_2]}$ are defined similarly.
We define $\mb{A}_t$ and $\lambda_{t,1},\ldots,\lambda_{t,l_t}$ for each $t$, with $rank(\mbs{\Sigma}_t) = l_t$ and $\mbs{\Sigma}_t = \mb{A}_t \cdot \textbf{ diag}(\lambda_{t,1}^2,\ldots,\lambda_{t,l_t}^2) \cdot \mb{A}_t^\top$.

In the following, we generalize Definition \ref{uncertainty-set-demand} for multi-period demand.

\begin{definition}[Multi-period Uncertainty Set]
	The uncertainty set for the demand at sink nodes $(\mb{d}_1,\ldots, \mb{d}_T) \in \mathbb{R}^{n \times T}$ over $T$ periods is 
	\begin{align*}
	\mU^T = \Bigg\{ &\left( \mb{d}_1,\ldots,\mb{d}_T \right) ~\bigg|~ \mb{d}_t = \mbs{\mu}_t + \mb{A}_t \mb{z}_t ~~\forall t=1,\ldots,T \\
	&\sum_{t=1}^T \sum_{i=1}^{l_t} \left| \frac{z_{t,i}}{\lambda_{t,i}} \right| \leq \Gamma,~ \sum_{i=1}^{l_t} \left| \frac{z_{t,i}}{\lambda_{t,i}} \right| \leq \Gamma_t,~ \left| \frac{z_{t,i}}{\lambda_{t,i}} \right| \leq \Gamma_B  ~~\forall i=1,\ldots,l_t,\; t=1,\ldots,T \Bigg\}.
	\end{align*}
\label{assumption:multi_period}
\end{definition}
In this set, the additional constraint controls the absolute deviation over nodes and time periods.
It prevents the demand to take extreme values in every period $t$, which reduces the conservatism over time.
This definition can also describe seasonality of demands, which applies to many areas. 
When there is an explicit time-dependence between the periods, $\mU^T$ can be expressed as a conic set~\citep{nohadani2017robust}.

For the multi-period newsvendor networks, we can express the dynamics of inventories and backlogged demands in (\ref{echelon-inventories}) with vectors and matrices as
\begin{align*}
& \mb{u}_{t}^{\mb{s}} = \mb{u}_{t-1}^{\mb{s}} + \mb{s}_{t} - \mb{R}_S \mb{x}_{t}
= \sum_{\tau=1}^{t} \Big( \mb{s}_\tau - \mb{R}_S \mb{x}_\tau \Big) \\
& \mb{u}_{t}^\mb{d} = \mb{u}_{t-1}^\mb{d} + \mb{d}_{t} - \mb{R}_D \mb{x}_{t}
= \sum_{\tau=1}^{t} \Big( \mb{d}_\tau - \mb{R}_D \mb{x}_\tau \Big),
\end{align*}
and model a multi-stage robust optimization problem as
\begin{equation}
	\begin{split}
	\max_{\agg{\mb{s}}_{t}(\mb{D}_{[1:t-1]})}
	\min_{\mb{D}_{[1:T]} \in \mU^T}
	\max_{\agg{\mb{x}}_t \in \mathcal{P}(\agg{\mb{s}}_t, \agg{\mb{d}}_t, \agg{\mb{x}}_{t-1} )} ~ \Bigg[
	- \mb{c}^\top_S \agg{\mb{s}}_T(\mb{D}_{[1:T-1]}) 
	 - \mb{c}^\top_H \sum_{t=1}^{T} \bigg[ \agg{\mb{s}}_t(\mb{D}_{[1:t-1]}) - \mb{R}_S \agg{\mb{x}}_t(\mb{D}_{[1:t]}) \bigg] \\
	- \mb{c}^\top_P \sum_{t=1}^{T} \bigg[ \agg{\mb{d}}_t - \mb{R}_D \agg{\mb{x}}_t(\mb{D}_{[1:t]}) \bigg] + \mb{r}^{\top} \agg{\mb{x}}_T(\mb{D}_{[1:T]}) \Bigg],
	\end{split}
\label{eq:multi_obj}
\end{equation}
where $\mb{D}_{[1:0]} = \mb{0}$,  $\agg{\mb{x}}_0 = \mb{0}$.
Note that $\agg{\mb{x}}_t$ is determined after $\agg{\mb{s}}_t$ and $\agg{\mb{d}}_t$, within a set
\[
	\mathcal{P}( \agg{\mb{s}}_t, \agg{\mb{d}}_t, \agg{\mb{x}}_{t-1} )
	= \Bigg\{  \agg{\mb{x}}_t \in \mathbb{R}_+^p ~ \Bigg| ~	
	\mb{R}_S \agg{\mb{x}}_t \leq \agg{\mb{s}}_t ,~ 
	\mb{R}_D \agg{\mb{x}}_t \leq \agg{\mb{d}}_t ,~
	\agg{\mb{x}}_t \geq \agg{\mb{x}}_{t-1}
	\Bigg\},
\]
which is defined for $\agg{\mb{x}}_t$ to maximize profit, where the last constraint requires non-negative network activities. 
The main difference between single-period and multi-period models is that the order quantities are not \emph{static}.
That means, in order to obtain an optimal solution, one should find $\agg{\mb{s}}_t$ as a function of $\mb{D}_{[1:t-1]}$ so that they are \textit{fully-adjustable} to all previous demands.
Such policies also need to be \emph{non-anticipative}, i.e., adjustable decisions should only be based on realized uncertainties.

Even for $T=1$, the problem (\ref{eq:multi_obj}) is a two-stage robust optimization problem and shown to be NP-hard~\citep{ben2004}.
For multi-period setting, the complexity only worsens and, to our knowledge, no tractable algorithm has been proposed to exactly solve the general problem in (\ref{eq:multi_obj}).
Because of this, restrictions to specific policies have been considered.
In particular, \textit{affine policies} have been proposed, where adaptive decisions are assumed to be affine functions of realized uncertainties.

\begin{definition}[Affine Policy]
	\label{def:affine_policy}
	A policy is called an affine policy, if there exist ${\{ \mb{w}_t \in \mathbb{R}^m: ~ 1 \leq t \leq T\}}$ and $\{ \mb{W}_{\tau,t} \in \mathbb{R}^{m \times n} :~ 1 \leq \tau \leq t-1,~ 1 \leq t \leq T \}$ such that
	\begin{equation}
	\mb{s}_1 = \mb{w}_1 ,~~~ \mb{s}_t = 
	\mb{w}_t + \displaystyle \sum_{\tau=1}^{t-1} \mb{W}_{\tau,t} \mb{d}_\tau 	~~~~ t = 2,\ldots,T.
	\label{eq:def_fully_affine}
	\end{equation}
\end{definition}

Affine policies have exhibited excellent performance in many real-world applications.
With such policies, the multi-period problem \eqref{eq:multi_obj} converts to determining the affine weights.
These policies force non-anticipativity of $\mb{s}_t$ and one can reformulate (\ref{eq:multi_obj}) as a two-stage adaptive linear optimization problem
\begin{eqnarray}
\max\limits_{\mb{w}_t, \mb{W}_{\tau,t}}
\min\limits_{\mb{d_t}}
&&\max\limits_{\agg{\mb{x}}_t} \Bigg[
- \mb{c}^\top_S \agg{\mb{s}}_T 
- \mb{c}^\top_H  \sum_{t=1}^T \Big( \agg{\mb{s}}_t - \mb{R}_S \agg{\mb{x}}_t \Big)
- \mb{c}^\top_P \sum_{t=1}^{T} \Big( \agg{\mb{d}}_t - \mb{R}_D \agg{\mb{x}}_t \Big) + \mb{r}^{\top} \agg{\mb{x}}_T \Bigg]
\label{eq:obj_multi_FA}\\
\text{s.t.} &&
\begin{rcases}
\mb{w}_{t} + \displaystyle\sum_{\tau=1}^{t-1} \mb{W}_{\tau,t}\mb{d}_\tau \geq \mb{0} ,~~~~~~\mb{w}_1 \geq \mb{0}
\end{rcases}
~~\forall~ t=2,\ldots, T, ~\forall~ (\mb{d}_1,\ldots, \mb{d}_T) \in \mU^T
\label{eq:constr_multi_FA_first}\\
&&\begin{rcases}
\agg{\mb{s}}_1 = \mb{w}_1 \\
\agg{\mb{s}}_t = \mb{w}_1 + \displaystyle \sum_{j=2}^t \Big( \mb{w}_j + \sum_{\tau=1}^{j-1} \mb{W}_{\tau,j} \mb{d}_\tau \Big) \\
\mb{R}_S \agg{\mb{x}}_t ~\leq~ \agg{\mb{s}}_t 
~~~~~~~~~~~~~~~~~\\
\mb{R}_D \agg{\mb{x}}_t ~\leq~ \agg{\mb{d}}_t \\
\agg{\mb{x}}_T \geq \agg{\mb{x}}_{T-1} \geq \cdots \geq \agg{\mb{x}}_1 \geq \mb{0}
\end{rcases}
~~\forall t=1,\ldots, T.
~~~~~
\label{eq:constr_multi_FA_second}
\end{eqnarray}
Constraint (\ref{eq:constr_multi_FA_first}) implies that the order quantities are non-negative for any realizations of past demands, and constraint (\ref{eq:constr_multi_FA_second}) affects the inner maximization problem, which determines the processing activities.

\begin{proposition}
	Finding an optimal affine policy for a multi-period newsvendor network in~(\ref{eq:obj_multi_FA}--\ref{eq:constr_multi_FA_second}) is a convex optimization problem.
	\label{prop:multi_feasible}
\end{proposition}

\begin{remark}
The network activities $\mb{x}_t$ maximize the net profit over the entire horizon, not just at time $t$, i.e., we relax non-anticipativity of $\mb{x}_t$ in the optimization problem~(\ref{eq:obj_multi_FA}--\ref{eq:constr_multi_FA_second}).
However, we claim that this relaxation will not be loose, because penalty cost and holding cost force $\mb{x}_t$ to maximize profit in the corresponding period.
As a special case, one can show that in single-station models, $\mb{x}_t$ maximizes the overall profit if and only if it maximizes the profit at time $t$.
This relaxation facilitates generality, as problem (\ref{eq:obj_multi_FA}--\ref{eq:constr_multi_FA_second}) is defined for any polyhedral uncertainty sets, whereas in the stochastic case optimal strategies are only available for restricted cases (demands are i.i.d. over time as in~\cite{mieghem2002newsvendor}).
\end{remark}

	Since the inner minimization problem over $\mb{d}_t$ in (\ref{eq:obj_multi_FA}--\ref{eq:constr_multi_FA_second}) is non-convex, the overall problem is solved with cut generation.
	If the uncertainty set $\mU^T$ is a polyhedron, then it is guaranteed to find an optimal solution within finite number of iterations (cuts).
	Therefore, our method only requires a polyhedral structure of $\mU^T$.
	Note that the solution procedure does not exploit a specific structure of our uncertainty sets in Definition \ref{assumption:multi_period}, and the main purpose of using the budgeted uncertainty sets is to reduce conservatism.
	Even though an optimal solution can be obtained within finite iterations, the problem is still NP-hard \citep{ben2004} and the computation grows significantly as $T$ increases.

	While affine policies assume affine dependence of order quantities to realized demands, one can think of another class of policies, for which both ordering decisions and network activities are given as affine functions.
	We call these policies as \textit{affine-approximation policies}, which are defined below.

\begin{definition}[Affine-approximation Policy]
	\label{def:affine_approx}
	An adaptive policy is called an affine-approximation policy (Aff-approx) if there are ${\{ \mb{w}_t \in \mathbb{R}^m: ~ 1 \leq t \leq T\}}$ and $\{ \mb{W}_{\tau,t} \in \mathbb{R}^{m \times n} :~ 1 \leq \tau \leq t-1,~ 1 \leq t \leq T \}$, ${\{ \mb{y}_t \in \mathbb{R}^p: ~ 1 \leq t \leq T\}}$ and $\{ \mb{Y}_{\tau,t} \in \mathbb{R}^{p \times n} :~ 1 \leq \tau \leq t,~ 1 \leq t \leq T \}$ such that
		\begin{equation}
		\begin{aligned}
			&\mb{s}_1 = \mb{w}_1 ,~~~ \mb{s}_t = 
			\mb{w}_t + \displaystyle \sum_{\tau=1}^{t-1} \mb{W}_{\tau,t} \mb{d}_\tau 	\quad t = 2,\ldots,T \\
			&\mb{x}_t = \mb{y}_t + \sum_{\tau=1}^t \mb{Y}_{\tau,t} \mb{d}_\tau \quad\quad\quad\quad\quad\quad t=1,\ldots,T.
		\end{aligned}
		\label{eq:affine_approx}
		\end{equation}
\end{definition}

	Note that $\mb{x}_t$ depends on $\mb{D}_{[1:t]}$, while $\mb{s}_t$ is a function of $\mb{D}_{[1:t-1]}$, as network activities are assigned after the demand is realized at each time.
	We provide two observations for affine-approximation policies.
	\begin{enumerate}[label = (\roman*)]
		\item Aff-approx policies find affine parameters $\mb{y}_t$ and $\mb{Y}_{\tau,t}$, so that $\mb{x}_t$ is feasible to the innermost max operator.
		Thus, they give lower bounds to the affine policies.
		\item When Aff-approx policies are used, the problem (\ref{eq:obj_multi_FA}--\ref{eq:constr_multi_FA_second}) converts to a max-min problem, where all the affine weights (both for $\mb{s}_t$ and $\mb{x}_t$) are determined in the outer max operator.
		Using \eqref{eq:affine_approx}, both the objective function and the constraints can be expressed as functions of $\mb{w}_t$, $\mb{W}_{\tau,t}$, $\mb{y}_t$, $\mb{Y}_{\tau,t}$, and $\mb{d}_t$, which possibly include bilinear terms between the affine weights and $\mb{d}_t$.
		This type of problem is referred to as a (static) robust linear optimization problem, and they can be reformulated to a linear program, whenever $\mU^T$ is a polyhedron.
		Therefore, Aff-approx policies are tractable.
	\end{enumerate}

So far, we formulated multi-period robust newsvendor network problems and introduced two policies.
Affine policies convert the multi-period problem into a two-stage problem, which is computationally intractable.
On the contrary, Aff-approx policies solve a tractable linear program, but it only provides suboptimal solutions to affine policies  (we will numerically study their suboptimality in Section \ref{sec:experiments}).
Our main contribution is motivated by taking an alternative approach to these policies, as presented in the next section.

%%%%%%%%%%%%%%%%%
\section{Periodic-affine policies for single-station models}
\label{sec:single_prod}

As discussed, affine policies face computational difficulties when a decision maker has a larger number of resources and products over an extended period of time.
We propose a new solution concept, denoted as \emph{periodic-affine policies} (PA), where the overall time horizon is separated into subperiods, that are interconnected by the preceding surplus to become the proceeding demand.
In this approach, the order quantities are determined as an affine function of past demands realized only within its subperiod, as opposed to affine and affine-approximation policies where all previous demands are considered.
This scheme reduces the number of decision variables and consequently the computation time.
Our framework constructs this policy by first formulating a dynamic programming (DP) problem, where each stage corresponds to a subperiod.
We also propose an algorithm to compute such periodic-affine policies and show that they are computationally more tractable than affine policies. 
In addition, we present a sufficient condition that this algorithm provides the optimal solution to the DP problem.
We first consider $T$-period single-station models in this section.
However, our framework is naturally extended to multi-station networks which we discuss in the subsequent section.

\subsubsection*{Notation.}
We use same notations for all cost parameters $c_S, c_H, c_P$ with revenue per item, $r$, and we may assume that $R_S = R_D = 1$ without loss of generality in single-station models.
In this section $\mb{d}=(d_1,\ldots,d_T)$, $\mb{s}=(s_1,\ldots,s_T)$, and $\mb{x}=(x_1,\ldots,x_T)$.
We denote $\pi(w_t,W_{\tau,t})$ as an affine policy with affine parameters $\{w_t, W_{\tau,t}: 1\leq \tau \leq t-1, 1 \leq t \leq T\}$.
Furthermore, the problem of a $T$-period single-station newsvendor model is denoted as $\Phi(s_0,d_0)$ for an uncertainty set $\mU^T$ with initial input $s_0 \geq 0$ and demands $d_0 \geq 0$.

\subsubsection*{Analysis of initial input and demand.}
We first study the role of initial input and demand for the optimal affine policy in the multi-period model, given by 
\begin{equation*}
\Phi(s_0, d_0) := \underset{\pi}{\text{max}}~\underset{\mb{d} \in \mU^T}{\text{min}}~ \underset{\mb{x},\mb{s}\in \mathcal{X}(\pi, \mb{d}, s_0, d_0)}{\text{max}} P\Big(\pi(w_t,W_{\tau,t}),\mb{d},\mb{x}; s_0,d_0\Big),
\end{equation*}
where the profit during the period is
\begin{equation*}
\begin{aligned}
P\Big(\pi(w_t,W_{\tau,t}),\mb{d},\mb{x}; s_0,d_0\Big) = &-c_S \Bigg( \sum_{t=1}^T s_t \Bigg) -c_H \sum_{t=1}^T \Bigg( s_0 + \sum_{\tau=1}^t ( s_\tau - x_\tau) \Bigg) \\
&~~~~~~~~~~~~-c_P \sum_{t=1}^T \Bigg( d_0 + \sum_{\tau=1}^t ( d_\tau - x_\tau) \Bigg) + r \Bigg( \sum_{t=1}^T x_t \Bigg)
\end{aligned}
\end{equation*}
and the feasible set $\mathcal{X}(\pi, \mb{d},s_0, d_0)$ is given by
\begin{equation}
\mathcal{X}(\pi, \mb{d}, s_0, d_0) = \left\{ \mb{s}, \mb{x} \geq \mb{0} ~~ \middle|  ~~
\begin{aligned}
&s_1 = w_1,~s_t = w_t + \displaystyle \sum_{\tau=1}^{t-1} W_{\tau,t}d_\tau ~~~~\forall~ t=2,\ldots,T \\
&\displaystyle\sum_{\tau=1}^t x_\tau \leq s_0 + \sum_{\tau=1}^{t} s_\tau ~~~~~~~~~~~~~~~~\forall~ t=1,\ldots,T  \\
&\displaystyle\sum_{\tau=1}^t x_\tau \leq d_0 + \sum_{\tau=1}^{t} d_\tau ~~~~~~~~~~~~~~~~\forall~ t=1,\ldots,T
\end{aligned}
\right\}.
\label{eq:constr_single_fa}
\end{equation}
The result is intuitive and plays a key role in establishing periodic-affine policies.
\begin{proposition}
	For an optimal affine policy $\pi^*(w^*_t,W^*_{\tau,t})$ of $\Phi(0,0)$ with no initial input and demand, if $s_0 \leq w^*_1$, then:
	\begin{enumerate}[label = (\roman*)]
		\item An optimal affine policy $\lbar{\pi} = \lbar{\pi} (\lbar{w}_t^*,\lbar{W}_{\tau,t}^*)$ of $\Phi(s_0,d_0)$ is characterized as
		\begin{equation}
		\left\{
		\begin{aligned}
		&\lbar{w}^*_{1} = w_{1}^*- s_0 + d_0 \\
		&\lbar{w}_t^* = w_t^*~~~~~~~~~~~~~~~~~~~~~~~~~~~~\forall~ t=2,\ldots,T \\
		&\lbar{W}^*_{\tau,t} = W^*_{\tau,t}~~~~~~~~~~~~~~~~~~~~~~~~ \forall~ \tau=1, \ldots,t-1,~\forall~ t=1,\ldots,T.
		\end{aligned}
		\right.
		\label{eq:initial_single}
		\end{equation}
		\item There exists a single worst-case demand $\mb{d}^* \in \mU^T$ for both $\Phi(0,0)$ and $\Phi(s_0,d_0)$.
	\end{enumerate}
	\label{lemma:initial_single}
\end{proposition}
This result implies that for small enough connecting inventories, the subperiods become effectively decoupled.

\subsection{Model Formulation}
We now introduce the DP formulation for a multi-period newsvendor network.

\subsubsection*{Notation.}
For a $T$-period single-station model, we partition the time period into $N$ subperiods sorted as ${0 = t_0 < t_1 < \cdots <t_{N-1} < t_N = T}$.
In interval $I_j = \{t_{j-1}+1,\ldots, t_j\}$, the uncertainty set $\mU^j \in \mathbb{R}_+^{|I_j|}$ for every $j=1,\ldots,N$.
The amount of on-hand input stock and backlogged demands after time $t$ are $u_t^s$ and $u_t^d$.
A class of affine policies for $j^\text{th}$ subperiod is denoted by $\Pi_\text{aff}(\mU^j, \Xi_{j-1})$ on the uncertainty set $\mU^j$, where the state $\Xi_{j-1}$ contains all past information at the beginning of $j^\text{th}$ period with $\Xi_0 = 0$.
In this section, $\mb{d}_j= (d_{t_{j-1}+1},\ldots,d_{t_j}) \in \mathbb{R}_+^{|I_j|}$ denotes the demand at the $j$-th superiod for $j=1,\ldots,N$. We proceed similarly for $\mb{s}_j$ and $\mb{x}_j$.

\subsubsection*{DP formulation.}
We consider an $N$-stage robust DP problem, where each stage corresponds to each subperiod.
At the beginning of the $j^\text{th}$ subperiod, a decision maker obtains an affine policy $\pi_j \in \Pi_{\text{aff}}(\mU^j, \Xi_{j-1})$ to make adaptive ordering decisions for the current subperiod.
This can be formulated as
\begin{align}
\max\limits_{\pi_1 \in \Pi_{\text{aff}}(\mU^1,\Xi_0)} &
\Bigg[ \min\limits_{\mb{d}_1 \in \mU^1} 
\max\limits_{\mb{x}_1, \mb{s}_1 \in \mathcal{X}_1 }
\Bigg[ P_1 \Big( \pi_1,\mb{d}_1, \mb{x}_1; 0,0 \Big) \nonumber 
+ \!\!\!\!\!\!\!\!\!\max\limits_{\pi_2 \in \Pi_{\text{aff}}(\mU^2,\Xi_1)}  
\Bigg[ \min\limits_{\mb{d}_2 \in \mU^2}
\max\limits_{\mb{x}_2, \mb{s}_2 \in \mathcal{X}_2}
\bigg[ P_2 \Big(\pi_2, \mb{d}_2, \mb{x}_2; u_{t_1}^s,u_{t_1}^d \Big) \nonumber \\
&\cdots+ \max\limits_{\pi_N \in \Pi_{\text{aff}}(\mU^N,\Xi_{N-1})} 
\bigg[ \min\limits_{\mb{d}_N \in \mU^N}
\max\limits_{\mb{x}_N, \mb{s}_N\in \mathcal{X}_N}
P_N \Big(\pi_N ,\mb{d}_N, \mb{x}_N; u_{t_{N-1}}^s, u_{t_{N-1}}^d \Big) \bigg] \cdots \bigg] \Bigg] \Bigg] \Bigg], \label{eq:obj_PA_DP}
\end{align}
where $P_j \Big(\pi_j, \mb{d}_j, \mb{x}_j; u_{t_{j-1}}^s, u_{t_{j-1}}^d \Big)$ is a profit generated during the $j^\text{th}$ subperiod with an initial input and demand
\begin{equation}
\label{eq:profit}
\begin{split}
P_j \Big(\pi_j, \mb{d}_j, \mb{x}_j; u_{t_{j-1}}^s, u_{t_{j-1}}^d \Big)
= &-c_S \Bigg( \sum_{t \in I_j} s_t \Bigg) -c_H \sum_{t \in I_j} \Bigg( u_{t_{j-1}}^s + \sum_{\tau=t_{j-1}+1}^t ( s_\tau - x_\tau) \Bigg) \\
&-c_P \sum_{t \in I_j} \Bigg( u_{t_{j-1}}^d + \sum_{\tau=t_{j-1}+1 }^t ( d_\tau - x_\tau) \Bigg) + r \Bigg( \sum_{t \in I_j} x_t \Bigg),
\end{split}
\end{equation}
and $\mb{s}_j$ and $\mb{x}_j$ are determined within a feasible set $\mathcal{X}_j = \mathcal{X}(\pi_j, \mb{d}_j, u_{t_{j-1}}^s, u_{t_{j-1}}^d)$ from~\eqref{eq:constr_single_fa}.
In~\eqref{eq:profit}, $s_t$ denotes the order quantity at time $t$, if ordering decision is made by $\pi_j$ and $\mb{d}_j$ is realized.
The only constraint for $\pi_1,\ldots,\pi_N$ is to ensure non-negative order quantities for any demand realizations.

Recall that affine and Aff-approx policies are obtained before any demand realizations in (\ref{eq:obj_multi_FA}--\ref{eq:constr_multi_FA_second}), and the parameters are fixed over all time periods.
In the DP formulation \eqref{eq:obj_PA_DP}, the affine parameters are chosen dynamically at each subperiod, depending on the past information.
Specifically, the affine parameters $w_t^{(j)}(\Xi_{j-1})$ and $W_{\tau,t}^{(j)}(\Xi_{j-1})$ of $\pi_j$ can be any function of the past realization $\Xi_{j-1}$.
Hence, we can write the order quantities in \eqref{eq:obj_PA_DP} as
\begin{equation*}
s_{t}(\mb{d}_j,\Xi_{j-1}) = 
\begin{cases}
\displaystyle w_1^{(j)}(\Xi_{j-1}) \quad\quad\quad\quad\quad\quad\quad\quad\quad\quad\quad~~ t = t_{j-1} + 1 \\
\displaystyle w_i^{(j)}(\Xi_{j-1}) + \sum_{\tau=1}^{i-1} W_{\tau,i}^{(j)}(\Xi_{j-1})d_{t_{j-1}+\tau} \quad t = t_{j-1}+i,~ i \geq 2,~ t \in I_j.
\end{cases}
\end{equation*}
Every feasible affine and Aff-approx policy ensures non-negative orders, and hence, is also feasible in \eqref{eq:obj_PA_DP}.
This implies that by solving \eqref{eq:obj_PA_DP}, one can propose policies that have better worst-case profit than an optimal affine policy.

\subsubsection*{Periodic-affine policy formulation.}
With initial input $u^s$ bounded above with $w_1^{(j)}$, we define \emph{affine-IBS} policies by modifying the initial period of an affine policy $\pi_j(w_t^{(j)},W_{\tau,t}^{(j)})$.
\begin{definition}[Affine-IBS]
\label{def:AIBS}
For $j^\text{th}$ subperiod, the \textit{affine Initial Base-Stock policy} $\lbar{\pi}_j(w_t^{(j)}, W_{\tau,t}^{(j)})$ associated with 
an affine policy $\pi_j(w_t^{(j)}, W_{\tau,t}^{(j)})$ determines order quantity by
\begin{equation*}
s_t(u^s,u^d, \mb{d}_j) = 
\begin{cases}
	w_1^{(j)} - u^s + u^d ~~~~~~~~~~~~~~~ t = t_{j-1}+1 \\
	w_i^{(j)} + \displaystyle \sum_{\tau=1}^{i-1} W_{\tau,i}^{(j)}d_{t_{j-1}+\tau} ~~~~ t=t_{j-1}+i,~ i \geq 2,~ t \in I_j.
\end{cases}
\end{equation*}
\end{definition}
Note that at each subperiod, affine-IBS policies adapt to initial input and demand by adjusting the order quantity at the first period.
From the second period, affine-IBS and its associated affine policies are equivalent.
 
We now consider a sequence of affine-IBS policies $\lbar{\pi} = (\lbar{\pi}_1,\ldots, \lbar{\pi}_N)$, where each $\lbar{\pi}_j = \lbar{\pi}_j(w_t^{(j)}, W_{\tau,t}^{(j)})$ is for $j^\text{th}$ subperiod.
Note that this policy may not be well-defined for each subperiod because it does not guarantee that every order quantity is non-negative.
That means, if an input stock after $t_j$ is greater than $w_1^{(j+1)}$, then the policy would not be feasible.
To account for this, we impose
\begin{equation*}
w_{1}^{(j+1)} ~\geq~ u_{t_j}^s = u_{t_j}^s \Big( \lbar{\pi}_j,\mb{d}_j \Big)~~~~\forall \mb{d}_j \in \mU^j~ \forall j=1,\ldots,N-1.
\end{equation*}
By Definition \ref{def:AIBS}, the right-hand side is equivalent to
\begin{align*}
u_{t_j}^s \Big( \lbar{\pi}_j, \mb{d}_j \Big)
&= \text{max} \Bigg( 0,~ \bigg( u_{t_{j-1}}^s + \sum_{t \in I_j} s_t(u_{t_{j-1}}^s,u_{t_{j-1}}^d,\mb{d}_j) \bigg) - \bigg( u_{t_{j-1}}^d + \sum_{t \in I_j} d_t   \bigg) \Bigg)\\
&= \text{max} \Bigg( 0,~\bigg( \sum_{t = 1}^{t_j-t_{j-1}} w_t^{(j)} + \sum_{t=2}^{t_j-t_{j-1}} \sum_{\tau=1}^{t -1} W_{\tau,t}^{(j)} d_{t_{j-1}+\tau} \bigg) - \sum_{t = 1}^{t_j - t_{j-1}} d_{t_{j-1} + t} \Bigg).
\end{align*}
Since $w_1^{(j+1)}$ has to be non-negative for any demand realization, the periodic-affine policy is well-defined if
\begin{equation}
\label{prob:PA_algorithm_sub}
w_{1}^{(j+1)} ~\geq~  \theta_j^* := \underset{\mb{d}_j \in \mU^{j}}{\text{max}}~ \Bigg[
\sum_{t = 1}^{t_j-t_{j-1}} w_t^{(j)} + \sum_{t=2}^{t_j-t_{j-1}} \sum_{\tau=1}^{t -1} W_{\tau,t}^{(j)} d_{t_{j-1}+\tau} - \sum_{t = 1}^{t_j - t_{j-1}} d_{t_{j-1} + t} \Bigg]
\end{equation}
for every $j=0,\ldots,N-1$, where $\theta_j^*$ denotes the maximum leftover input after $j^\text{th}$ subperiod with $\theta_0^* = 0$.
Now we can define periodic-affine policies.

\begin{definition}[Periodic-affine Policy]
A \textit{periodic-affine policy} ${\lbar{\pi}_{\text{PA}} := (\lbar{\pi}_1,\ldots,\lbar{\pi}_N)}$ is an affine-IBS policy $\lbar{\pi}_i$  satisfying~\eqref{prob:PA_algorithm_sub} for affine policies $\pi_i$.
\label{def:PA}
\end{definition}

For a periodic-affine policy $\lbar{\pi}_{\text{PA}} = (\lbar{\pi}_1,\ldots, \lbar{\pi}_N)$, where $\lbar{\pi}_j = \lbar{\pi}_j(w_t^{(j)}, W_{\tau,t}^{(j)})$, order quantities at time $t$ are determined as follows:
\begin{equation}
\tag{PA}
\label{eq:PA_order_quantity}
s_t = s_t(\lbar{\pi}_{\text{PA}}) =
\begin{cases}
w_{1}^{(j)} - u_{t_{j-1}}^s + u_{t_{j-1}}^d ~~~~~~~~~~~ \forall~ t = t_{j-1} + 1 \\
w_{i}^{(j)} + \displaystyle \sum_{\tau=1}^{i-1} W_{\tau,i}^{(j)} d_{t_{j-1}+\tau} ~~~~~~~\forall~ t=t_{j-1} + i,~ i \geq 2,~ t \in I_j,~ \\
\end{cases}
\end{equation}
for every $j=1,\ldots, N$.
Since $\lbar{\pi}_j \in \Pi_{\text{aff}}(\mU^j, \Xi_{j-1})$ and $\lbar{\pi}_\text{PA}$ satisfies (\ref{prob:PA_algorithm_sub}), every periodic-affine policy is a feasible solution to the DP in~\eqref{eq:obj_PA_DP}.

In the next section, we present our algorithm to compute periodic-affine policies.

\subsection{Periodic-affine algorithm}
\label{subsec:PA_algorithm}

Our algorithm obtains affine-IBS policies for each subperiod by solving smaller subproblems.
However, since affine-IBS policies take initial input and demands into account, we construct the objective function to account for leftover resources and demands. % which will be backlogged to the next subperiods.
We identify such objective functions from the DP problem~\eqref{eq:obj_PA_DP}.
We first show that if initial input is small, an affine-IBS policy will be optimal among $\Pi_\text{aff}(\mU^N, \Xi_{N-1})$.
The proof is similar to the Proposition~\ref{lemma:initial_single} and is omitted.

\begin{corollary}
	\label{corr:AIBS_property}
	Let $\pi_N(w_t^{(N)},W_{\tau,t}^{(N)})$ be an optimal affine policy with zero initial input and demands.
	If $u_{t_{N-1}}^s \leq w^{(N)}_1$ for any realization of $u_{t_{N-1}}^s$, then its associated affine-IBS policy $\lbar{\pi}_N$ is an optimal solution among $\Pi_\emph{aff}(\mU^N, \Xi_{N-1})$.
	Moreover,
	\begin{eqnarray}
	&&\underset{\pi \in \Pi_{\emph{aff}}(\mU^N, \Xi_{N-1})}{\emph{max}} ~ 
	\underset{\mb{d}_N \in \mU^N}{\emph{min}} ~ 
	\underset{\mb{x}_N\in\mathcal{X}^N}{\emph{max}}~
	\bigg[ P_N \big( \pi, \mb{d}_N, \mb{x}_N; u_{t_{N-1}}^s, u_{t_{N-1}}^d \big) \bigg] \nonumber\\
	&=& c_S u_{t_{N-1}}^s + (r-c_S) u_{t_{N-1}}^d 
	+ \underset{\pi \in \Pi_{\emph{aff}}(\mU^N, 0)}{\emph{max}} ~ 
	\underset{\mb{d}_N \in \mU^N}{\emph{min}} ~ 
	\underset{\mb{x}_N\in\mathcal{X}^N}{\emph{max}}~
	\bigg[ P_N \big( \pi, \mb{d}_N, \mb{x}_N; 0,0 \big) \bigg].
	\label{eq:corollary}
	\end{eqnarray}
\end{corollary}

Using this Corollary, we reformulate an optimality condition for the last stage as
\begin{align*}
	\max_{\pi_{ N} \in \Pi_{\text{aff}}(\mU^N,\Xi_{N-1})}
	&\min_{\mb{d}_N \in \mU^N}  \max_{\mb{x}_N \in\mathcal{X}_N} 
	\Bigg[ P_N \Big( \pi_{N}, \mb{d}_N, \mb{x}_N; u_{t_{N-1}}^s, u_{t_{N-1}}^d \Big) \Bigg]  \\
	&= c_S  u_{t_{N-1}}^s + (r-c_S) u_{t_{N-1}}^d + 
	\max_{\pi_N \in \Pi_{\text{aff}}(\mU^N,0)}
	\min_{\mb{d}_N \in \mU^{N}}
	\max_{\mb{x}_N \in\mathcal{X}_N}
	\bigg[ P_N \Big( \pi_N, \mb{d}_N, \mb{x}_N; 0,0 \Big) \bigg].
\end{align*}
In single-station cases, $u_{t_{N-1}}^s$ and $u_{t_{N-1}}^d$ can be rewritten as
\begin{equation*}
	u_{t_{N-1}}^s = u_{t_{N-2}}^s + \sum_{t \in I_{N-1} } \Big( s_t - x_t \Big), ~~
	u_{t_{N-1}}^d = u_{t_{N-2}}^d + \sum_{t \in I_{N-1} } \Big( d_t - x_t \Big).
	\label{eq:PA_link}
\end{equation*}
This can be incorporated with $P_{N-1}\Big( \pi_{N-1},\mb{d}_{N-1},\mb{x}_{N-1}; u_{t_{N-2}}^s,u_{t_{N-2}}^d \Big)$ as
\begin{align}
&P_{N-1}\Big( \pi_{N-1},\mb{d}_{N-1},\mb{x}_{N-1}; u_{t_{N-2}}^s,u_{t_{N-2}}^d \Big) + c_S u_{t_{N-1}}^s + (r - c_S) u_{t_{N-1}}^d \nonumber \\
\begin{split}
	= -c_S \Bigg( \sum_{t \in I_{N-1}} s_t \Bigg) -c_H \sum_{t \in I_{N-1}} \Bigg( u_{t_{N-2}}^s + \sum_{\tau = t_{N-2}+1}^{t} ( s_\tau - x_\tau ) \Bigg) -c_P \sum_{t \in I_{N-1}} \Bigg( u_{t_{N-2}}^d+ \sum_{\tau = t_{N-2}+1}^{t} ( d_\tau - x_\tau ) \Bigg)
	\\
	+ ~ r \Bigg( \sum_{t \in I_{N-1}} x_t \Bigg) + c_S \Bigg( u_{t_{N-2}}^s + \sum_{t \in I_{N-1}} (s_t - x_t) \Bigg) + (r-c_S) \Bigg(u_{t_{N-2}}^d + \sum_{t \in I_{N-1}} (d_t - x_t) \Bigg)
\end{split}
\nonumber \\
\begin{split}
	=~ c_S u_{t_{N-2}}^s + (r - c_S) u_{t_{N-2}}^d -c_H \sum_{t \in I_{N-1}} \Bigg( u_{t_{N-2}}^s + \sum_{\tau = t_{N-2}+1}^{t} ( s_\tau - x_\tau ) \Bigg) 
	~~~~~~~~~~~~~~~~~~~~~~~~~~~~~~ \\
	-c_P \sum_{t \in I_{N-1}} \Bigg( u_{t_{N-2}}^d+ \sum_{\tau = t_{N-2}+1}^{t} ( d_\tau - x_\tau ) \Bigg) + (r-c_S) \Bigg( \sum_{t \in I_{N-1}} d_t \Bigg).
\end{split} \label{eq:reformulate}
\end{align}
We define a modified objective function $\widetilde{P}_{N-1}$ as 
\begin{equation}
\begin{split}
\widetilde{P}_{N-1} \Big( \pi_{N-1}, \mb{d}_{N-1}, \mb{x}_{N-1}; & u_{t_{N-2}}^s, u_{t_{N-2}}^d \Big)
= -c_H \sum_{t \in I_{N-1}} \Bigg( u_{t_{N-2}}^s + \sum_{\tau = t_{N-2}+1}^{t} ( s_\tau - x_\tau ) \Bigg) 
\\
&-c_P \sum_{t \in I_{N-1}} \Bigg( u_{t_{N-2}}^d+ \sum_{\tau = t_{N-2}+1}^{t} ( d_\tau - x_\tau ) \Bigg) + (r-c_S) \Bigg( \sum_{t \in I_{N-1}} d_t \Bigg).
\end{split}
\end{equation}
If we assume that both $u_{t_{N-2}}^s$ and $u_{t_{N-1}}^s$ are small, we can rewrite
\begin{align}
\begin{split}
& \max_{\pi_{N-1} \in \Pi_{\text{aff}}(\mU^{N-1},\Xi_{N-2})}  \min_{\mb{d}_{N-1} \in \mU^{N-1}} \max_{\mb{x}_{N-1} \in\mathcal{X}_{N-1}}
\Bigg[ P_{N-1}\Big( \pi_{N-1},\mb{d}_{N-1},\mb{x}_{N-1}; u_{t_{N-2}}^s,u_{t_{N-2}}^d \Big) 
\\
&~~~~~~~~~~~~~~~~~~~~~~~~~~~~~~~~~~~~~~~~
+ \max_{\pi_{ N} \in \Pi_{\text{aff}}(\mU^N,\Xi_{N-1})} \min_{\mb{d}_N \in \mU^N} \max_{\mb{x}_N \in\mathcal{X}_N}
	\bigg[ P_N \Big( \pi_N, \mb{d}_N,\mb{x}_N; u_{t_{N-1}}^s, u_{t_{N-1}}^d \Big) \bigg] \Bigg]
\end{split} 
\nonumber \\
\begin{split}
	&= \max_{\pi_{N-1} \in \Pi_{\text{aff}}(\mU^{N-1},\Xi_{N-2})} 
	\min_{\mb{d}_{N-1} \in \mU^{N-1}} \max_{\mb{x}_{N-1} \in\mathcal{X}_{N-1}}
	\Bigg[ P_{N-1} \Big( \pi_{N-1},\mb{d}_{N-1}, \mb{x}_{N-1},u_{t_{N-2}}^s,u_{t_{N-2}}^d \Big) \\
	&~~~~~~~~~~~~~~~~~~
	+ c_S u_{t_{N-1}}^s + (r - c_S) u_{t_{N-1}}^d
	+ \max_{\pi_N \in \Pi_{\text{aff}}(\mU^N,0)} 
	\min_{\mb{d}_N \in \mU^{N}} \max_{\mb{x}_N \in\mathcal{X}_N}
	\bigg[ P_N \Big( \pi_N,\mb{d}_N, \mb{x}_N; 0,0 \Big) \bigg] \Bigg]
\end{split}
\nonumber \\
\begin{split}
&= \max_{\pi_{N-1} \in \Pi_{\text{aff}}(\mU^{N-1},\Xi_{N-2})} 
\min_{\mb{d}_{N-1} \in \mU^{N-1}} \max_{\mb{x}_{N-1} \in\mathcal{X}_{N-1}}
\Bigg[ c_S u_{t_{N-2}}^s + (r - c_S) u_{t_{N-2}}^d \\
&~~~~~~~~+ \widetilde{P}_{N-1}\Big( \pi_{N-1},\mb{d}_{N-1}, \mb{x}_{N-1}; u_{t_{N-2}}^s,u_{t_{N-2}}^d \Big) \Bigg]
+ \max_{\pi_N \in \Pi_{\text{aff}}(\mU^N,0)} 
\min_{\mb{d}_N \in \mU^{N}} \max_{\mb{x}_N \in\mathcal{X}_N}
\bigg[ P_N \Big( \pi_N, \mb{d}_N, \mb{x}_N; 0,0 \Big) \bigg]
\end{split}
\nonumber \\
\begin{split}
&=~ c_S u_{t_{N-2}}^s + (r - c_S) u_{t_{N-2}}^d +\!\!\!\!\!\!\!\!
\max_{\pi_{N-1} \in \Pi_{\text{aff}}(\mU^{N-1},0)}  \min_{\mb{d}_{N-1} \in \mU^{N-1}} \max_{\mb{x}_{N-1} \in\mathcal{X}_{N-1}} \bigg[ \widetilde{P}_{N-1} \Big( \pi_{N-1},\mb{d}_{N-1},\mb{x}_{N-1}; 0,0 \Big) \bigg] \\
&\quad\quad + \max_{\pi_{N} \in \Pi_{\text{aff}}(\mU^{N},0)} 
\min_{\mb{d}_N \in \mU^{N}} \max_{\mb{x}_N \in\mathcal{X}_N} \bigg[ P_N \Big( \pi_N,\mb{d}_N, \mb{x}_N; 0,0 \Big) \bigg] .
\end{split} 
\label{eq:reformulate_final}
\end{align}
Note that the second equality comes from (\ref{eq:reformulate}), and one can verify similarly from Corollary \ref{corr:AIBS_property} that the last equality holds.
This reformulation shows that if leftover input after every subperiod is small enough, we can solve the DP problem in~\eqref{eq:obj_PA_DP} by solving smaller subproblems. 
These subproblems are defined with modified objective function $\widetilde{P}_j$ with no backlogged input and demand, hence we can solve them independently.
Proceeding iteratively, we define an objective $P_j^\text{PA} \Big( \pi_j, \mb{d}_j, \mb{x}_j \Big)$ as
\begin{eqnarray}
P_j^\text{PA} \Big( \pi_j, \mb{d}_j, \mb{x}_j \Big) =
\begin{cases}
\begin{split}
-c_H \displaystyle\sum_{t \in I_{j}} \Big( \sum_{\tau=t_{j-1}+1}^{t}(s_\tau - x_\tau) \Big)
-c_P \sum_{t \in I_{j}} \Big(\sum_{\tau=t_{j-1}+1}^{t}(d_\tau-x_\tau)	\Big) \\
+ (r-c_S) \sum_{t \in I_j} d_t,
\end{split}
~~~~j=1,\ldots,N-1
\\
\begin{split}
-c_S \displaystyle\sum_{t \in I_j} s_t
-c_H \displaystyle\sum_{t \in I_{j}} \Big( \sum_{\tau=t_{j-1}+1}^{t}(s_\tau - x_\tau) \Big) 
~~~~~~~~~~~~~~~~~~\quad\quad\\
-c_P \sum_{t \in I_{j}} \Big(\sum_{\tau=t_{j-1}+1}^{t}(d_\tau-x_\tau)	\Big)
+ r \sum_{t \in I_j} x_t,
\end{split}
~~~~j=N.
\end{cases}
\label{eq:obj_PA_new}
\end{eqnarray}

%\subsubsection*{Periodic-affine algorithm.}
We now propose the periodic-affine algorithm.
This algorithm (i) ensures that the solution is well-defined, and (ii) exploits the modified objective functions $P_j^\text{PA}$.
The $j^\text{th}$ subproblem can be solved by the following optimization problem
\begin{equation}
\begin{aligned}
&\underset{\pi_j \in \Pi_\text{aff}(\mU^j,0)}{\text{max}} ~ \underset{\mb{d}_j \in \mU^j}{\text{min}} ~ \underset{\mb{x}_j, \mb{s}_j}{\text{max}} ~~ P_j^\text{PA} \Big( \pi_j, \mb{d}_j, \mb{x}_j \Big) \\
&~~~~\text{s.t. }~~~~ (\mb{x}_j, \mb{s}_j) \in \mathcal{X}(\pi_j, \mb{d}_j, 0, 0) \\
& ~~~~~~~~~~~~~w_1^{(j)} \geq \theta_{j-1}^* .\\
\end{aligned}
\label{prob:PA_algorithm_main}
\end{equation}
The last constraint ensures that periodic-affine policy is well-defined, where $\theta_{j-1}^*$ is the maximum amount of on-hand input after $(j-1)^\text{th}$ subperiod, computed by~\eqref{prob:PA_algorithm_sub}.
The overall procedure solves (\ref{prob:PA_algorithm_main}) and (\ref{prob:PA_algorithm_sub}) iteratively, as summarized in Algorithm \ref{alg:PA_single}.

\begin{algorithm}[H]
\caption{Periodic-affine algorithm for single-station problems}
	\emph{Given.} time indices $0 = t_0 < t_1 < \cdots < t_N = T$, uncertainty set $\mU = \mU^1 \times \cdots \times \mU^N$, $j=1$, and $\theta_0^*=0$. \\
	\emph{Step 1.} Solve (\ref{prob:PA_algorithm_main}) to obtain $\pi_j$ for the $j^\text{th}$ subperiod. \\
	\emph{Step 2.} Using $\pi_j$, compute the maximum leftover input $\theta_j^*$ by (\ref{prob:PA_algorithm_sub}). \\
	\emph{Step 3.} If $j=N$, return $\lbar{\pi}_\text{PA}= (\lbar{\pi}_1, \ldots, \lbar{\pi}_N)$ and STOP. Otherwise, $j \leftarrow j+1$ and go to \emph{Step 1}.
	\label{alg:PA_single}
\end{algorithm}

\noindent We present an additional useful property of the periodic-affine algorithm, namely that the worst-case scenario can be obtained from each iteration.
\begin{proposition}
	\label{prop:PA_worst_scenario}
	Let $\lbar{\pi}_\text{PA}$ be a solution of periodic-affine algorithm and $\mb{d}^*_j$ be a worst-case scenario from the $j^\text{th}$ subproblem.
	Then $\lbar{\pi}_\text{PA}$ has a worst-case scenario $(\mb{d}^*_1,\ldots,\mb{d}^*_N)$.
\end{proposition}

So far, we discussed a single-station, multi-period robust newsvendor model, where the uncertainty set over time periods is defined as a Cartesian product of pre-specified uncertainty sets for each subperiod.
We formulated a dynamic programming problem, where each stage corresponds to each subperiod.
Motivated from this formulation, we developed an algorithm to find a periodic-affine policy, by defining the modified objective functions $P_j^\text{PA}$'s.

\subsubsection*{Computational Advantages.}
We provide two theoretical evidences that periodic-affine policies are computationally advantageous over affine policies.
First, the number of affine parameters in PA grows linearly in $T$, while it increases quadratically in affine policies.
Second, solving subproblems in PA requires significantly less computational effort than solving one large problem in affine policies.
This is because an optimal affine policy solves  (\ref{eq:obj_multi_FA}--\ref{eq:constr_multi_FA_second}) by cut generation.
The number of extreme points in multi-period uncertainty sets (Definition \ref{assumption:multi_period}) grows exponentially in $T$, and this results in a very large number of iterations for affine policies.
In PA, however, one can arbitrarily choose the length of subperiods.
Hence, the number of iterations grows at most linearly in $T$.
Strong empirical evidence in Section \ref{sec:experiments} demonstrate these advantages.

Also note that PA does not suffer from the curse of dimensionality.
Even though PA is motivated from the DP formulation \eqref{eq:obj_PA_DP}, an optimal solution is not obtained by solving the Bellman equation; our algorithm decomposes the overall problem into smaller subproblems in order to achieve tractable solutions.

In the next section, we present theoretical properties of periodic-affine policies, where we provide a sufficient condition for the algorithm to have an optimal solution to the DP \eqref{eq:obj_PA_DP}.

\subsection{Optimality of periodic-affine policies}

In this section, we present theoretical properties of PA by analyzing the effect of base-stock levels on the worst-case performance.
Specifically, we provide a sufficient condition under which the periodic-affine algorithm solves the DP problem (\ref{eq:obj_PA_DP}).
To compare the worst-case performance of PA with affine policies, we consider affine policies under the rectangular uncertainty set $\mU = \mU^1 \times \cdots \times \mU^N$ so that both policies are defined equivalently.
Moreover, we present an analytical approximation for the suboptimality of PA.
Note that this is a \emph{posterior} approximation, i.e., it is computed during the algorithm.

Let the worst-case profit of the two policies be $V_\text{PA}^*$ and $V_\text{Aff}^*$, which are evaluated by \eqref{eq:obj_PA_DP}, and $V_\text{DP}^*$ as an optimal value of the DP problem \eqref{eq:obj_PA_DP}.

The following assumption guarantees the optimality of PA policies.

\begin{assumption}
	\label{assumption:optimality_PA}
	For a solution of the PA algorithm $\lbar{\pi}_{\emph{PA}}$, assume that the maximum leftover input level after each subperiod $\theta_j^*$ satisfies the last constraint in (\ref{prob:PA_algorithm_main}).
\end{assumption}
In other words, for a solution of the PA algorithm, the last constraint in \eqref{prob:PA_algorithm_main} is not active at every iteration.
	Since an optimal periodic-affine policy maximizes the overall profit, it tends to have lower leftovers at each time period and hence, is likely to satisfy Assumption \ref{assumption:optimality_PA}.
	This assumption is not too restrictive as evidenced in the numerical experiments in Section~\ref{sec:experiments}.
\begin{remark}
We also suggest that it is possible to enforce the Assumption~\ref{assumption:optimality_PA} to hold, or reduce the suboptimality of PA, by a proper choice of partitioning the time horizon.
For example, consider a problem with $T=12$, where the nominal mean is 10 for the first 6 periods and decreases to 2 for the last 6 periods.
Then, Assumption \ref{assumption:optimality_PA} may not hold if PA is obtained by partitioning 12 periods into two (of each 6 periods).
However, one can obtain a better result by partitioning the periods into 3 subperiods of each 4 periods.	
\end{remark}

\begin{theorem}
\label{thm:optimality_PA}
For a single-station network, if Assumption \ref{assumption:optimality_PA} holds, then $V_\emph{Aff}^* \leq V_\emph{PA}^* = V_\emph{DP}^*$.
\end{theorem}

Since affine policies are defined on a subset of every feasible solution of the DP, the worst-case profit of PA is guaranteed to be greater than or equal to that of affine policies under Assumption \ref{assumption:optimality_PA}.
We show in the following proposition that in single-station problems, their worst-case performance is indeed equal.

\begin{proposition}
	\label{proposition:PA_vs_FA}
	For single-station models under Assumption \ref{assumption:optimality_PA}, $V_\emph{Aff}^* = V_\emph{PA}^* = V_\emph{DP}^*$.
\end{proposition}

\subsubsection*{Remark.}
All previous theoretical properties of PA rely on Assumption \ref{assumption:optimality_PA}.
In other words, the worst-case performance of PA may not match that of affine policies without the assumption.

We now relax Assumption \ref{assumption:optimality_PA} and provide a suboptimality bound for PA.
For this, we compare the solutions of the following optimization problems

\begin{minipage}{.5\linewidth}
	\begin{equation*}
	\begin{aligned}
	\mathcal{P}_j := &\!\!\!\!\!\!\max_{\pi_j \in \Pi_\text{aff}(\mU^j,0)} \min_{\mb{d}_j \in \mU^j} \max_{\mb{x}_j,\mb{s}_j} ~~ P_j^\text{PA} \Big( \pi_j, \mb{d}_j, \mb{x}_j \Big) \\
	&~~~~\text{s.t. }~~~~ (\mb{x}_j, \mb{s}_j) \in \mathcal{X}(\pi_j, \mb{d}_j, 0, 0) \\
	& ~~~~~~~~~~~~~w_1^{(j)} \geq \theta_{j-1}^* \\
	\end{aligned}
	\label{prob:with_S}
	\end{equation*}
\end{minipage}%
\begin{minipage}{.5\linewidth}
	\begin{equation*}
	\begin{aligned}
	\widetilde{\mathcal{P}}_j := &\!\!\!\!\!\!\max_{\pi_j \in \Pi_\text{aff}(\mU^j,0)} \min_{\mb{d}_j \in \mU^j} \max_{\mb{x}_j,\mb{s}_j} ~~ P_j^\text{PA} \Big( \pi_j, \mb{d}_j, \mb{x}_j \Big) \\
	&~~~~\text{s.t. }~~~~ (\mb{x}_j, \mb{s}_j) \in \mathcal{X}(\pi_j, \mb{d}_j, 0, 0), \\
	\\
	\end{aligned}
	\label{prob:no_S}
	\end{equation*}
\end{minipage}
where $\mathcal{P}_j$ is the $j^\text{th}$ subproblem during the PA algorithm and $\widetilde{\mathcal{P}}_j$ solves a subproblem with assuming no leftover input from the previous subperiods.

\begin{theorem}
	\label{thm:suboptimal_PA}
	For any single-station newsvendor networks with an objective value $\widetilde{f}_j^*$ of $\widetilde{\mathcal{P}}_j$, 
	\begin{equation*}
	V_{\emph{PA}}^* ~\leq~ V_{\emph{Aff}}^* ~\leq~ V_{\emph{DP}}^* ~\leq~ \widetilde{V}_\text{PA}^* := \sum_{j=1}^N \widetilde{f}_j^*.
	\label{eq:suboptimal_PA}
	\end{equation*}
\end{theorem}

Theorem \ref{thm:suboptimal_PA} provides a tight bound.
All the inequalities hold with equalities if Assumption \ref{assumption:optimality_PA} holds.
Note that one may not need to resolve a problem $\widetilde{\mathcal{P}}_j$; recall that every subproblem in the periodic-affine algorithm is solved by generating cuts.
Once a solution of $\mathcal{P}_j$ is obtained, one can relax the last constraint in $\mathcal{P}_j$ and continue cut generation in order to solve $\widetilde{\mathcal{P}}_j$.
This requires fewer iterations, since (i) an optimal solution of $\mathcal{P}_j$ serves as a warm-start initial point for additional cuts, and (ii) the previously generated cuts are still valid without any modifications to $\widetilde{\mathcal{P}}_j$.

So far, we introduced theoretical properties of periodic-affine policies.
We showed that under mild condition, the algorithm finds an optimal solution to the DP problem, and thus achieve a worst-case performance equal to that of an optimal affine policy for single-station problems.
If this assumption does not hold, we provided a tight bound that can measure the gap between periodic and affine policies.
Moreover, this gap can be computed by minimally modifying the algorithm with similar computational requirements.
In the proceeding section, we extend this framework to general multi-station newsvendor networks and infinite-horizon problems.

%%%%%%%%%%%%%%%%%
\section{Extensions of Periodic-affine Policies}
\label{sec:RNN_multi_prod}
We extend our approach in Section \ref{sec:single_prod} to general multi-station networks, where a decision-maker intends to satisfy customers' demand at multiple locations.
First, we extend Algorithm~\ref{alg:PA_single} for multi-station networks. 
Then, we develop periodic-affine policies for infinite horizon problems.

\subsection{Multi-station networks}

Here, we follow the flow of Section \ref{sec:single_prod}. 
To set this up, we define a matrix that plays a key role in implementing periodic-affine policies.

\begin{definition}[Basic Matrix]
\label{def:basic_matrix}
Let $\ell_k \in \mathcal{L}$ be an optimal solution to satisfy a unit demand at a sink node $k \in \mathcal{S}$.
A \emph{basic matrix} $\mb{R}_B \in \mathbb{R}^{p \times n}$ is given by
	\begin{equation*}
		\mb{R}_B(\ell,k) = 1~ \text{if } \ell = \ell_k \; \forall k, \text{ and } 
			~0 ~\text{otherwise.}
	\end{equation*}
\end{definition}

Using such a basic matrix, we obtain a closed-form expressions of ordering quantities and network activities if the demand is deterministic.
In particular, for any $\mb{d} \in \mathbb{R}^n_+$, an optimal decision is given by $\mb{s} = \mb{R}_S \mb{R}_B \mb{d}$ and $\mb{x} = \mb{R}_B \mb{d}$.

Recall that we have defined a modified objective functions for each stage of the DP problem (\ref{eq:obj_PA_DP}) to separate the overall problem into subproblems.
In single-station models, the values of on-hand products and backlogged demands at the beginning of the $(j+1)^{\text{th}}$ subperiod are expressed as
\begin{align*}
c_S u_{t_j}^s &~=~ c_S \cdot \sum_{t \in I_{j}} \Big( s_t - x_t \Big) \\
(r-c_S) u_{t_j}^d & ~=~ (r-c_S) \cdot \sum_{t \in I_{j}} \Big( d_t - x_t \Big),
\end{align*}
which are taken into $j^\text{th}$ subproblem.
After the $j^\text{th}$ subperiod, $\mb{u}_{t_j}^\mb{s}$ and $\mb{u}_{t_j}^\mb{d}$ are deterministic and hence, their values can be expressed by using the basic matrix $\mb{R}_B$ as
\begin{eqnarray*}
&\displaystyle \mb{c}^\top_S  \mb{u}_{t_j}^\mb{s}  = \mb{c}^\top_S ~\sum\limits_{t \in I_{j}} \Big( \mb{s}_t - \mb{R}_S \mb{x}_t \Big) \\ 
&\displaystyle \Big( \mb{R}_B^\top \mb{r} - \mb{R}_B^\top \mb{R}_S^\top \mb{c}_S \Big)^\top \mb{u}_{t_j}^\mb{d}
= \Big( \mb{R}_B^\top \mb{r} - \mb{R}_B^\top \mb{R}_S^\top \mb{c}_S \Big)^\top \sum\limits_{t \in I_j} \Big( \mb{d}_t - \mb{R}_D \mb{x}_t \Big).
\end{eqnarray*}
Note that the value of $\mb{u}_{t_j}^\mb{d}$ is determined by ordering $\mb{R}_B \mb{u}_{t_j}^\mb{d}$ and processing $\mb{R}_S \mb{R}_B \mb{u}_{t_j}^\mb{d}$.
This allows us to extend the definition of affine-IBS policies as follows.
\begin{definition}[Affine-IBS for Multi-station]
	\label{def:AIBS_multi}
	For $j^\text{th}$ subperiod, the affine-IBS policy $\lbar{\mbs{\pi}}_j(\mb{w}_t^{(j)}, \mb{W}_{\tau,t}^{(j)})$ associated with an affine policy $\mbs{\pi}_j(\mb{w}_t^{(j)}, \mb{W}_{\tau,t}^{(j)})$ determines order quantity by
	\begin{equation*}
	\mb{s}_t(\mb{u}^\mb{s},\mb{u}^\mb{d}, \mb{D}_j) = 
	\begin{cases}
	\mb{w}_1^{(j)} - \mb{u}^\mb{s} + \mb{R}_S \mb{R}_B \mb{u}^\mb{d} ~~~~~~~~ t = t_{j-1}+1 \\
	\mb{w}_i^{(j)} + \displaystyle \sum_{\tau=1}^{i-1} \mb{W}_{\tau,i}^{(j)}\mb{d}_{t_{j-1}+\tau} ~~~~~~ t=t_{j-1}+i,~ i \geq 2,~ t \in I_j.
	\end{cases}
	\end{equation*}
\end{definition}

\subsubsection*{Period-affine policy for multi-station networks.} 
As in Definition \ref{def:PA}, periodic-affine policies are defined as a sequence of affine-IBS policies.
Eq. (\ref{prob:PA_algorithm_sub}), which is required for periodic-affine policies to be well-defined, is readily extended by replacing with a vector inequality.
With this generalization, all the arguments in Section \ref{subsec:PA_algorithm} can be repeated in multi-station network setting.

As a result, the objective function for each subproblem is given by
\begin{eqnarray}
\begin{split}
\mb{P}_j^\text{PA} \Big( \mbs{\pi}_j, \mb{D}_{j}, \mb{X}_{j} \Big) =
\end{split}
\begin{cases}
\begin{split}
&-\mb{c}^\top_H \displaystyle\sum_{t \in I_{j}} \Big( \sum_{\tau=t_{j-1}+1}^{t}(\mb{s}_\tau - \mb{R}_S \mb{x}_\tau) \Big) \\
&~~~~~~~~-\mb{c}^\top_P \sum_{t \in I_{j}} \Big(\sum_{\tau=t_{j-1}+1}^{t}(\mb{d}_\tau- \mb{R}_D \mb{x}_\tau)	\Big)
 + \mb{v}^\top_\mb{d} \sum_{t \in I_j} \mb{d}_t + \mb{v}^\top_\mb{x} \sum_{t \in I_j} \mb{x}_t,
\end{split}
~~~ j \leq N-1
\\
\begin{split}
&-\mb{c}^\top_S \displaystyle\sum_{t \in I_j} \mb{s}_t
-\mb{c}^\top_H \displaystyle\sum_{t \in I_{j}} \Big( \sum_{\tau=t_{j-1}+1}^{t}(\mb{s}_\tau - \mb{R}_S \mb{x}_\tau) \Big) \\
& ~~~~~~~~-\mb{c}^\top_P \sum_{t \in I_{j}} \Big(\sum_{\tau=t_{j-1}+1}^{t}(\mb{d}_\tau-\mb{R}_D \mb{x}_\tau)	\Big)
+ \mb{r}^\top \sum_{t \in I_j} \mb{x}_t,
\end{split}
~~~~~~~~~~~~~~~~~~~j=N
\end{cases}
\label{eq:obj_PA_multi}
\end{eqnarray}
where $\mb{v}_\mb{d} = \mb{R}^\top_B \mb{r} - \mb{R}^\top_B \mb{R}^\top_S \mb{c}_S$ and $\mb{v}_\mb{x} = \mb{r} - \mb{R}^\top_S \mb{c}_S - \mb{R}^\top_D \mb{R}^\top_B \mb{r} + \mb{R}^\top_D \mb{R}^\top_B \mb{R}^\top_S \mb{c}_S$.

\subsubsection*{Period-affine algorithm for multi-station networks.} 
As in Section~\ref{sec:single_prod}, Problem~\eqref{prob:PA_algorithm_main} can be readily converted into multidimensional form by replacing the last inequality with a vector inequality.
However, it is challenging to obtain the multi-station version of (\ref{prob:PA_algorithm_sub}), which computes the maximum amount of leftover resources.
Therefore, we incorporate these into a single robust two-stage optimization problem, as follows:
\begin{eqnarray*}
\max_{\mbs{\theta}_j, \mbs{\pi}_j \in \mbs{\Pi}_\text{aff}(\mU^j,\mb{0})} && \min_{\mb{D}_{j} \in \mU^j} \max_{\mb{S}_{j},\mb{X}_{j}} ~~ \mb{P}_j^\text{PA} \bigg( \mbs{\pi}_j, \mb{D}_{j}, \mb{X}_{j} \bigg) - \delta \cdot \mb{1}^\top \mbs{\theta}_j \\
&&
\begin{aligned}
\!\!\!\!\!\!\!\!\!\!\!\!\!\!\!\!\!\!\!\text{s.t.} ~~~~
& (\mb{S}_{j}, \mb{X}_{j}) \in \mathcal{X}(\mbs{\pi}_j, \mb{D}_{j},\mb{0},\mb{0}) \nonumber \\
& \mb{w}_1^{(j)} ~\geq~ \mbs{\theta}_{j-1}^* \\
&\sum_{t \in I_j} \left( \mb{s}_t - \mb{R}_S \mb{x}_t \right) ~\leq~ \mbs{\theta}_j
\end{aligned}
 \nonumber\\
%\label{eq:constr_multi_PA_second}
\end{eqnarray*}
where $\delta > 0$ is a small real number.
In this way, the periodic-affine algorithm for multi-station networks proceeds by iteratively solving subproblems, similar to Algorithm \ref{alg:PA_single}.

\subsubsection*{Properties of periodic-affine for multi-station networks.}
We now generalize the theoretical properties for the multi-station networks.
We consider the DP problem (\ref{eq:obj_PA_DP}) by replacing every single-dimensional quantity by multi-dimensional quantities.
Assumption \ref{assumption:optimality_PA} is extended with vector inequalities, each of which is for each source node.
We use $V_\text{Aff}^*$, $V_\text{PA}^*$, and $V_\text{DP}^*$ for the worst-case objective values for affine, PA, and the DP problem, and define $\widetilde{V}_\emph{PA}^*$ similar to single-station problems.
\begin{theorem}
\label{corollary:optimality_PA_multi}
For multi-station networks, if Assumption \ref{assumption:optimality_PA} holds, then $V_\emph{Aff}^* \leq V_\emph{PA}^* = V_\emph{DP}^*$.
Otherwise, $V_{\emph{PA}}^* ~\leq~ V_{\emph{DP}}^* ~\leq~ \widetilde{V}_\text{PA}^*$.
\end{theorem}

Note that Proposition \ref{proposition:PA_vs_FA} cannot be extended to multi-station networks.
In other words, PA policies for multi-station networks are not necessarily be an affine policy.
Theorem \ref{corollary:optimality_PA_multi} implies that an optimal PA policy has a worst-case performance not less than an optimal affine policy.
However, for multi-station networks, we cannot compare the two policies without Assumption \ref{assumption:optimality_PA}.

\subsection{Infinite horizon problems}

So far, PA policies are based on multi-period problems of finite horizon.
In this section, we extend these PA policies to infinite horizon problems with a discount factor of $\beta < 1$. 
For this, we assume that nominal means and covariances of demands have periodicity with the period $T_0 \geq 1$.
We then define an uncertainty set ${\mU^{T_0} \in \mathbb{R}^{n \times {T_0}}}$ to describe demand uncertainties for each period.
This framework models settings where demand has stationary mean and covariance along the periods.

\begin{definition}[infinite-horizon uncertainty set] 
The \emph{infinite-horizon demand uncertainty set} is a Cartesian product of $\mU^{T_0}$ via
\[
\mU^\infty = \mU^{T_0} \times  \mU^{T_0} \times \mU^{T_0} \times \cdots.
\]
\end{definition}

We implement PA policies for infinite horizon problems by \textit{replicating} policies over the periods.
We construct a PA policy of period ${T_0}$ by solving a single problem of duration ${T_0}$.
As in previous sections, we define the objective function $\mb{P}_\infty^\text{PA}(\mbs{\pi}, \mb{d}, \mb{x})$, by taking leftover inventories, unsatisfied demands, and the discount factor into account as 
\begin{align*}
	\mb{P}_\infty^\text{PA}&\Big( \mbs{\pi},\mb{D}_{[1:{T_0}]},\mb{X}_{[1:{T_0}]} \Big) \\
	&= - \mb{c}^\top_H \displaystyle\sum_{t=1}^{T_0} \beta^t\Big(\sum_{\tau=1}^{t}(\mb{s}_\tau - \mb{R}_S \mb{x}_\tau) \Big) 
	-\mb{c}^\top_P \sum_{t=1}^{T_0} \beta^t \Big(\sum_{\tau=1}^{t}(\mb{d}_\tau- \mb{R}_D \mb{x}_\tau)	\Big) + \mb{r}^\top \sum_{t=1}^{T_0} \beta^t \mb{x}_t\\
	&~~~- \mb{c}^\top_S \sum_{t=1}^{T_0} \beta^t \mb{s}_t + \beta^{{T_0}+1} \mb{c}^\top_S \sum_{t=1}^{T_0} ( \mb{s}_t - \mb{R}_S \mb{x}_t) + \beta^{{T_0}+1} \Big( \mb{R}^\top_B \mb{r}- \mb{R}^\top_B \mb{R}^\top_S \mb{c}_S \Big)^\top \sum_{t=1}^{T_0} (\mb{d}_t - \mb{R}_D \mb{x}_t) \\
&= - \sum_{t=1}^{T_0} \beta^t ~ \mb{c}^{*\top}_{S,t} \mb{s}_t - \displaystyle\sum_{t=1}^{T_0} \beta^t ~ \mb{c}^{*\top}_{H,t} \Big(\sum_{\tau=1}^{t}(\mb{s}_\tau - \mb{R}_S \mb{x}_\tau) \Big)  - \sum_{t=1}^{T_0} \beta^t ~ \mb{c}^{*\top}_{P,t} \Big(\sum_{\tau=1}^{t}(\mb{d}_\tau- \mb{R}_D \mb{x}_\tau)	\Big) + \sum_{t=1}^{T_0} \beta^t ~ \mb{r}^{*\top}_{t} \mb{x}_t,
\end{align*}
where $\mb{c}_{S,t}^* = \mb{c}_S$, $\mb{r}_{t}^* = \mb{r}$ for $t=1,\ldots,{T_0}$, and 
\begin{equation*}
\mb{c}_{H,t}^* = 
\begin{cases}
\mb{c}_H ~~~~~~~~~~~~ 1 \leq t \leq {T_0}-1 \\
\mb{c}_H - \beta \mb{c}_S ~~~~ t = {T_0}
\end{cases}\quad
\mb{c}_{P,t}^* = 
\begin{cases}
\mb{c}_P ~~~~~~~~~~~~~~~~~~~~~~~~~~~~~~~~ 1 \leq t \leq {T_0}-1 \\
\mb{c}_P - \beta \mb{R}^\top_B \mb{r} + \beta \mb{R}^\top_B \mb{R}^\top_S \mb{c}_S ~~~~ t = {T_0}.
\end{cases}
\end{equation*}
As a result, an optimal PA policy is obtained by solving a single optimization problem
\begin{eqnarray}
\max_{\mbs{\theta}, \mbs{\pi}} && \min_{\mb{D}_{[1:{T_0}]} \in \mU^{T_0}} \max_{{\mb{S}_{[1:{T_0}]},\mb{X}_{[1:{T_0}]}}} ~~ \mb{P}_\infty^\text{PA} \Big( \mbs{\pi},\mb{D}_{[1:{T_0}]},\mb{X}_{[1:{T_0}]} \Big) \label{eq:obj_infty_PA} \\
&&
\begin{aligned}
\!\!\!\!\!\!\!\!\!\!\!\!\!\!\!\!\!\!\!\text{s.t.} ~~~~
& (\mb{S}_{[1:{T_0}]}, \mb{X}_{[1:{T_0}]}) \in \mathcal{X}(\mbs{\pi}, \mb{D}_{[1:{T_0}]},\mb{0},\mb{0}) \nonumber \\
& \mb{w}_1 ~\geq~ \mbs{\theta} \\
\end{aligned}
\end{eqnarray}
where the last constraint ensures that the solution is replicable over time periods.
Based on the solution of (\ref{eq:obj_infty_PA}), an infinite PA policy determines order quantity as
\begin{equation}
\label{policy:PA_infty}
\mb{s}_t = 
\begin{cases}
\mb{w}_1 - \mb{u}_{t-1}^\mb{s} + \mb{R}_S \mb{R}_B \mb{u}_{t-1}^\mb{d} \quad\quad t=n{T_0}+1, ~n=0,1,2,\ldots \\
\mb{w}_i + \displaystyle \sum_{\tau=1}^{i-1} \mb{W}_{\tau,i}\mb{d}_{n{T_0}+\tau} \quad\quad\quad\quad\;\; t=n{T_0} + i, 1 < i < T_0, ~n=0,1,2,\ldots.
\end{cases}
\end{equation}
We next present our main result for infinite-horizon cases.
\begin{theorem}
	\label{thm:optimality_PA_infty}
For a infinite-horizon multi-station network and the uncertainty set $\mU^\infty$, if Assumption~\ref{assumption:optimality_PA} holds, then the infinite periodic-affine policy (\ref{policy:PA_infty}) is optimal to the DP in~\eqref{eq:obj_PA_DP}.
\end{theorem}

In summary, we generalized the periodic-affine policies into multi-station networks and infinite horizon problems.
In both these cases, we presented periodic-affine algorithms and showed that the theoretical properties hold.
We next discuss a numerical case study to demonstrate the practical applicability of these findings and the performance of the proposed policies.

%%%%%%%%%%%%%%%%%

\section{Discussion: Insights and Implications}
\label{sec:experiments}
In this section, we present various implications of our modeling and solution approach and demonstrate the following advantages of our modeling approach and solution algorithm:
\begin{itemize}
\item \emph{Practical Relevance}: The relevance of an approach and corresponding algorithms hinges on the ability to model features in a real-world setting and provide implementable solutions. We demonstrate that our approach is able to achieve this. In particular, we show that PA policies perform well in large-scale and data-driven environments by studying the case of a major pharmacy retailer in India. We also demonstrate that we are able to model the service level guarantees by using the robust optimization approach. This ensures that the demands are satisfied for all scenarios in an uncertainty set. We also demonstrate the robustness to mis-specification and study the performance for a spectrum of various cost parameters.

\item \emph{Generalizability and Extendability}: It is also important for the approach to be generalizable and extendable in order to accommodate higher-dimensional versions of the problem and newer types of demand information. We achieve this by modifying the uncertainty set based on the available demand information, and by showing that our approach naturally extends to multi-dimensional settings. In particular, we incorporate correlation information in computing the optimal PA policy, and demonstrate our algorithm on the high-dimensional real-world case study.

\item \emph{Computational Tractability}: An algorithm suited for real applications needs to be tractable and implementable. We demonstrate tractability of the PA policies by presenting empirical evidence on the computational times on simulated data. % as well as on the data from the pharmacy retailer. 

\end{itemize}
Next we elaborate on each of these advantages.

\subsection{Practical Relevance: Case study of a Pharmacy retailer}
\label{subsec:pharma-single}
We analyze the sales data of a leading pharmacy retailer in India to probe the performance of the policies in a real-world setting.
A common problem in forecasting demands is that sales records do not necessarily imply customers' demands, because product shortage is not reflected in sales data.
However, since pharmacy retailers in India face a prohibitively large penalty for unmet demand, we can interpret the sales records as demands for this numerical study.

The data consists of more than 1.5 million transactions over 40 days (end of September to early November of 2016) for 228 different products.
To reduce the problem size, we analyzed the 20 most-popular products, comprising nearly 80\% of all transactions.
Hierarchical clustering \citep{package_cluster} is used to bundle the products into groups, within which demands are highly correlated.
Moving averages and residuals are extracted from the sales records, and used as nominal means and variances to define data-driven uncertainty sets for each group.
Penalty and holding costs are not revealed in the sales records. 
Therefore, we fix penalty and holding cost rates, and compute penalty and holding costs as a product of the corresponding rates and net profit per unit.

\subsubsection*{Uncertainty sets.}
We defined the multi-period uncertainty sets following Definition \ref{assumption:multi_period}, where the variability parameters are defined as $\Gamma = 2\sqrt{nT}$, $\Gamma_t = 2\sqrt{n}$, and $\Gamma_B = 2$.
The factor 2 is inspired from the fact that $\mathbb{P}(|Z| \leq 2) \simeq 0.95$ for a standard Normal random variable $Z$. 
Such a setting is typical in the robust optimization literature for describing random variables bounded by uncertainty sets \citep{Bertsimas2017}.
The parameters $\Gamma_t$ and $\Gamma$ are chosen to be proportional to $\sqrt{n}$ and $\sqrt{nT}$, motivated by the Central Limit Theorem \citep{bertsimas2004price}.

\subsubsection*{Performance of PA policies.}
We compute ordering policies for each product group. We begin by assuming that the sales of different product groups are independent of each other, and later consider the more realistic case of correlated sales. 
We compare the performance of three policies: PA, Aff-approx, and base-stock policies for these product groups.
We implemented two different base-stock policies.
Myopic base-stock policies are implemented, where the order-up-to level at each period is determined myopically using the nominal means and variances, assuming normally distributed demand.
We also compute the base-stock levels with Sample Average Approximation (SAA), using the approach discussed in \cite{Bertsimas2017}.

For pre-computed policies, we generated random demand samples to evaluate the two policies.
The samples are generated independently over time with Normal distributions, in which nominal mean and variances are used.
For given cost and revenue parameters, profit is calculated for each sample.
We compared 5\%, 25\% quantiles, and median for the policies.

\begin{table}[t]
	\centering
	\resizebox{\textwidth}{!}{
		\begin{tabular}{c@{\hskip 0.2in}c@{\hskip 0.45in} c@{\hskip 0.4in} c@{\hskip 0.3in} c@{\hskip 0.3in} c@{\hskip 0.4in} c@{\hskip 0.3in}}
			\toprule
			Penalty cost rate & Policy & Worst & 5\% quantile & 25\% quantile & Median & Historical \\
			\midrule
			\multirow{4}{*}{0.2} & Aff-approx & 2.68 & 3.50 & 3.64 & 3.74 & 3.81 \\
			& PA & 2.72 & 3.64 & 3.77 & 3.85 & 3.97 \\
			& SAA & N/A & {3.73} & {4.02} & {4.13} & {4.41} \\
			& {Myopic} & N/A & {3.80} & {4.16} & {4.37} & {4.67} \\
			\midrule
			\multirow{4}{*}{1.0} & Aff-approx & 1.96 & 2.87 & 3.27 & 3.51 & 3.86 \\
			& PA & 2.49 & 3.30 & 3.41 & 3.48 & 3.56 \\
			& SAA & N/A & {3.43} & {3.53} & {3.81} & {3.96} \\
			& {Myopic} & N/A & {3.27} & {3.66} & {3.96} & {4.25} \\
			\midrule
			\multirow{4}{*}{5.0} & Aff-approx & 1.49 & 3.31 & 3.60 & 3.79 & 4.04 \\
			& PA & 1.78 & 3.37 & 3.64 & 3.82 & 3.83 \\
			& SAA & N/A & {3.21} & {3.58} & {3.61} & {3.52} \\
			& {Myopic} & N/A & {2.85} & {3.27} & {3.58} & {3.68} \\
			\midrule
			\multirow{4}{*}{10.0} & Aff-approx & 1.26 & 3.10 & 3.45 & 3.68 & 4.15 \\
			& PA & 1.39 & 3.28 & 3.62 & 3.84 & 3.77 \\
			& SAA & N/A & {3.12} & {3.39} & {3.64} & {3.61} \\
			& {Myopic} & N/A & {2.76} & {3.19} & {3.49} & {3.40} \\
			\bottomrule
			\multicolumn{7}{l}{\footnotesize Note: all values are multiplied by $10^6$.} \\
			\bottomrule
	\end{tabular}}
	%	\vspace{-0.2in}
	\caption{Performance of policies for different penalty cost rates.
		Percentiles and medians are calculated from 1000 samples. The last column is from the historical sales data.}
	\label{tbl:single_pharma_result}
\end{table}

\autoref{tbl:single_pharma_result} displays performance of various policies for different values of the penalty cost rate.
We observe that PA performs better than Aff-approx in terms of worst-case performance.
In Section~\ref{subsec:random}, we will show that this is also consistent with synthetic data.
On the other hand, if the penalty cost is low, both PA and Aff-approx are not effective and the base-stock policy outperforms them.
However, PA yields better lower percentile performance than the other policies for increased penalty cost.
This is because PA maximizes the worst-case profit.
We also observe that under significant penalty costs, PA not only protects the worst-case performance and lower percentiles (improves by 19\% over base-stock at 5th percentile for cost of 10) but also leads to better average profit and historical backtesting than the other policies.
We also notice that while SAA improves over the Myopic policy, both Aff-approx and PA outperform it.

\subsubsection*{Sensitivity to model misspecification.}
\begin{figure}[h!]
	\centering
	\includegraphics[width=55mm]{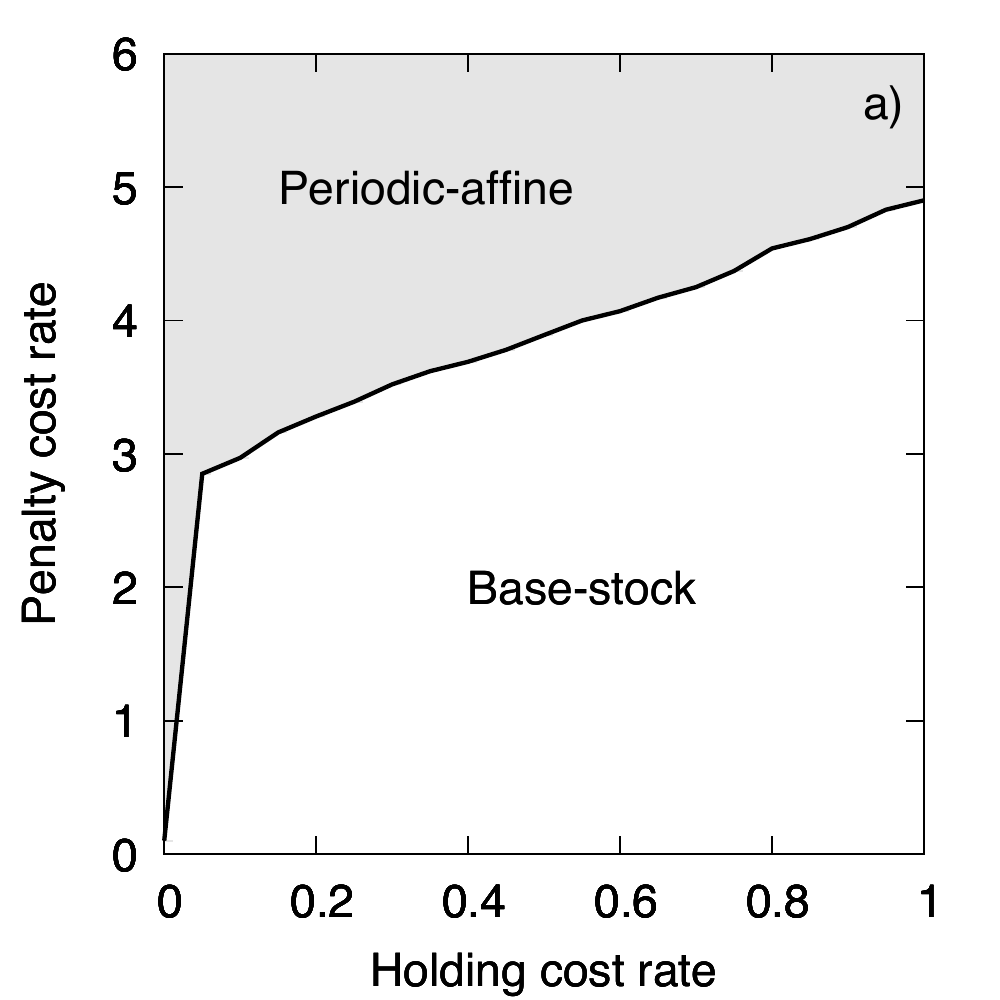}
	\hspace{-3mm}
	\includegraphics[width=55mm]{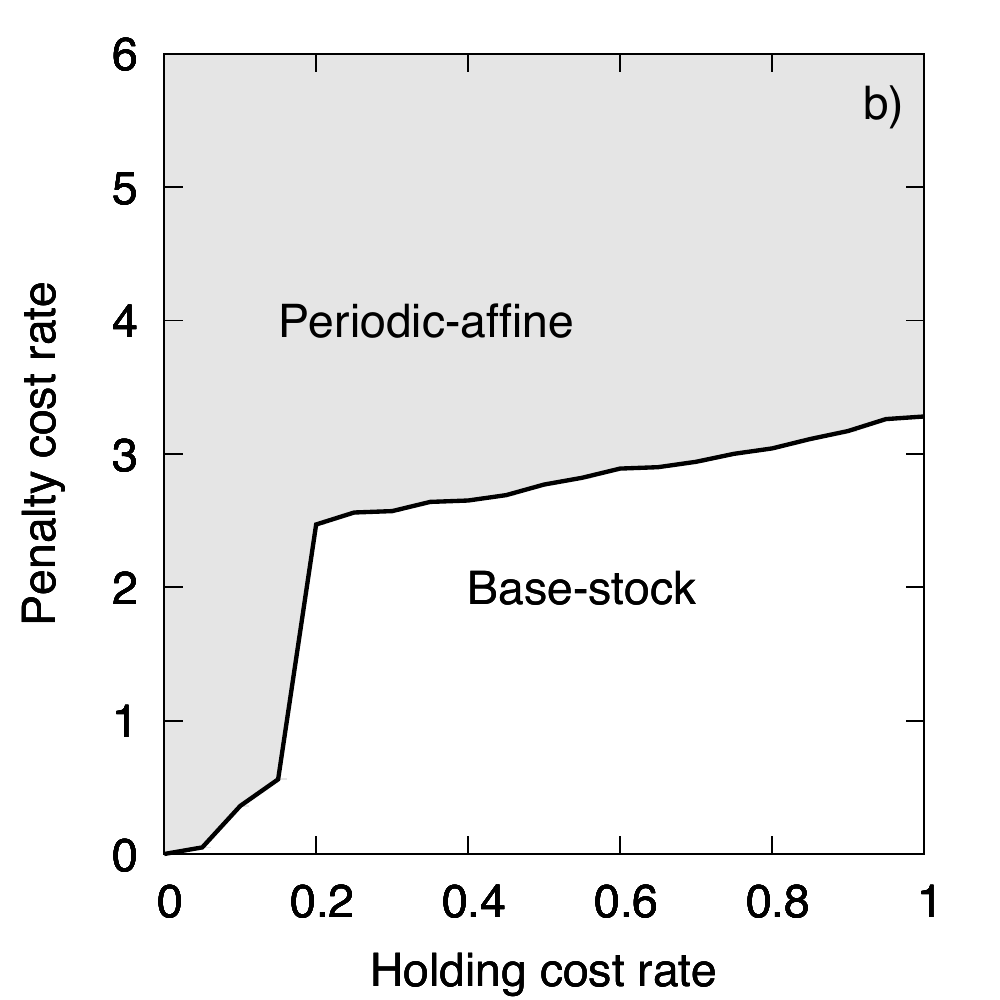}
	\hspace{-3mm}
	\includegraphics[width=55mm]{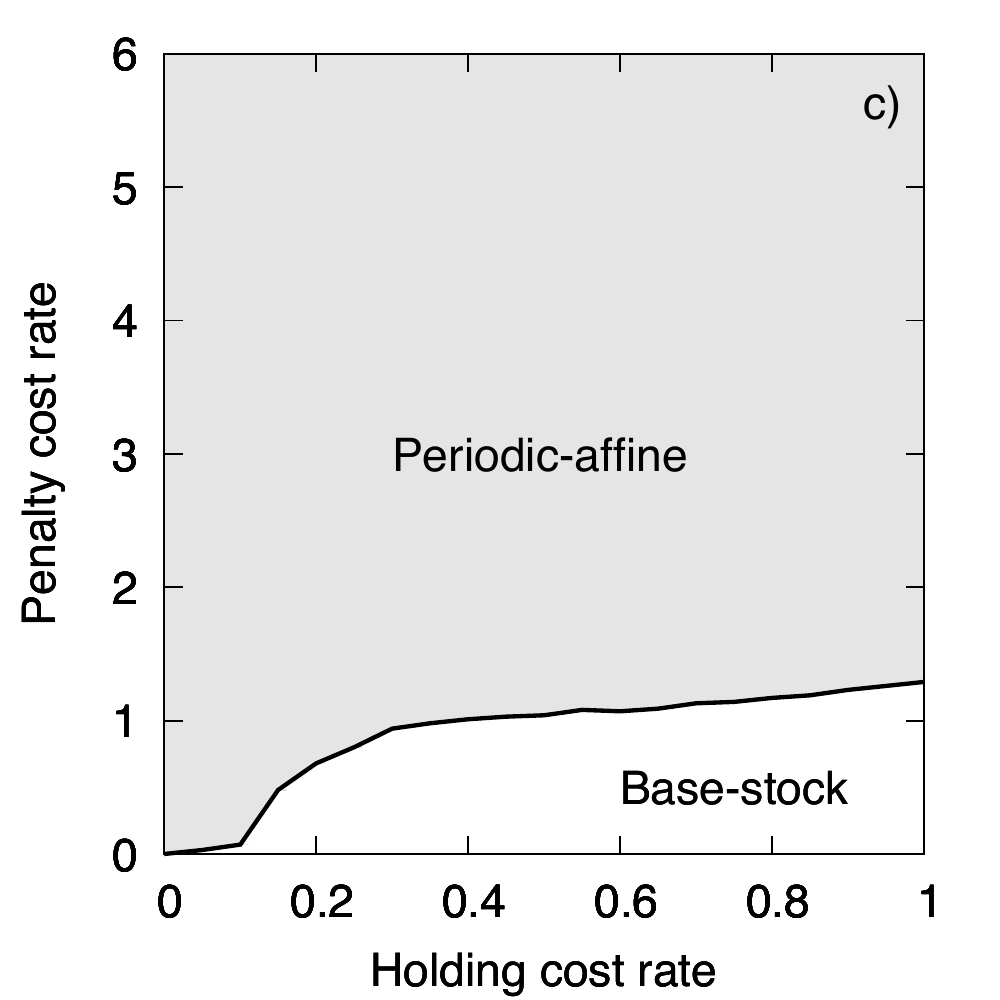}
	\caption{Phase diagram of periodic-affine and the base-stock policies. The realized demand means are (a) increased by 5\%, (b) not changed, and (c) decreased by 5\% from the nominal values.}
	\label{fig:region}
\end{figure}

Given that all these policies are implemented using the nominal mean and variance inferred from past records, it is important to measure their robustness to errors in model calibration.
For this, we consider demand realizations to have mean greater than (\autoref{fig:region}a), same as (\autoref{fig:region}b), or less than (\autoref{fig:region}c) their nominal values for varying holding and penalty costs.
We observe that when the realized demand distribution differs from the assumed one, the region of parameters (\emph{phase}), for which PA outperforms base-stock policy changes.
For example for a holding cost of 0.1 and penalty of 1, the PA policy outperforms base-stock, if the mean of the assumed demand coincides with the realized one (see \autoref{fig:region}b).
However, only 5\% increase of the means is sufficient for the base-stock to prevail (see \autoref{fig:region}a).

\subsubsection*{Performance dependence on holding and penalty cost.}
As the costs vary, we observe a \emph{phase transition} between a phase where PA outperforms the base-stock policies and a phase in which the base-stock policy outperforms.
The phase diagrams in \autoref{fig:region} allow the decision maker to select the policy based on the given cost and demand structure.
In fact, for pharmacy retailers, who face a substantial penalty with unsatisfied demand, it shows that our proposed PA policy is preferable.
On the other hand, if the decision maker is committed to a certain ordering policy (e.g. contractually), the phase diagrams in \autoref{fig:region} can suggest suitable changes to the cost structure in order to make the policy superior.

\subsubsection*{Impact of high penalty cost.}
\begin{figure}[h!]
	\centering
	\includegraphics[width=80mm]{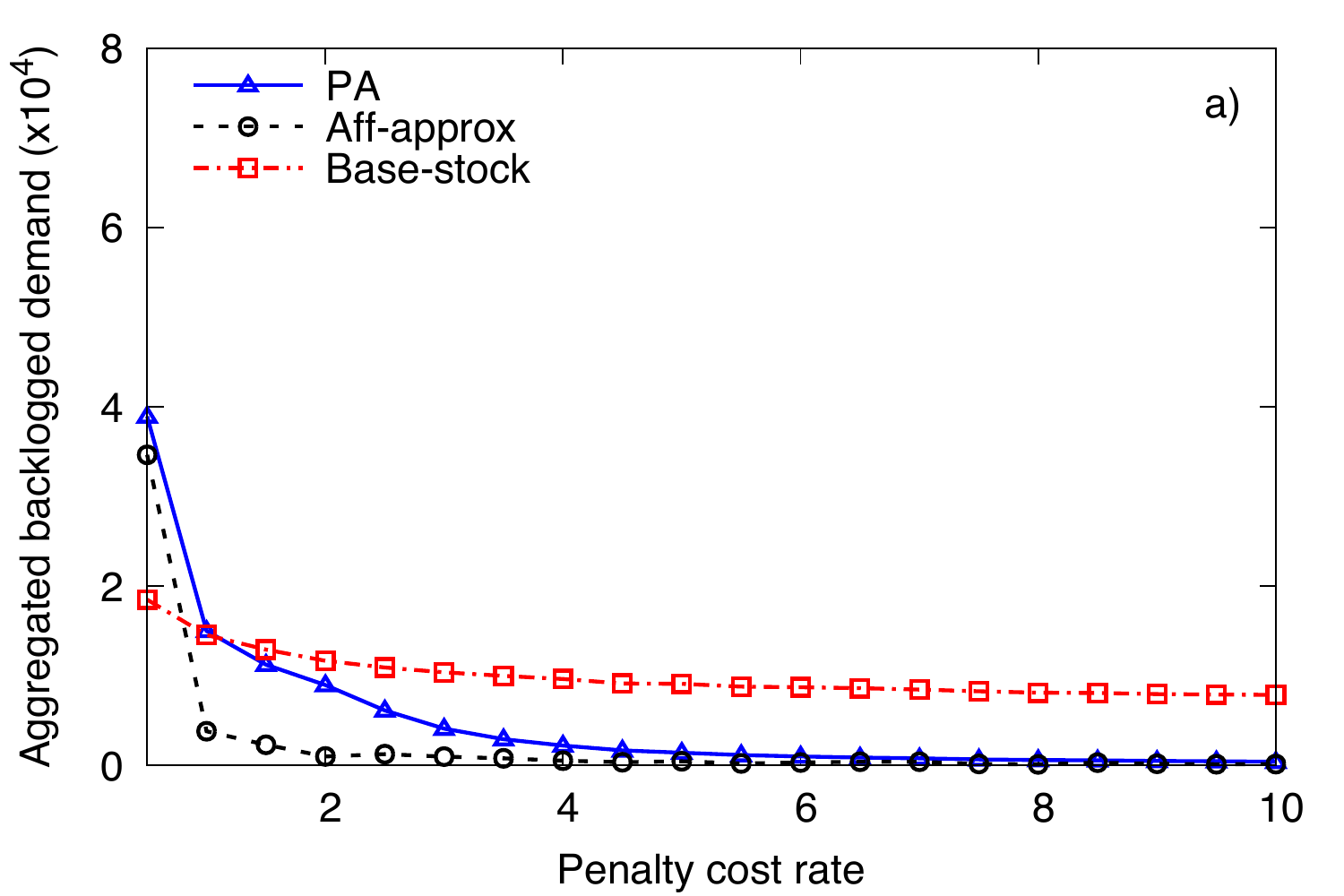}
	\includegraphics[width=80mm]{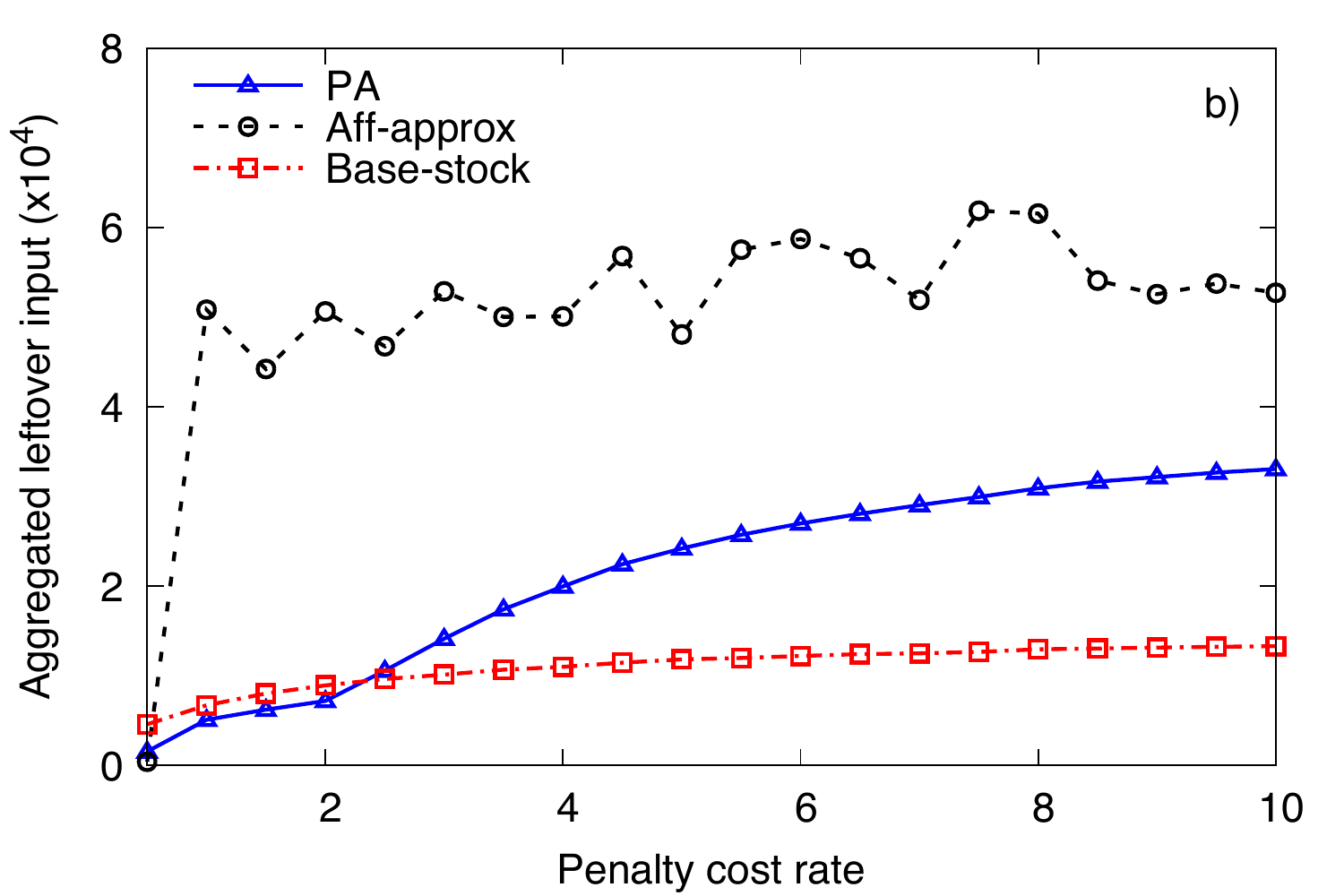}
	\caption{Impact of penalty cost on (a) backlogged demands, and (b) leftover resources on random samples.}
	\label{fig:backlog_history}
\end{figure}

When analyzing backlogged demands and inputs for different values of the penalty cost, \autoref{fig:backlog_history} shows that the three policies react differently for high penalty costs.
First, the base-stock policy does not effectively control the backlogged demands.
Although the penalty cost is accounted for in the newsvendor quantile to avoid high backlogs, increasing it slightly decreases the amount of backlogged demands.
On the other hand, Aff-approx determines order quantities more conservatively than PA.
Under high penalty cost, Aff-approx satisfies nearly all customer demands by ordering an excessive amount of input.
This causes a significantly larger holding cost, and leads to less profit than PA.
Finally, PA controls both leftover input and backlogged demands.
As the penalty cost increases, PA not only reduces backlogged demands (same as Aff-approx), but also maintains much lower input levels than Aff-approx. 
This leads to a higher profit than the other policies.

\subsubsection*{Summary.}
In summary, the case-study with pharmacy retailer's data demonstrates sizable increase in performance for the periodic-affine policies.
In particular, in the presence of high penalty cost, PA improves the lower-quantile performance by more than 10\% than the base-stock policies.
It also allows decision makers to identify the optimal policy based on their respective cost and demand structures.

\subsection{Generalizability and Extentability: Modeling Correlation and Solving multi-station newsvendor problems}
To adequately discuss the performance of PA under a multi-station setting, we take demand correlation information into account.
 We model correlations in demand using the \emph{correlated uncertainty sets} presented in (\ref{def:unc_set_corr}) by using a factor model approach.
For the numerical analysis, we consider two products over 15 time periods with subperiods of length 5.

\begin{figure}[b]
	\centering
	\includegraphics[width = \textwidth]{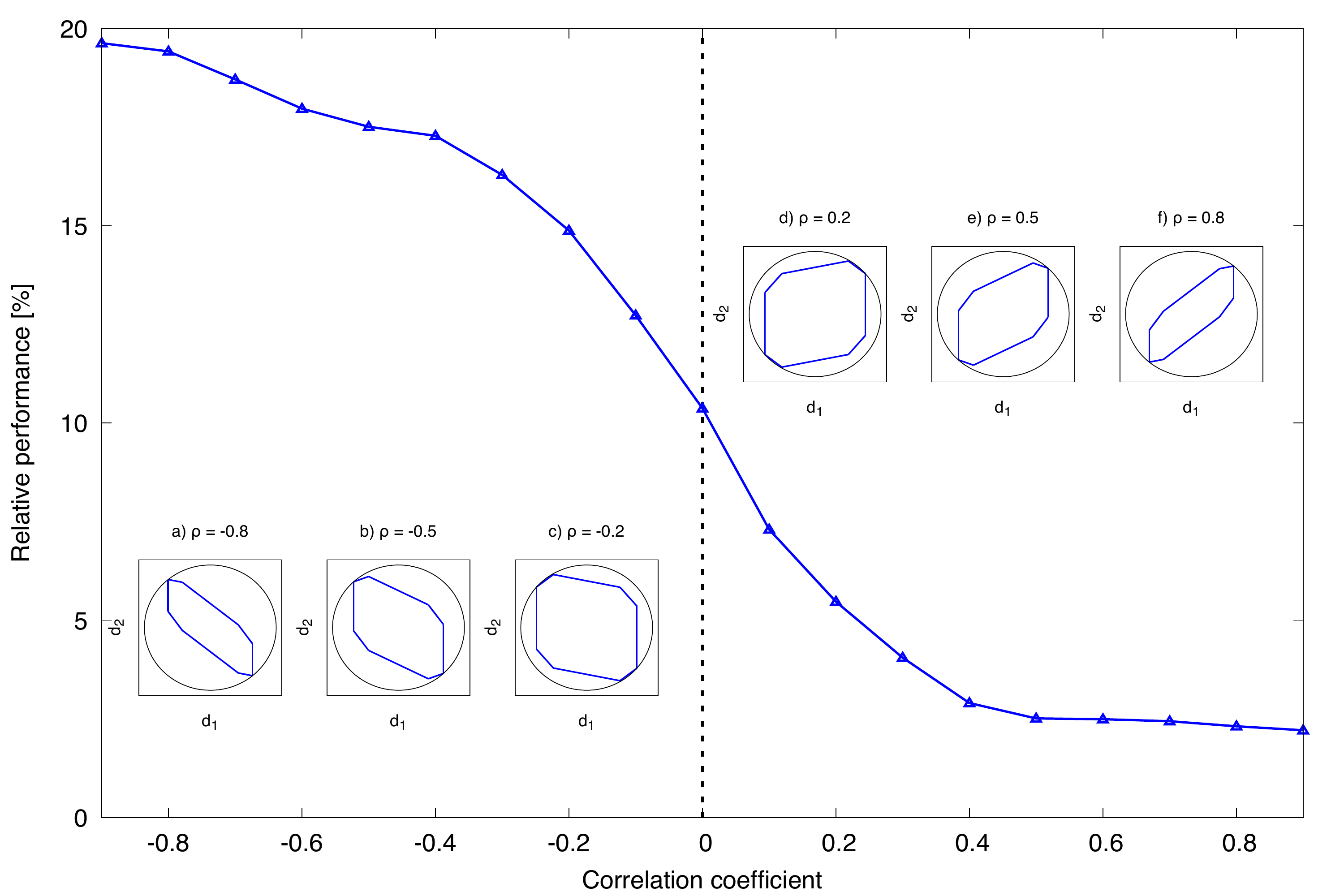}
	\vspace{-7mm}
	\caption{Relative performance for different correlations $\rho$. Inserts are the corresponding uncertainty sets.}
	\label{fig:impact_corr_rho}
\end{figure}
The benefits of modeling correlation become apparent, when comparing the following two policies: single-station PA (PA-single) for each product using marginal mean and variance, and multi-station PA (PA-multi) using the correlated uncertainty set.
For comparison, we compute the \emph{relative performance} (RP) of PA-multi over PA-single as
\begin{equation*}
	\text{RP}~ = ~ \frac{ \text{profit of PA-multi} - \text{profit of PA-single}}{\text{profit of PA-single}}.
\end{equation*}
After the two policies are implemented, we generate random demand with nominal mean and covariance and evaluate the relative performance for each sample.
For each correlation $\rho$, correlated uncertainty sets are defined by substituting
\begin{equation*}
	A = \begin{bmatrix}
	\rho & \sqrt{1-\rho^2} \\
	0 & 1
	\end{bmatrix}
\end{equation*}
into \eqref{def:unc_set_corr}, and the same matrix is used to generate multivariate normal random demands.
\autoref{fig:impact_corr_rho} displays this relative performance for different correlation coefficient $\rho$.
We observe that the median RP for every $\rho$ is positive.
However, the behavior differs for positively or negatively correlated products.
While for highly correlated products the RP slightly decreases with growing $\rho$, significant improvements are made for negatively correlated products.
In fact, more RP increases by more than 17\%, when the products have a strong negative correlation ($\rho = -0.9$).

This observation can be interpreted by the structure of the uncertainty sets, as shown in the inserts of \autoref{fig:impact_corr_rho}.
Positively correlated products lead to sets that allow both uncertain demands $d_1$ and $d_2$ to be concurrently at their maximum or minimum value.
As $\rho$ decreases, the area of the polyhedron shrinks, however the extreme points are unaffected.
However, when $\rho$ becomes negative, i.e. the products are negatively correlated, if one of the uncertain demands can take its maximum value, the other is forced to its lowest, and vice versa.
This effect forces an increase in RP as $\rho\rightarrow -1$.
\begin{figure}[t]
	\centering
	\includegraphics[width=\textwidth/3]{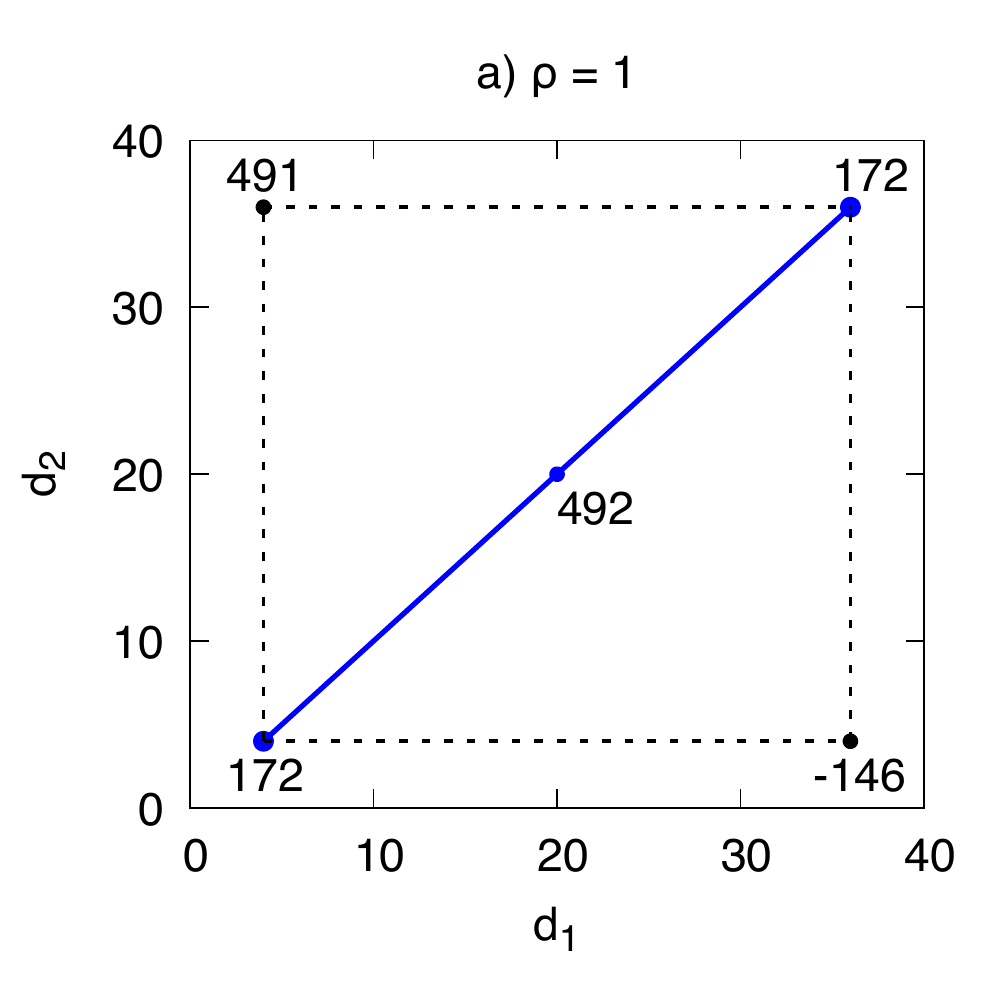}
	\hspace{-3mm}
	\includegraphics[width=\textwidth/3]{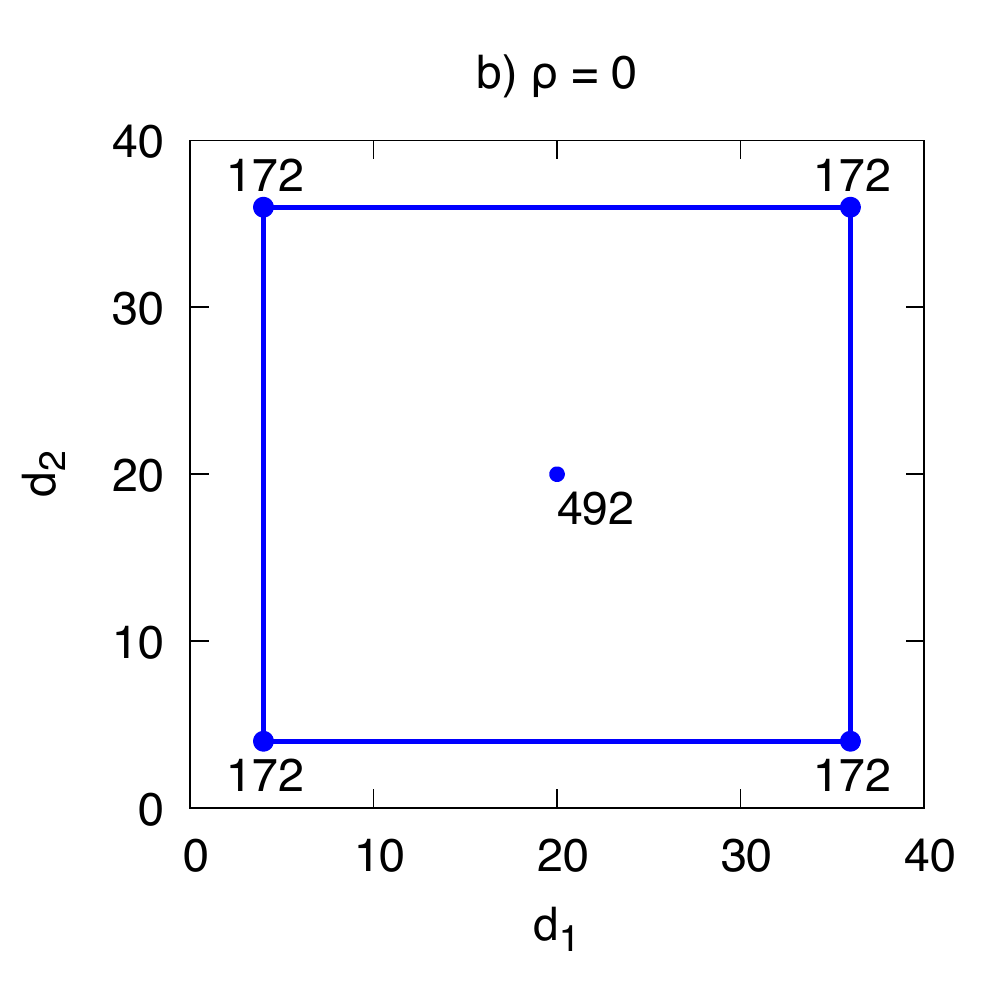}
	\hspace{-3mm}
	\includegraphics[width=\textwidth/3]{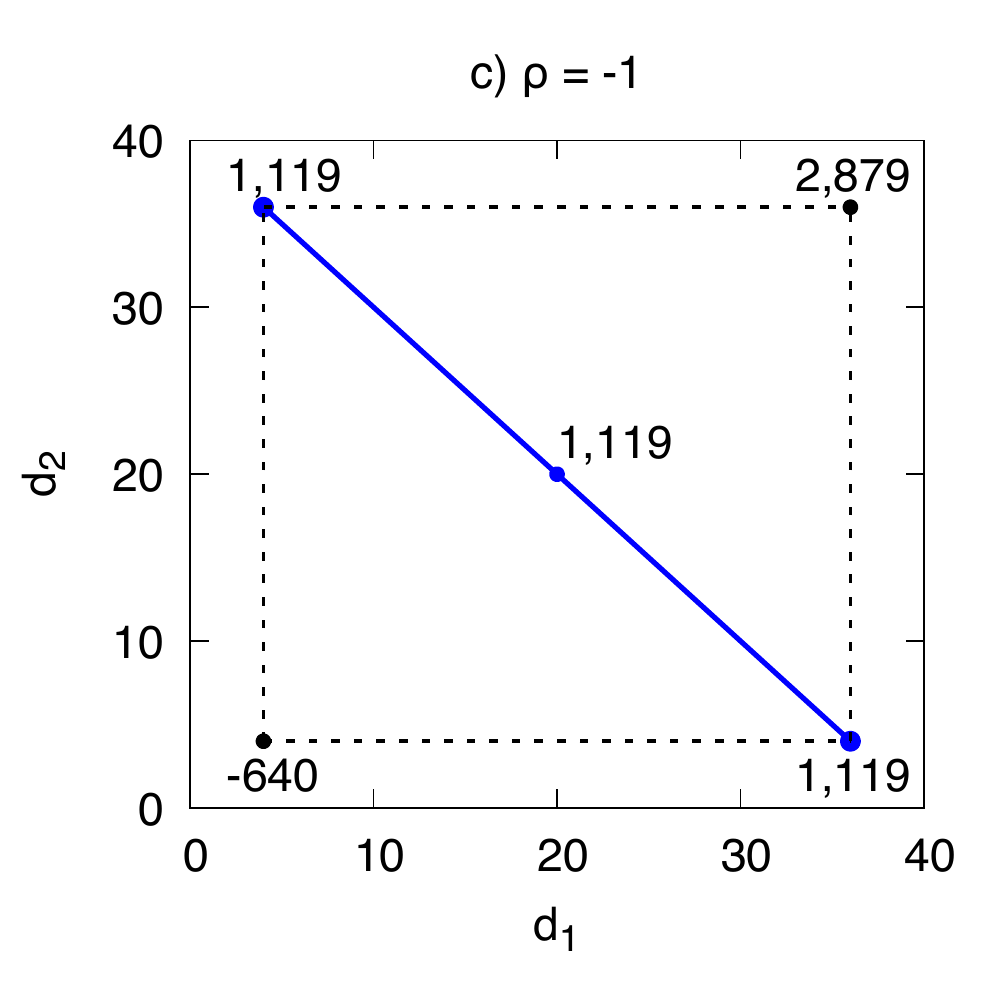}
	\vspace{-3mm}
	\caption{An illustrative example for correlation and multi-station modeling:
		Two products are (a) perfectly correlated, (b) uncorrelated, and (c) perfectly anti-correlated.
		Uncertainty sets are defined as within the blue contours along with the profit of optimal policies for extreme and nominal points.}
	\label{fig:example_corr}
\end{figure}
The extreme cases are illustrated in the example of \autoref{fig:example_corr}.
For perfectly correlated products, \autoref{fig:example_corr}a shows that even though the uncertainty set is dramatically shrunk, the worst-case profit cannot not improved over the uncorrelated demand setting (\autoref{fig:example_corr}b), because the worst-case is often captured when both demands are high or low.
However, for negatively correlated demands, the uncertainty set does not contain this region (high/high or low/low), allowing for substantial improvement in worst-case profit, as shown in \autoref{fig:example_corr}c.

\subsubsection*{Case study of a pharmacy retailer -- revisited.}
Here, we account for correlation amongst the product groups and compare the following five policies in \autoref{tbl:multi_pharma_result}: PA-multi, PA-single, Aff-approx with correlated uncertainty set (Aff-approx-multi), Aff-approx for each group (Aff-approx-single), SAA base-stock, and Myopic policies.
\begin{table}[h!]
	\centering
	\resizebox{\textwidth}{!}{
		\begin{tabular}{c@{\hskip 0.2in}c@{\hskip 0.4in} c@{\hskip 0.4in} c@{\hskip 0.3in} c@{\hskip 0.3in} c@{\hskip 0.4in} c@{\hskip 0.1in}}
			\toprule
			Penalty cost rate& Policy & Worst & 5\% quantile & 25\% quantile & Median & Historical \\
			\midrule
			\multirow{6}{*}{0.2} & Aff-approx-single & 6.20 & 9.09 & 9.65 & 10.03 & 10.60 \\
			& Aff-approx-multi & 7.03 & 8.65 & 9.72 & 10.31 & 9.51 \\
			& PA-single & 6.32 & 8.59 & 8.97 & 9.21 & 9.66 \\
			& PA-multi & 7.17 & 9.24 & 10.24 & 10.71 & 11.12 \\
				& SAA & N/A & {9.12} & {10.04} & {10.51} & {10.93} \\
			& {Myopic} & N/A & {8.50} & {9.97} & {10.74} & {10.81} \\
			\midrule
			\multirow{6}{*}{1.0} & Aff-approx-single & 4.07 & 6.44 & 7.59 & 8.15 & 8.39 \\
			& Aff-approx-multi & 5.29 & 8.22 & 8.83 & 9.34 & 10.02 \\
			& PA-single & 5.69 & 7.76 & 8.07 & 8.31 & 8.97 \\
			& PA-multi & 6.64 & 8.65 & 9.30 & 9.70 & 10.54 \\
			& SAA & N/A & {7.49} & {7.82} & {9.35} & {9.13} \\
			& {Myopic} & N/A & {6.70} & {8.10} & {9.21} & {8.98} \\
			\midrule
			\multirow{6}{*}{5.0} & Aff-approx-single & 2.63 & 5.19 & 6.82 & 7.90 & 7.21 \\
			& Aff-approx-multi & 4.20 & 7.44 & 8.44 & 9.01 & 9.37 \\
			& PA-single & 3.44 & 7.57 & 8.33 & 8.77 & 9.75 \\
			& PA-multi & 4.95 & 7.82 & 8.73 & 9.43 & 9.25 \\
			& SAA & N/A & {6.41} & {7.23} & {8.61} & {8.92} \\
			& {Myopic} & N/A & {5.04} & {6.61} & {7.76} & {7.03} \\
			\midrule
			\multirow{6}{*}{10.0} & Aff-approx-single & 1.69 & 5.06 & 6.67 & 7.64 & 6.88 \\
			& Aff-approx-multi & 3.46 & 7.07 & 8.17 & 8.85 & 9.88 \\
			& PA-single & 2.11 & 7.26 & 8.24 & 8.81 & 9.01 \\
			& PA-multi & 3.98 & 7.04 & 8.25 & 8.98 & 8.87 \\
			& SAA & N/A & {6.01} & {7.12} & {7.83}& {7.61} \\
			& {Myopic} & N/A & {4.61} & {6.27} & {7.40} & {6.42} \\
			\midrule
			\multicolumn{7}{l}{\footnotesize Note: all values are multiplied by $10^6$.} \\
			\bottomrule
	\end{tabular}}
	\caption{Performance of policies for two correlated product groups for different penalty cost rates.
		Percentiles and medians are calculated from 1000 samples. The last column is from the historical sales data.}
	\label{tbl:multi_pharma_result}
\end{table}
For lower quantiles, we observe that the base-stock policy performs poorly compared to the single-station models.
Both PA-single and PA-multi yield significantly greater profit in lower quantiles than the base-stock policy.
The multi-station framework and correlated demand uncertainty sets offer better performance than single-station framework.
In particular, PA-multi achieves at least 7\% more profit than PA-single for moderate choice of the penalty cost. 
However, for extremely high penalty cost rates (e.g. 10), PA-single performs slightly better than PA-multi for lower quantiles, even though the worst case objective value is lower.
This is due to samples that are generated outside of the correlated uncertainty sets. 
Note that SAA and the Myopic policies do not take correlations into account and hence underperform.

\subsubsection*{Summary.}
In summary, \autoref{tbl:multi_pharma_result} demonstrates that capturing correlation using the periodic-affine policies in multi-station setting outperforms the base-stock policies for moderate penalties.
Since the demand of the studied pharmaceutical product groups is positively correlated, variability of the base-stock policy increases, causing a sizable degradation of the profit when compared to periodic-affine policies.

\subsection{Computational Tractability}
\label{subsec:random}
In order to focus on computational performance in a sterile environment, we consider the following simulation environment. We simulate three cases with duration $T\in\{10,15,20\}$ with a subperiod consisting of 5 time periods.
We randomly generate 100 instances of single-station newsvendor problems for each $T$. 
Nominal means are generated by autoregression processes AR(1) and nominal coefficient of variations are chosen uniformly in (0.3, 0.5).
Unless modified, cost parameters are $c_S = 20$, $r = 120$, $c_H = c_P = 20$, and all variability parameters are set similarly as in Section \ref{subsec:pharma-single}. We then compare the PA policies with affine and Aff-approx policies.

While we evaluate the worst-case performances of PA using \eqref{eq:obj_PA_DP} with Proposition \ref{prop:PA_worst_scenario}, those of affine policies are calculated using (\ref{eq:obj_multi_FA}--\ref{eq:constr_multi_FA_second}) to overcome tractability issues.
Note that for any non-anticipative policies, (\ref{eq:obj_multi_FA}--\ref{eq:constr_multi_FA_second}) is an upper bound for the DP \eqref{eq:obj_PA_DP}.
Hence, the optimality gap between the affine policies and PA, evaluated in \eqref{eq:obj_PA_DP}, is closer than the gaps presented in Table~\ref{tbl:tractability} and Figure~\ref{fig:PA_worst}.

For the three policies, \autoref{tbl:tractability} displays the computation times and worst objective values on the same uncertainty set.
We observe that the computation of PA is significantly faster than affine policies, because PA is tractable.
However, the worst objective values of PA are very close (within 0.1\%) to affine policies, while Aff-approx consistently deviates by $\ge10\%$ from the others.
Indeed, only 13 out of the 300 artificial instances have greater worst objective values in affine policies. 
Moreover, there is only one instance in which PA loses more than 1\% of optimality.
This implies that Assumption~\ref{assumption:optimality_PA}, which is a sufficient optimality condition of PA, holds for fairly general settings.
This means that PA is as competitive as affine policies in worst-case values.

\begin{table}[h!]
	\centering
	\resizebox{\textwidth}{!}{
	\begin{tabular}{ c@{\hskip 0.4in}c@{\hskip 0.1in} c@{\hskip 0.1in} c@{\hskip 0.1in} c@{\hskip 0.2in} c@{\hskip 0.1in} c@{\hskip 0.1in} c@{\hskip 0.1in} c@{\hskip 0.2in} c@{\hskip 0.1in} c@{\hskip 0.1in} c@{\hskip 0.1in}}
		\toprule
		Policy & \multicolumn{3}{c}{Affine} & & \multicolumn{3}{c}{PA} & & \multicolumn{3}{c}{Aff-approx} \\
		\midrule
		Time periods $T$& $10$ & $15$ & $20$ & & $10$ & $15$ & $20$ & & $10$ & $15$ & $20$ \\
		\midrule
		Computation time [sec] & 3.8 & 88.5 & 1388.3 & & 0.18 & 0.24 & 0.29 & & 0.02 & 0.09 & 0.22 \\
		Worst objective value & 11,035 & 16,833 & 22,440 & & 11,033 & 16,831 & 22,438 & & 9,523 & 14,601 & 19,479 \\
		Difference to Affine [\%] &  &  &  & & -0.015 & -0.011 & -0.006 & & -13.39 & -13.03  & -12.96  \\
		\bottomrule
	\end{tabular}}
	\vspace{0.1in}
	\caption{Average performance of three policies on randomly generated instances.}
	\label{tbl:tractability}
\end{table}

\subsubsection*{Comparison of PA and Aff-approx policies on synthetic data.}
\begin{figure}[h!]
	\centering
	\includegraphics[width = \textwidth/2]{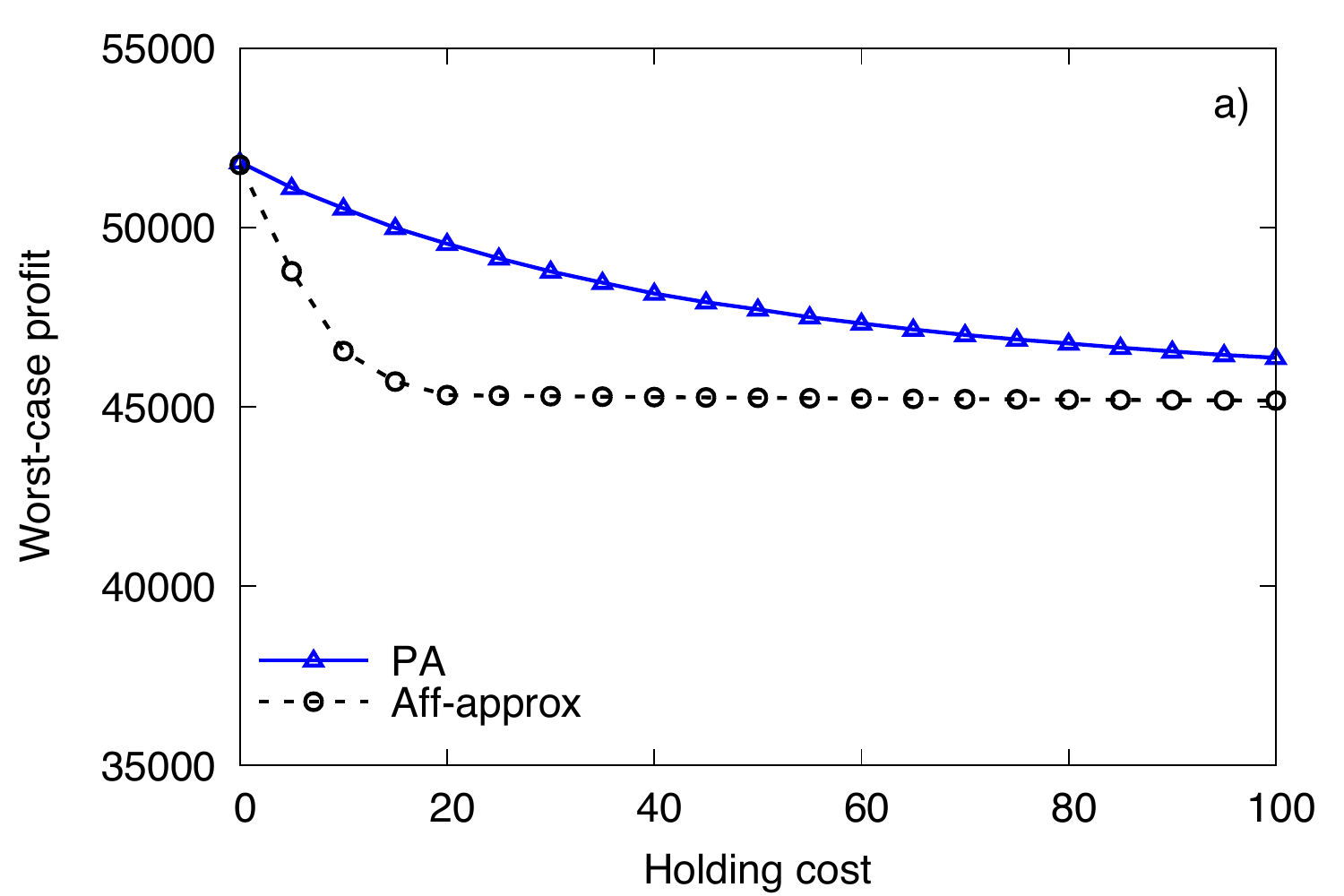}
	\hspace{-3mm}
	\includegraphics[width = \textwidth/2]{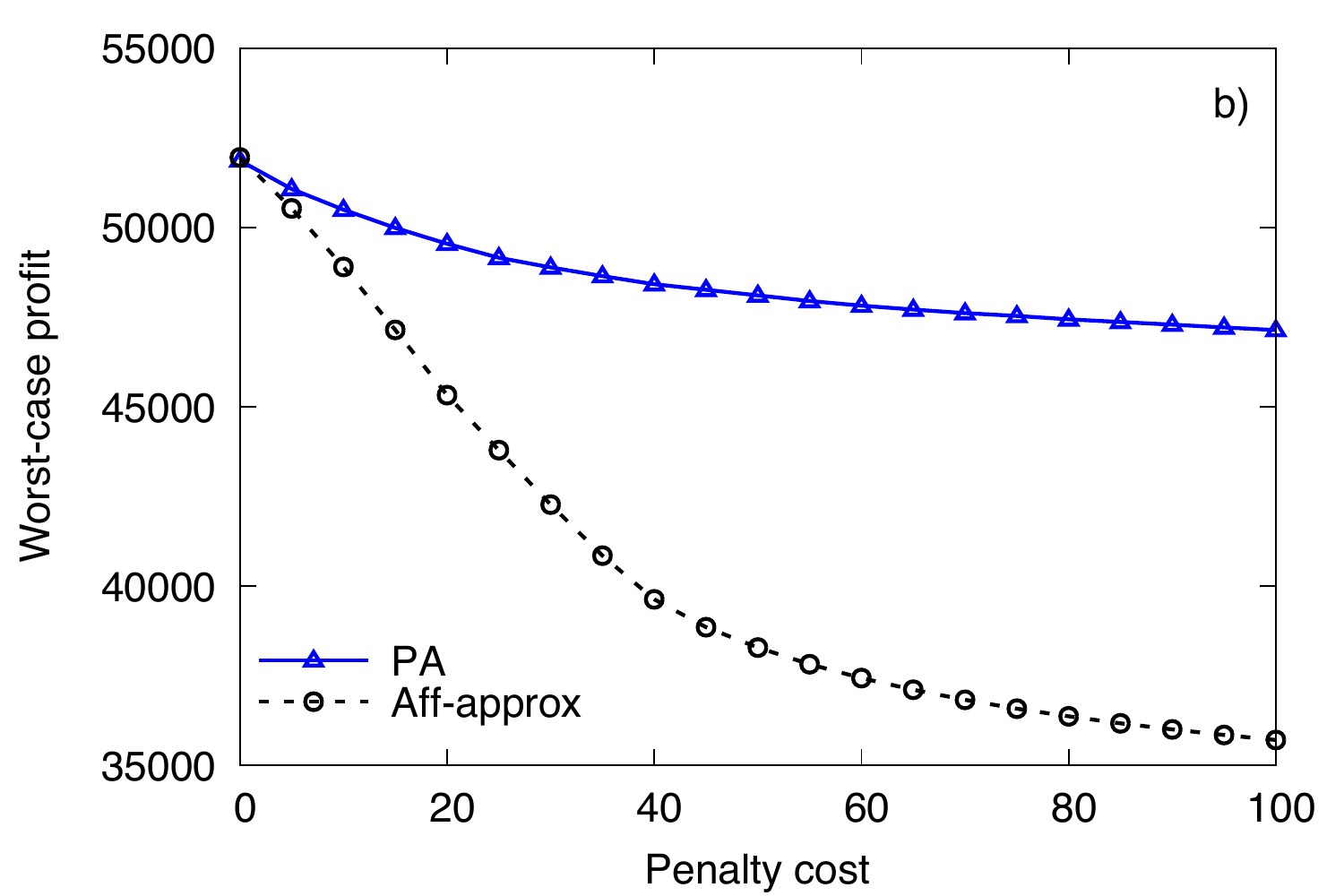}
	\vspace{-7mm}
	\caption{Impact of (a) holding cost, and (b) penalty cost on worst-case performance.}
	\label{fig:PA_worst}
\end{figure}

We next compare the performance of PA and Aff-approx policies for a spectrum of parameters.
\autoref{fig:PA_worst} shows that the gap between the two policies are different for holding and for penalty cost.
\autoref{fig:PA_worst}a shows that the gap between PA and Aff-approx decreases as holding cost increases.
However, PA protects the worst-case profit significantly better than Aff-approx as penalty cost increases, shown in \autoref{fig:PA_worst}b.
These results suggest that depending of the holding and penalty costs, PA orders the proper amount of input to meet demands, i.e., both on-hand input and unsatisfied demands are well-controlled.
Therefore, its worst-case performance is far less affected by larger cost parameters than Aff-approx.
On the other hand, \autoref{fig:PA_worst}b indicates that Aff-approx may not manage the leftover resources and backlogged demands as well as PA, which yields poor performance for higher penalty cost.

The significant relative decrease of worst objective values between the two policies of up to 35\%, as shown in~\autoref{fig:PA_worst} and also observed in~\autoref{tbl:tractability}, implies that the degree of suboptimality of Aff-approx may render it inferior to PA.
Although theoretical bounds of affine approximation have been proposed \citep{bertsimas2012power, bertsimas2015performance}, these results indicates that the affine approximation cannot be successfully applied to general settings, despite its computational advantage of tractability.

\subsubsection*{Summary.}
In summary, the experiment on artificial data reveals the optimality and tractability of periodic-affine policies, while the other methods did not achieve both properties.
In addition, we observe that periodic-affine policies perform substantially better for larger $T$ than Aff-approx in terms of protecting its worst-case performance at high penalty cost.
This implies that PA is particularly useful when a decision maker faces a massive penalty for unsatisfied demands, such as in the case of the pharmacy retailer in India.

%%%%%%%%%%%%%%%%%
\section{Conclusions}
\label{sec:concl}
In this paper, we consider the problem of optimal control in multi-period and multi-stage newsvendor networks. To this end, we introduce a new class of adaptive policies called \emph{periodic-affine policies}.
These policies are data-driven and incorporate the correlation amongst products, which is an instrumental feature of real-world settings.
These policies also remain robust to parameter mis-specification. 
For this, we model the uncertain demand via sets, which can incorporate correlations, and can be generalized to multi-product settings and extended to multi-period problems.
This approach offers a natural framework to study current competing policies of base-stock, affine, and approximative approaches with respect to their profit, sensitivity to parameters and assumptions, and computational scalability.
We showed that the periodic-affine policies are sustainable, i.e. time consistent, because they warrant optimality both within subperiods and over the entire planning horizon.

We presented empirical evidence of tractability and robustness which makes our approach well-suited for real-world applications.
We demonstrate the advantages of our approach by considering the problem of one of India's largest pharmacy retailers using their sales data. 
We show that the periodic-affine policies are capable to increase the profits by up to 17\% over base-stock policies.
This study reveals capturing the demand correlation can sizably affect the performance.
Furthermore, we offer a phase diagram for managers to select the optimal policy based on their cost and demand structures.

In future, we intend to incorporate time-dependent uncertainty sets~\citep{nohadani2017robust} to more accurately model seasonal demand.
This step forward will lend itself well to incorporate returns, i.e. feedback from satisfied demand that can guide the next period's decisions.

\ACKNOWLEDGMENT{We are grateful to Jan A. Van Mieghem for insightful comments.
We also thank two anonymous referees for their thorough and suggestive reviews.}

\FloatBarrier
%%%%%%%%%%%%%%%%%
\bibliographystyle{informs2014}
\bibliography{ref}

\begin{thebibliography}{49}
\providecommand{\natexlab}[1]{#1}
\providecommand{\url}[1]{\texttt{#1}}
\providecommand{\urlprefix}{URL }

\bibitem[{Arrow et~al.(1951)Arrow, Harris, \protect\BIBand{}
  Marschak}]{Arrow1951}
Arrow KJ, Harris T, Marschak J (1951) Optimal inventory policy.
  \emph{Econometrica} 250--272.

\bibitem[{Baboli et~al.(2011)Baboli, Fondrevelle, Tavakkoli-Moghaddam,
  \protect\BIBand{} Mehrabi}]{baboli2011replenishment}
Baboli A, Fondrevelle J, Tavakkoli-Moghaddam R, Mehrabi A (2011) A
  replenishment policy based on joint optimization in a downstream
  pharmaceutical supply chain: centralized vs. decentralized replenishment.
  \emph{The Internat. J. Adv. Manufacturing Tech.} 57(1-4):367--378.

\bibitem[{Bandi \protect\BIBand{} Bertsimas(2012)}]{bandi2012tractable}
Bandi C, Bertsimas D (2012) Tractable stochastic analysis in high dimensions
  via robust optimization. \emph{Math. Programming} 134(1):23--70.

\bibitem[{Ben-Tal et~al.(2013)Ben-Tal, Den~Hertog, De~Waegenaere, Melenberg,
  \protect\BIBand{} Rennen}]{ben2013robust}
Ben-Tal A, Den~Hertog D, De~Waegenaere A, Melenberg B, Rennen G (2013) Robust
  solutions of optimization problems affected by uncertain probabilities.
  \emph{Management Sci.} 59(2):341--357.

\bibitem[{Ben-Tal et~al.(2005)Ben-Tal, Golany, Nemirovski, \protect\BIBand{}
  Vial}]{ben2005retailer}
Ben-Tal A, Golany B, Nemirovski A, Vial J (2005) Retailer-supplier flexible
  commitments contracts: a robust optimization approach. \emph{Manufacturing \&
  Service Oper. Management} 7(3):248--271.

\bibitem[{Ben-Tal et~al.(2004)Ben-Tal, Goryashko, Guslitzer, \protect\BIBand{}
  Nemirovski}]{ben2004}
Ben-Tal A, Goryashko A, Guslitzer E, Nemirovski A (2004) Adjustable robust
  solutions of uncertain linear programs. \emph{Math. Programming}
  99(2):351--376.

\bibitem[{Bertsimas \protect\BIBand{}
  Bidkhori(2015)}]{bertsimas2015performance}
Bertsimas D, Bidkhori H (2015) On the performance of affine policies for
  two-stage adaptive optimization: a geometric perspective. \emph{Math.
  Programming} 153(2):577--594.

\bibitem[{Bertsimas \protect\BIBand{} Goyal(2012)}]{bertsimas2012power}
Bertsimas D, Goyal V (2012) On the power and limitations of affine policies in
  two-stage adaptive optimization. \emph{Math. Programming} 134(2):491--531.

\bibitem[{Bertsimas et~al.(2017)Bertsimas, Gupta, \protect\BIBand{}
  Kallus}]{Bertsimas2017}
Bertsimas D, Gupta V, Kallus N (2017) Robust sample average approximation.
  \emph{Math. Programming} ISSN 1436-4646.

\bibitem[{Bertsimas et~al.(2010)Bertsimas, Iancu, \protect\BIBand{}
  Parillo}]{Bertsimas2010}
Bertsimas D, Iancu D, Parillo P (2010) Optimality of affine policies in
  multistage robust optimization. \emph{Math. Oper. Res.} 35(2):363--394.

\bibitem[{Bertsimas \protect\BIBand{} Sim(2004)}]{bertsimas2004price}
Bertsimas D, Sim M (2004) The price of robustness. \emph{Oper. Res.}
  52(1):35--53.

\bibitem[{Bertsimas \protect\BIBand{} Thiele(2006)}]{Bertsimas2006}
Bertsimas D, Thiele A (2006) Robust and data-driven optimization: Modern
  decision making under uncertainty. \emph{Tutorials in Oper. Res.} 95--122.

\bibitem[{Bienstock \protect\BIBand{} \"Ozbay(2008)}]{Bienstock2008}
Bienstock D, \"Ozbay N (2008) Computing robust base-stock levels.
  \emph{Discrete Optim.} 5(2):389--414.

\bibitem[{Clark \protect\BIBand{} Scarf(1960)}]{Clark1960}
Clark AJ, Scarf H (1960) Optimal policies for a multi-echelon inventory
  problem. \emph{Management Sci.} 6(4):475--490.

\bibitem[{Crawford \protect\BIBand{} Shum(2005)}]{crawford2005uncertainty}
Crawford GS, Shum M (2005) Uncertainty and learning in pharmaceutical demand.
  \emph{Econometrica} 73(4):1137--1173.

\bibitem[{Delage \protect\BIBand{} Ye(2010)}]{delage2010distributionally}
Delage E, Ye Y (2010) Distributionally robust optimization under moment
  uncertainty with application to data-driven problems. \emph{Oper. Res.}
  58(3):595--612.

\bibitem[{Federgruen \protect\BIBand{} Zipkin(1986)}]{Federgruen1986}
Federgruen A, Zipkin P (1986) An inventory model with limited production
  capacity and uncertain demands i. the average-cost criterion. \emph{Math.
  Oper. Res.} 11(2):193--207.

\bibitem[{Fu et~al.(2005)Fu, Glover, \protect\BIBand{} April}]{Fu2005}
Fu M, Glover F, April J (2005) Simulation optimization: a review, new
  developments, and applications. \emph{Proc. 2005 Winter Simulation Conf.},
  83--95.

\bibitem[{Fu \protect\BIBand{} Healy(1997)}]{Fu1997b}
Fu M, Healy K (1997) Techniques for simulation optimization: An experimental
  study on an $({s},{S})$ inventory system. \emph{IIE Trans.} 29(3):191--199.

\bibitem[{Fu(1994)}]{Fu1994b}
Fu MC (1994) Sample path derivatives for (s, {S}) inventory systems.
  \emph{Oper. Res.} 42(2):351--364.

\bibitem[{Gallego \protect\BIBand{} Moon(1993)}]{Gallego1993}
Gallego G, Moon I (1993) The distribution free newsboy problem: review and
  extensions. \emph{J. Oper. Res. Soc.} 44(8):825--834.

\bibitem[{Glasserman \protect\BIBand{} Tayur(1995)}]{Glasserman1995}
Glasserman P, Tayur S (1995) Sensitivity analysis for base-stock levels in
  multiechelon production-inventory systems. \emph{Management Sci.}
  41(2):263--281.

\bibitem[{Graves \protect\BIBand{} Willems(2000)}]{Graves2000}
Graves S, Willems S (2000) Optimizing strategic safety-stock placement in
  supply chains. \emph{Manufacturing \& Service Oper. Management} 2(1):68--83.

\bibitem[{Guerrero et~al.(2013)Guerrero, Yeung, \protect\BIBand{}
  Gu{\'e}ret}]{guerrero2013joint}
Guerrero WJ, Yeung T, Gu{\'e}ret C (2013) Joint-optimization of inventory
  policies on a multi-product multi-echelon pharmaceutical system with batching
  and ordering constraints. \emph{Eur. J. Oper. Res.} 231(1):98--108.

\bibitem[{Huh \protect\BIBand{} Janakiraman(2008)}]{Huh2008}
Huh W, Janakiraman G (2008) A sample-path approach to the optimality of echelon
  order-up-to policies in serial inventory systems. \emph{Oper. Res. Lett.}
  36(5):547--550.

\bibitem[{Iyengar(2005)}]{iyengar2005robust}
Iyengar GN (2005) Robust dynamic programming. \emph{Math. Oper. Res.}
  30(2):257--280.

\bibitem[{Kapuscinski \protect\BIBand{} Tayur(1999)}]{Kapiscinski1999}
Kapuscinski R, Tayur S (1999) Optimal policies and simulation-based
  optimization for capacitated production inventory systems. \emph{Quantitative
  Models for Supply Chain Management}, 7--40 (Springer).

\bibitem[{Karlin(1960)}]{Karlin1960}
Karlin S (1960) Dynamic inventory policy with varying stochastic demands.
  \emph{Management Sci.} 6(3):231--258.

\bibitem[{Kasugai \protect\BIBand{} Kasegai(1961)}]{Kasugai1961}
Kasugai H, Kasegai T (1961) Note on minimax regret ordering policy-static and
  dynamic solutions and a comparison to maximin policy. \emph{J. Oper. Res.
  Soc. Japan} 3(4):155--169.

\bibitem[{Kerrigan \protect\BIBand{} Maciejowski(2004)}]{Kerrigan2004}
Kerrigan E, Maciejowski J (2004) Properties of a new parametrization for the
  control of constrained systems with disturbances. \emph{Proc. 2004 Amer.
  Control Conf.} 5:4669--4674.

\bibitem[{Kuhn et~al.(2011)Kuhn, Wiesmann, \protect\BIBand{}
  Georghiou}]{Kuhn2011}
Kuhn D, Wiesmann W, Georghiou A (2011) Primal and dual linear decision rules in
  stochastic and robust optimization. \emph{Math. Programming} 130(1):177--209.

\bibitem[{Langenhoff \protect\BIBand{} Zijm(1990)}]{Langenhoff1990}
Langenhoff L, Zijm W (1990) An analytical theory of multi-echelon
  production/distribution systems. \emph{Statistica Neerlandica}
  44(3):149--174.

\bibitem[{L\"ofberg(2003)}]{Lofberg2003}
L\"ofberg J (2003) Approximations of closed-loop minimax mpc. \emph{Proc. 42nd
  IEEE Conf. Decision Control} 2:1438--1442.

\bibitem[{Maechler et~al.(2016)Maechler, Rousseeuw, Struyf, Hubert,
  \protect\BIBand{} Hornik}]{package_cluster}
Maechler M, Rousseeuw P, Struyf A, Hubert M, Hornik K (2016) \emph{cluster:
  Cluster Analysis Basics and Extensions}. R package version 2.0.5.

\bibitem[{Mohajerin~Esfahani \protect\BIBand{} Kuhn(2017)}]{esfahani2015data}
Mohajerin~Esfahani P, Kuhn D (2017) Data-driven distributionally robust
  optimization using the wasserstein metric: performance guarantees and
  tractable reformulations. \emph{Math. Programming} ISSN 1436-4646.

\bibitem[{Morton(1978)}]{Morton1978}
Morton TE (1978) The nonstationary infinite horizon inventory problem.
  \emph{Management Sci.} 24(14):1474--1482.

\bibitem[{Muharremoglu \protect\BIBand{} Tsitsiklis(2008)}]{Muharremoglu2008}
Muharremoglu A, Tsitsiklis JN (2008) A single-unit decomposition approach to
  multiechelon inventory systems. \emph{Oper. Res.} 56(5):1089--1103.

\bibitem[{Nohadani \protect\BIBand{} Roy(2017)}]{nohadani2017robust}
Nohadani O, Roy A (2017) Robust optimization with time-dependent uncertainty in
  radiation therapy. \emph{IISE Trans. Healthcare Sys. Eng.} 7(2):81--92.

\bibitem[{Nohadani \protect\BIBand{} Sharma(2018)}]{nohadani2018}
Nohadani O, Sharma K (2018) Optimization under decision-dependent uncertainty.
  \emph{SIAM J. Optim, arXiv preprint arXiv:1611.07992} .

\bibitem[{Rikun(2011)}]{Rikun2011}
Rikun A (2011) \emph{Applications of robust optimization to queueing and
  inventory systems}. Ph.D. thesis, Massachusetts Institute of Technology.

\bibitem[{Rosling(1989)}]{Rosling1989}
Rosling K (1989) Optimal inventory policies for assembly systems under random
  demands. \emph{Oper. Res.} 37(4):565--579.

\bibitem[{Scarf(1958)}]{Scarf1958}
Scarf H (1958) \emph{Studies in the Mathematical Theory of Inventory and
  Production (eds. K.J. Arrow and S. Karlin and H. Scarf)}, chapter A Min-Max
  Solution of An Inventory Problems (Stanford University Press, Stanford, CA).

\bibitem[{Scarf(1960)}]{Scarf1960}
Scarf H (1960) \emph{Mathematical Methods in the Social Sciences (eds. K.J.
  Arrow and S. Karlin and P. Suppes)}, chapter The Optimality of $({s},{S})$
  policies in the dynamic inventory problem (Stanford University Press,
  Stanford, CA).

\bibitem[{Sethi \protect\BIBand{} Cheng(1997)}]{Sethi1997}
Sethi S, Cheng F (1997) Optimality of $({s},{S})$ policies in inventory models
  with markovian demand. \emph{Oper. Res.} 45(6):931--939.

\bibitem[{Uthayakumar \protect\BIBand{}
  Priyan(2013)}]{uthayakumar2013pharmaceutical}
Uthayakumar R, Priyan S (2013) Pharmaceutical supply chain and inventory
  management strategies: optimization for a pharmaceutical company and a
  hospital. \emph{Oper. Res. Health Care} 2(3):52--64.

\bibitem[{Van~Mieghem \protect\BIBand{} Rudi(2002)}]{mieghem2002newsvendor}
Van~Mieghem JA, Rudi N (2002) Newsvendor networks: Inventory management and
  capacity investment with discretionary activities. \emph{Manufacturing \&
  Service Oper. Management} 4(4):313--335.

\bibitem[{Van~Parys et~al.(2016)Van~Parys, Kuhn, Goulart, \protect\BIBand{}
  Morari}]{van2016distributionally}
Van~Parys BP, Kuhn D, Goulart PJ, Morari M (2016) Distributionally robust
  control of constrained stochastic systems. \emph{IEEE Trans. Autom. Control}
  61(2):430--442.

\bibitem[{Webman(2012)}]{pmref2}
Webman E (2012) How much is it worth?
  \url{http://www.ncpa.co/pdf/APNOV12-HowMuchWorth.pdf}.

\bibitem[{Webman(2016)}]{pmref1}
Webman E (2016) Value added.
  \url{http://www.ncpa.co/pdf/webman-pharmacy-valuation-article.pdf}.

\end{thebibliography}
\newpage
\APPENDIX{Appendix}
Proofs of Auxiliary Results.
\section{Proof of Proposition \ref{prop:van_mieghem}.}
First suppose that the nominal mean $\mbs{\mu}$ increases by $\mbs{\Delta \mu} \geq \mb{0}$.
Then there exists $\mbs{\Delta} \mb{s}^* \geq \mb{0}$ and $\mbs{\Delta}\mb{x}^* \geq \mb{0}$ that solves the nominal problem with deterministic demand $\mbs{\Delta \mu}$, with non-negative profit $P \geq 0$.
Thus the worst-case profit with $\mbs{\mu} + \mbs{\Delta \mu}$ increases at least by $P$.
On the other hand, one can show that the set $\mU$ increases in any of $\lambda_1,\ldots,\lambda_l$, followed by that the objective value decreases.
\hfill\Halmos \endproof

\section{Proof of Proposition \ref{prop:multi_feasible}.}
For fixed $\{\mb{w}_t, \mb{W}_{\tau,t}\}$ and $\mb{D}_{[1:T]}$, the inner maximization in problem (\ref{eq:obj_multi_FA}) with respect to ${\mb{X}}_{[1:T]}$ is a linear program in which $\{\mb{w}_t, \mb{W}_{\tau,t}\}$ are on the right-hand side.
It follows that for any ${\mb{D}}_{[1:T]}$, the inner maximization problem (\ref{eq:obj_multi_FA}) is concave in $\{\mb{w}_t, \mb{W}_{\tau,t}\}$.
Hence the objective function in (\ref{eq:obj_multi_FA}) is concave as well, because it is a pointwise infimum of concave functions.
Moreover, the problem is always feasible with assigning zero vectors and matrices to $\{\mb{w}_t, \mb{W}_{\tau,t}\}$.
Finally, applying strong duality to constraint (\ref{eq:constr_multi_FA_first}) shows that a feasible set of $\{\mb{w}_t, \mb{W}_{\tau,t}\}$ is a polyhedron and hence, convex.
\hfill\Halmos \endproof

\section{Proof of Proposition \ref{lemma:initial_single}.} 
$\Phi(s_0,d_0) := \underset{\pi}{\text{max}}~ V \Big( \pi(w_t, W_{\tau,t});s_0,d_0 \Big)$, where $V\Big(\pi(w_t, W_{\tau,t});s_0,d_0\Big)$ is defined as
\begin{align*}
V \Big(\pi(w_t, W_{\tau,t});s_0,d_0 \Big)
&:= \min_{{\mb{d}} \in \mU^T} \max_{{\mb{x}},{\mb{s}}}~ P\Big(\pi(w_t,W_{\tau,t}),{\mb{d}},{\mb{x}};s_0,d_0\Big) 
\\
\begin{split}
~= \min_{{\mb{d}} \in \mU^T}~ \max_{{\mb{x}},{\mb{s}}}~ \Bigg[
-c_S \Bigg( \sum_{t=1}^T s_t \Bigg) -c_H \sum_{t=1}^T \Bigg( s_0 + \sum_{\tau=1}^t ( s_\tau - x_\tau) \Bigg) 
~~~~~~~~~\\
-c_P \sum_{t=1}^T \Bigg( d_0 + \sum_{\tau=1}^t ( d_\tau - x_\tau) \Bigg) + r \Bigg( \sum_{t=1}^T x_t \Bigg) \Bigg],
\end{split}
\end{align*}
where the inner maximization problem has a feasible set $\mathcal{X}(\pi, {\mb{d}}, s_0, d_0)$.

We first claim that the following equation holds for any affine policy $\pi(w_t, W_{\tau,t})$ such that $s_0 \leq w_1$,
\begin{equation}
V\Big(\lbar{\pi}(\lbar{w}_t, \lbar{W}_{\tau,t}); s_0, d_0 \Big) = V \Big( \pi(w_t, W_{\tau,t}); 0, 0 \Big) + c_S s_0 + (r-c_S) d_0,
\label{eq:cost_reform}
\end{equation}
where $\lbar{\pi}(\lbar{w}_t, \lbar{W}_{\tau,t})$ is defined as (\ref{eq:initial_single}).
For any demand realization ${\mb{d}} \in \mU^T$, let $\agg{s}_t = \agg{s}_t(\pi, {\mb{d}})$ and $\agg{x}_t = \agg{x}_t(\pi, {\mb{d}})$ be an aggregated order quantity and the corresponding optimal processing activity at time $t$ with zero initial input and demand, i.e., it solves the inner maximization problem of $V(\pi;0,0)$ with ${\mb{d}}$.
Let $\agg{\lbar{s}}_t = \agg{\lbar{s}}_t(\lbar{\pi},{\mb{d}})$ be an order quantity by the policy $\lbar{\pi}$ for $V(\lbar{\pi};s_0,d_0)$, and define $\agg{\lbar{x}}_t = \agg{\lbar{x}}_t({\mb{d}}) = \agg{x}_t({\mb{d}}) + d_0$.
Then by (\ref{eq:initial_single}), $\agg{\lbar{s}}_t = \agg{s}_t - s_0 + d_0$ for every $t=1,\ldots,T$, and thus
\begin{align*}
P\Big(&\pi(w_t,W_{\tau,t}), {\mb{d}}, {\mb{x}};0,0\Big) + c_S s_0 + (r-c_S) d_0 \\
& = \Big( -c_S \agg{s}_T - c_H \sum_{t=1}^{T} (\agg{s}_t - \agg{x}_t) - c_P \sum_{t=1}^{T} (\agg{d}_t - \agg{x}_t) + r \agg{x}_T \Big) + (r-c_S) d_0 + c_S s_0 \\
& = -c_S (\agg{s}_T-s_0 +d_0) - c_H \sum_{t=1}^{T} (\agg{s}_t + d_0) - c_P \sum_{t=1}^{T} (d_0+\agg{d}_t) + (c_H + c_P) \sum_{t=1}^{T} (\agg{x}_t + d_0) + r (\agg{x}_T + d_0) \\
%& = -c_S \lbar{S}_T - c_H \sum_{t=1}^{T} (s_0 + \lbar{S}_t) - c_P \sum_{t=1}^{T} (d_0 + D_t) + (c_H + c_P) \sum_{t=1}^{T} \lbar{X}_t + r \lbar{X}_T \\
& = -c_S \agg{\lbar{s}}_T - c_H \sum_{t=1}^{T} (s_0 + \agg{\lbar{s}}_t - \agg{\lbar{x}}_t) - c_P \sum_{t=1}^{T} (d_0 + \agg{d}_t - \agg{\lbar{x}}_t) + r \agg{\lbar{x}}_T \\
& = P\Big(\lbar{\pi}(\lbar{w}_t,\lbar{W}_{\tau,t}), {\mb{d}}, {\lbar{\mb{x}}}; s_0, d_0 \Big).
\end{align*}
Since ${\lbar{\mb{x}}}$ is a feasible solution of the inner maximization problem in $V( \lbar{\pi}; s_0, d_0)$, the LHS of (\ref{eq:cost_reform}) is greater than or equal to the RHS.
Similar argument can be made to show the contrary and this concludes the proof of (\ref{eq:cost_reform}).

Now suppose $\lbar{\pi}^*(\lbar{w}_t^*,\lbar{W}_{\tau,t}^*)$ is not optimal, and let $\lbar{\varphi}^*(\lbar{v}_t^*,\lbar{V}_{\tau,t}^*)$ be an optimal solution of $\Phi(s_0,d_0)$.
Now one can easily check that by (\ref{eq:cost_reform}),
\begin{align*}
\Phi(s_0, d_0) = V(\lbar{\varphi}^*; s_0,d_0) 
&= V(\varphi^*; 0, 0) + c_S s_0 + (r-c_S) d_0 \\
&\leq V(\pi^*; 0, 0) + c_S s_0 + (r-c_S) d_0
= V(\lbar{\pi}^*;s_0,d_0) < \Phi(s_0,d_0),
\end{align*}
which makes a contradiction.
Finally, the proof of Eq. (\ref{eq:cost_reform}) directly shows that both $\Phi(0,0)$ and $\Phi(s_0,d_0)$ shares a common worst-case scenario among $\mU^T$.
\hfill\Halmos \endproof

\section{Proof of Proposition \ref{prop:PA_worst_scenario}.}
Let $\lbar{\pi}_\text{PA} = (\lbar{\pi}_1, \ldots, \lbar{\pi}_N)$, where each $\lbar{\pi}_j$ solves the $j^\text{th}$ subproblem.
Then the worst-case scenario ${(\mb{d}^*_1,\ldots,\mb{d}^*_N)}$ solves an optimization problem
\begin{equation}
\label{eq:prop3_proof}
	\min_{{\mb{d}_1,\ldots,\mb{d}_N}} ~\max_{{\mb{x}_1,\ldots,\mb{x}_N}}~ \bigg[
	\sum_{j=1}^N P_j(\lbar{\pi}_j,{\mb{d}_j},{\mb{x}_j};u_{t_{j-1}}^s,u_{t_{j-1}}^d) \bigg].
\end{equation}
Since $\lbar{\pi}_\text{PA}$ is a well-defined periodic-affine policy, (\ref{eq:prop3_proof}) is rewritten as
\begin{align*}
&\min_{{\mb{d}_1,\ldots,\mb{d}_N}} ~\max_{{\mb{x}_1,\ldots,\mb{x}_N}}~ \bigg[
\sum_{j=1}^N P_j(\lbar{\pi}_j,{\mb{d}_j},{\mb{x}_j};u_{t_{j-1}}^s,u_{t_{j-1}}^d) \bigg] \\
&= \min_{{\mb{d}_1,\ldots,\mb{d}_N}} ~\max_{{\mb{x}_1,\ldots,\mb{x}_N}}~ \bigg[
\sum_{j=1}^N \widetilde{P}_j(\lbar{\pi}_j,{\mb{d}_j},{\mb{x}_j};u_{t_{j-1}}^s,u_{t_{j-1}}^d) \bigg] \\
&= \min_{{\mb{d}_1,\ldots,\mb{d}_N}} ~\max_{{\mb{x}_1,\ldots,\mb{x}_N}}~ \bigg[
\sum_{j=1}^N P_j^\text{PA} (\lbar{\pi}_j,{\mb{d}_j},{\mb{x}_j}) \bigg] \\
&=~ \sum_{j=1}^N  \bigg[ ~\min_{{\mb{d}_j} \in \mU^j} \max_{{\mb{x}_j}}~
 P_j^\text{PA} (\lbar{\pi}_j,{\mb{d}_j},{\mb{x}_j}) \bigg] \\
 &=~ \sum_{j=1}^N {\max_{{\mb{x}_j}}}~ P_j^\text{PA} (\lbar{\pi}_j,{\mb{d}^*_j},{\mb{x}_j}) ,
\end{align*}
where $\widetilde{P}_j$ and $P_j^\text{PA}$ are defined in (\ref{eq:reformulate_final}) and (\ref{eq:obj_PA_new}).
As a result, the overall objective function is separable for each subproblem, and hence, the worst-case scenario consists of those of the subperiods.
\hfill\Halmos \endproof

\section{Proof of Theorem \ref{thm:optimality_PA}.}
Since every affine policy is feasible to the DP problem (\ref{eq:obj_PA_DP}), proving $V_\text{{Aff}}^* \leq V_\text{DP}^*$ directly follows.
Thus it suffices to show $V_\text{PA}^* = V_\text{DP}^*$, and let $\lbar{\pi}_\text{PA} = (\lbar{\pi}_1, \ldots, \lbar{\pi}_N)$ be an output of the periodic-affine algorithm.
That is, $\lbar{\pi}_j$ is an affine-IBS policy associated with $\pi_j = \pi_j(w_t^{(j)},W_{\tau,t}^{(j)})$, where $\pi_j$ solves (\ref{prob:PA_algorithm_main}) for each $j=1,\ldots,N$.
Define $V_j(u^s, u^d)$ as a worst-case optimal profit from the $j^{\text{th}}$ subperiod to the last period, where $u^s$ and $u^d$ are current on-hand input and backlogged demand.
Our framework justifies using the (robust) optimality equation and a (worst-case) value function approach in robust dynamic programming scheme; we refer \cite{iyengar2005robust} to readers for technical details.
We will show for every $j=1,\ldots,N$,
\begin{itemize}
	\item[(a)] $V_j(u^s, u^d)$ is concave in $(u^s,u^d)$, and 
	\item[(b)] $V_j(u^s,u^d) = V_j(0,0) + c_S u^s + (r-c_S) u^d$ for every $0 \leq u^s \leq w_1^{(j)}$, $u^d \geq 0$.
\end{itemize}
by mathematical induction.

Now consider $j=N$ and suppose $u^s \leq w_1^{(N)}$. 
Then $V_N(u^s,u^d)$ can be written as 
\begin{equation*}
\begin{aligned}
V_N(u^s,u^d) := \!\!\! &\max_{\pi \in \Pi_\text{aff}(\mU^{N})}  \min_{{\mb{d}_N}} \max_{{\mb{x}_N},{\mb{s}_N}}~ P_N (\pi, {\mb{d}_N}, {\mb{x}_N}; u^s, u^d) \\
&~~~~\text{s.t. }~~~~ ({\mb{x}_N}, {\mb{s}_N}) \in \mathcal{X}(\pi, {\mb{d}_N}, u^s, u^d),
\end{aligned}
\end{equation*}
where $\pi = \pi(w_t, W_{\tau,t})$ and the constraints in $\mathcal{X}(\pi, {\mb{d}_N}, u^s, u^d)$ can be rearranged so that the right hand sides are linear in $(u^s,u^d,w_t,W_{\tau,t})$.
Since $P_N$ is concave in $({\mb{x}_N},w_t,W_{\tau,t},u^s,u^d)$ and $\mathcal{X}(\pi, {\mb{d}_N}, u^s, u^d)$ defines a polyhedron for any $\pi$, $u^s$, and $u^d$, the objective function within the min operator is concave in $(u^s, u^d, w_t, W_{\tau,t})$ by concavity preservation under maximization.
Since a pointwise infimum of concave functions are concave and applying concavity preservation under maximization again to the outermost max operator, we finally have that $V_N(u^s,u^d)$ is concave in $(u^s,u^d)$.
On the other hand, (b) follows directly from Proposition \ref{lemma:initial_single} for $j=N$.

Now suppose that both (a) and (b) hold for any $1 < j \leq N$, and let $u^s \leq w_1^{(j-1)}$ and $u^d \geq 0$.
Then from the optimality equation, we have
\begin{equation}
V_{j-1}(u^s,u^d)
=~ \!\!\!\max_{\pi \in \Pi_\text{aff}(\mU^{j-1})} ~ \min_{{\mb{d}_{j-1}}} ~ \max_{{\mb{x}_{j-1}}}~ \bigg[ P_{j-1} \Big( \pi, {\mb{d}_{j-1}}, {\mb{x}_{j-1}}; u^s, u^d \Big) + V_j(u_j^s, u_j^d) \bigg].
\label{eq:thm_bellman}
\end{equation}
Since $V_j(u_j^s, u_j^d)$ is concave in $(u_j^s, u_j^d)$ and $(u_j^s,u_j^d)$ can be expressed as affine functions of $(w_t, W_{\tau,t}, u^s, u^d)$, applying the above argument shows that $V_{j-1}(u^s,u^d)$ is also concave in $(u^s, u^d)$.
In addition, we have 
\begin{align*}
	V_{j-1}(u^s,u^d)
	&=\!\!\! \max_{\pi \in \Pi_\text{aff}(\mU^{j-1})} ~ \min_{{\mb{d}_{j-1}}} ~ \max_{{\mb{x}_{j-1}}}~ \bigg[ P_{j-1} \Big( \pi, {\mb{d}_{j-1}}, {\mb{x}_{j-1}}; u^s, u^d \Big) + V_j(u_j^s, u_j^d) \bigg] \\
	&\leq\!\!\! \max_{\pi \in \Pi_\text{aff}(\mU^{j-1})} ~ \min_{{\mb{d}_{j-1}}} ~ \max_{{\mb{x}_{j-1}}}~ \bigg[ P_{j-1} \Big( \pi,{\mb{d}_{j-1}}, {\mb{x}_{j-1}}; u^s, u^d \Big) + V_j(0,0) + c_S u_j^s + (r-c_S) u_j^d \bigg] \\
	&= \!\!\!\max_{\pi \in \Pi_\text{aff}(\mU^{j-1})} ~ \min_{{\mb{d}_{j-1}}} ~ \max_{{\mb{x}_{j-1}}}~ \bigg[ P_{j-1} \Big( \pi, {\mb{d}_{j-1}}, {\mb{x}_{j-1}}; u^s, u^d \Big) + c_S u_j^s + (r-c_S) u_j^d \bigg] + V_j(0,0) \\
	&= \!\!\!\max_{\pi \in \Pi_\text{aff}(\mU^{j-1})} ~ \min_{{\mb{d}_{j-1}}} ~ \max_{{\mb{x}_{j-1}}}~ \bigg[ \widetilde{P}_{j-1} \Big( \pi,{\mb{d}_{j-1}}, {\mb{x}_{j-1}}; u^s, u^d \Big) + c_S u^s + (r-c_S) u^d \bigg] + V_j(0,0) \\
	&= c_S u^s + (r-c_S) u^d + \max_{\pi \in \Pi_\text{aff}(\mU^{j-1})} ~ \min_{{\mb{d}_{j-1}}} ~ \max_{{\mb{x}_{j-1}}}~  \bigg[ \widetilde{P}_{j-1} \Big( \pi, {\mb{d}_{j-1}}, {\mb{x}_{j-1}}; 0, 0 \Big) \bigg] + V_j(0,0) \\
	&\leq V_{j-1}(u^s,u^d).
\end{align*}
The first inequality comes from that both (a) and (b) hold for $V_j$, and the third equality is from the definition of $\widetilde{P}_{j-1}$.
Finally, the last inequality is a worst-case profit from $j^\text{th}$ subperiod with a policy $\lbar{\pi}_\text{PA}$, by Assumption \ref{assumption:optimality_PA}.
Since $\lbar{\pi}_\text{PA}$ is a feasible policy to the DP, the last inequality follows.
This shows that whenever $u^s \leq w_1^{(j-1)}$, then the value function $V_{j-1}(u^s,u^d)$ is achieved with a policy $(\lbar{\pi}_{j-1},\ldots, \lbar{\pi}_N)$.
Finally we have
\begin{align*}
	V_{j-1}(u^s,u^d)
	&= c_S u^s + (r-c_S) u^d + \max_{\pi \in \Pi_\text{aff}(\mU^{j-1})} ~ \min_{{\mb{d}_{j-1}}} ~ \max_{{\mb{x}_{j-1}}}~  \bigg[ \widetilde{P}_{j-1} \Big( \pi, {\mb{d}_{j-1}}, {\mb{x}_{j-1}}; 0,0 \Big) \bigg] + V_j(0,0) \\
	&= V_{j-1}(0,0) + c_S u^s + (r-c_S)u^d,
\end{align*}
and thus (b) holds for $j-1$.
By Assumption \ref{assumption:optimality_PA}, $\lbar{\pi}_\text{PA} = (\lbar{\pi}_1,\ldots, \lbar{\pi}_N)$ satisfies $u_{t_j}^s \leq w_1^{(j+1)}$ for every $j=1,\ldots, N-1$ and every realization of demands, and hence, $\lbar{\pi}_\text{PA}$ is Bellman-optimal to the DP (\ref{eq:obj_PA_DP}) and this concludes with $V_\text{PA}^* = V_\text{DP}^*$.
\hfill\Halmos \endproof

\section{Proof of Proposition \ref{proposition:PA_vs_FA}.}
It suffices to show that an optimal periodic-affine policy is indeed an affine policy.
Using the same notations in Theorem \ref{thm:optimality_PA} and without loss of generality, we may assume that $N=2$ and let $\lbar{\pi}_\text{PA}^* = (\lbar{\pi}_1, \lbar{\pi}_2)$, where $\pi_j = (w_t^{(j)}, W_{\tau,t}^{(j)})$ for $j=1,2$.
By definition of periodic-affine policies, we only need to check that if an order quantity at time $t_1+1$ is affine in $\mU = \mU^1 \times \mU^2$.
Recall that $\lbar{\pi}_\text{PA}^*$ determines order quantity at $t_1+1$ as
\begin{align*}
w_{1}^{(2)} - u_{t_1}^s + u_{t_1}^d 
&= w_{1}^{(2)} - \text{max } \bigg( \sum_{t=1}^{t_1} (s_t - d_t) ,~0 \bigg)
+ \text{max } \bigg( \sum_{t=1}^{t_1} (d_t - s_t) ,~0 \bigg) \\
& = w_{1}^{(2)} + \bigg( \sum_{t=1}^{t_1} (d_t - s_t) \bigg).
\end{align*}
It is affine in ${\mb{d}_1}$, since $s_t$ is affine in ${\mb{d}_1}$, and this concludes the proof.
\hfill\Halmos \endproof

\section{Proof of Theorem \ref{thm:suboptimal_PA}.}
We use the value function $V_j(u^s,u^d)$ defined in Theorem \ref{thm:optimality_PA}.
From concavity $V_j(u^s,u^d)$ and using (b), we have 
\begin{equation*}
V_j(u^s, u^d) \leq V_j(0,0) + c_S u^s + (r-c_S) u^d
\end{equation*}
for every $u^s \geq 0$ and $u^d \geq 0$.
We will show that
\begin{equation}
V_j(u^s, u^d) \leq c_S u^s + (r-c_S) u^d + \sum_{k=j}^N \widetilde{f}_k^*~~~~\forall j=1,\ldots, N
\label{eq:subopt_induction}
\end{equation}
by induction, and plugging $j=1$ and $u^s=u^d=0$ into \eqref{eq:subopt_induction} concludes the proof.

From $V_N(0,0) = \widetilde{f}_N^*$, we have that \eqref{eq:subopt_induction} holds for $j=N$.
Now suppose $1 < j \leq N$.
Then from the optimality equation we have
\begin{equation*}
\begin{aligned}
V_{j-1}(u^s,u^d)
&= \!\!\!\max_{\pi \in \Pi_\text{aff}(\mU^{j-1})} ~ \min_{{\mb{d}_{j-1}}} ~ \max_{{\mb{x}_{j-1}}}~ \bigg[ P_{j-1} \Big( \pi, {\mb{d}_{j-1}}, {\mb{x}_{j-1}}; u^s, u^d \Big) + V_j(u_j^s, u_j^d) \bigg] \\
&\leq\!\!\! \max_{\pi \in \Pi_\text{aff}(\mU^{j-1})} ~ \min_{{\mb{d}_{j-1}}} ~ \max_{{\mb{x}_{j-1}}}~ \bigg[ P_{j-1} \Big( \pi, {\mb{d}_{j-1}}, {\mb{x}_{j-1}}; u^s, u^d \Big) + c_S u_j^s + (r-c_S) u_j^d + \sum_{k=j}^N \widetilde{f}_k^* \bigg] \\
%&=~ \underset{\pi \in \Pi_\text{aff}(\mU^{j-1})}{\text{max}} ~ \underset{d_{I_{j-1}}}{\text{min}} ~ \underset{x_{I_{j-1}} \in \mathcal{P}(\mU^{j-1})}{\text{max}}~ \bigg[ P_{j-1} \Big( \pi, d_{I_{j-1}}, x_{I_{j-1}}; u^s, u^d \Big) + c_S u_j^s + (r-c_S) u_j^d \bigg] + \sum_{k=j}^N \widetilde{f}_k^* \\
%&=~ \underset{\pi \in \Pi_\text{aff}(\mU^{j-1})}{\text{max}} ~ \underset{d_{I_{j-1}}}{\text{min}} ~ \underset{x_{I_{j-1}} \in \mathcal{P}(\mU^{j-1})}{\text{max}}~ \bigg[ \widetilde{P}_{j-1} \Big( \pi, d_{I_{j-1}}, x_{I_{j-1}}; u^s, u^d \Big) \bigg] + \sum_{k=j}^N \widetilde{f}_k^* \\
&\leq~ c_S u^s + (r-c_S) u^d +\!\!\! \max_{\pi \in \Pi_\text{aff}(\mU^{j-1})} ~ \min_{{\mb{d}_{j-1}}} ~ \max_{{\mb{x}_{j-1}}}~ \bigg[ \widetilde{P}_{j-1} \Big( \pi, {\mb{d}_{j-1}}, {\mb{x}_{j-1}}; 0, 0 \Big) \bigg] + \sum_{k=j}^N \widetilde{f}_k^* \\
&=~ c_S u^s + (r-c_S) u^d + \sum_{k=j-1}^N \widetilde{f}_k^*,
\end{aligned}
\end{equation*}
where the first inequality holds from the induction hypothesis and the third equality is from definition of $\widetilde{P}_{j-1}$.
One can show that the maximization problem in the third equality is concave in $u^s$ and $u^d$, as similar in the proof of Theorem \ref{thm:optimality_PA} and this concludes the proof.
\hfill\Halmos \endproof

\section{Proof of Theorem \ref{corollary:optimality_PA_multi}.}
All the proofs of Theorem \ref{thm:optimality_PA} and \ref{thm:suboptimal_PA} can be extended into multi-station networks, by using a basis matrix $\mb{R}_B$ to replace $c_S u^s$ and $(r-c_S) u^d$ terms in the proof with $\mb{c}^\top_S \mb{u}^s$ and $(\mb{R}^\top_B \mb{r} - \mb{R}^\top_B \mb{R}^\top_S \mb{c}_S)^\top \mb{u}^d$, respectively.
This expressions are still linear in $\mb{u}^s$ and $\mb{u}^d$, hence all the arguments in the proof are valid.
\hfill\Halmos\endproof

\section{Proof of Theorem \ref{thm:optimality_PA_infty}.}
It suffices to show for single-station cases, since it is straightforward to extend the result to general multi-station networks, as in Theorem \ref{corollary:optimality_PA_multi}.
Note that the optimality equation for the infinite horizon problem is written as
\begin{equation}
V_\infty(u^s, u^d) = \max_{\pi \in \Pi_\text{aff}(\mU^{{T_0}})} ~ \min_{{D}_{[1:{{T_0}}]}} ~ \max_{{X}_{[1:{{T_0}}]}}~ \bigg[ P \Big( \pi, {D}_{[1:{{T_0}}]}, {X}_{[1:{{T_0}}]}; u^s, u^d \Big) + \beta^{{T_0}} V_\infty(\lbar{u}^s, \lbar{u}_j^d) \bigg],
\label{proof:optimality_infty}
\end{equation}
where $\lbar{u}^s$ and $\lbar{u}^d$ denotes on-hand input and backorders after ${{T_0}}$ periods (one stage).

We impose mild conditions so that the optimality equation (\ref{proof:optimality_infty}) defines a contraction mapping and there exists $V_\infty$ which is the unique fixed point. (See \cite{iyengar2005robust} for details.)
Hence the value iteration algorithm is well-defined, and let $V_n(u^s,u^d)$ be a value function after $n$ iterations.
Recalling that $\lbar{u}^s$ and $\lbar{u}^d$ are expressed as linear functions of $u^s$ and $u^d$, one can show by applying concavity preservation under maximization as similar in Theorem \ref{thm:optimality_PA} that if $V_n(u^s,u^d)$ is concave, then so $V_{n+1}(u^s,u^d)$ is.
Since we can start with any bounded continuous function for the value iteration algorithm, we conclude that $V_\infty(u^s,u^d)$ is concave in $(u^s,u^d)$.

In this setting, there exists a stationary optimal policy $\pi_\infty = (\pi,\pi, \ldots)$ where $\pi = \pi(u^s,u^d)$ is defined for each subperiod of length ${{T_0}}$.
By Proposition \ref{lemma:initial_single}, we can see that $V_\infty(u^s,u^d) = V_\infty(0,0) + c_S u^s + (r-c_S) u^d$ for $u^s \leq w_1$ and with concavity of $V_\infty$, we have
\begin{equation*}
V_\infty(u^s,u^d) \leq V_\infty(0,0) + c_S u^s + (r-c_S) u^d
\end{equation*}
for every $u^s \geq 0$ and $u^d \geq 0$.

Let $V_\infty(\pi_\infty)$ be a worst-case objective value under policy $\pi_\infty$, and $V_\infty^*$ be an optimal value of the DP problem.
Since $\pi_\infty$ is feasible to the DP by Assumption \ref{assumption:optimality_PA}, we have $V_\infty(\pi_\infty) \leq V_\infty^*$.
On the other hand,
\begin{align*}
V_{\infty}^*
&=~ \max_{\pi \in \Pi_\text{aff}(\mU^{{T_0}})} ~ \min_{{D}_{[1:{{T_0}}]}} ~ \max_{{X}_{[1:{{T_0}}]}}~ \bigg[ P \Big( \pi, {D}_{[1:{{T_0}}]}, {X}_{[1:{{T_0}}]}; 0, 0 \Big) + \beta^{{T_0}} V_\infty(\lbar{u}^s, \lbar{u}^d) \bigg] \\
&\leq~ \max_{\pi \in \Pi_\text{aff}(\mU^{{T_0}})} ~ \min_{{D}_{[1:{{T_0}}]}} ~ \max_{{X}_{[1:{{T_0}}]}}~ \bigg[ P \Big( \pi, {D}_{[1:{{T_0}}]}, {X}_{[1:{{T_0}}]}; u^s, u^d \Big) + \beta^{{T_0}} \left( c_S \lbar{u}^s + (r-c_S) \lbar{u}^d +  V_\infty(0,0) \right) \bigg] \\
&=~ \max_{\pi \in \Pi_\text{aff}(\mU^{{T_0}})} ~ \min_{{D}_{[1:{{T_0}}]}} ~ \max_{{X}_{[1:{{T_0}}]}}~ P_\infty^\text{PA} \Big( \pi, {D}_{[1:{{T_0}}]}, {X}_{[1:{{T_0}}]} \Big) + \beta^{{T_0}} V_\infty(0,0) \\
&= V_\infty(\pi_{\infty}),
\end{align*}
by Assumption \ref{assumption:optimality_PA} (this step is similar to Theorem \ref{thm:optimality_PA}), and this implies that an optimal value to the DP is achieved by $\pi_\infty$, where the stationary policy $\pi$ is obtained by solving the optimization problem \eqref{eq:obj_infty_PA}.
\hfill\Halmos \endproof

\end{document}